\documentclass[UTF-8,reqno]{amsart}
\usepackage{enumerate, bbm, framed}
\usepackage{pgfplots}
\usetikzlibrary{plotmarks, pgfplots.groupplots}
\usepackage{amssymb,url,color, booktabs}
\usepackage{tikz}
\usepackage{pgfplots}
\pgfplotsset{compat=1.17}
\usetikzlibrary{intersections}
\usepgfplotslibrary{fillbetween}
\usepackage{algorithm}       
\usepackage{algpseudocode}   
\usepackage{caption}  
\usepackage{enumitem}  

\usepackage{caption}
\usepackage{float}

\usepackage{mathrsfs}
\usepackage{dutchcal}
\usepackage{threeparttable}

\usepackage{colortbl}
\usepackage{color}
\usepackage[colorlinks=true]{hyperref}
\hypersetup{
    linkcolor=blue,          
    citecolor=red,        
    filecolor=blue,      
    urlcolor=cyan
}

\usepackage{graphicx}
\usepackage{subcaption}
\usepackage{caption}
\usepackage{lipsum}

\usepackage{pgfplots}
\pgfplotsset{compat=newest}
\usepackage{color}
\usepackage{ulem}

 \usepackage{scalerel} 

\definecolor{MyDarkBlue}{cmyk}{0.8,0.3,0.8,0.4}
\definecolor{yellow}{rgb}{0.99,0.99,0.70}
\definecolor{white}{rgb}{1.0,1.0,1.0}
\definecolor{black}{rgb}{0.00,0.00,0.00}

\newcommand{\bsigma}{\boldsymbol{\sigma}}

\numberwithin{equation}{section}

\newcommand{\be}{\begin{eqnarray}}
\newcommand{\ee}{\end{eqnarray}}
\newcommand{\ce}{\begin{eqnarray*}}
\newcommand{\de}{\end{eqnarray*}}
\newtheorem{theorem}{Theorem}[section]
\newtheorem{lemma}[theorem]{Lemma}
\newtheorem{remark}[theorem]{Remark}
\newtheorem{definition}[theorem]{Definition}
\newtheorem{proposition}[theorem]{Proposition}
\newtheorem{Examples}[theorem]{Example}
\newtheorem{corollary}[theorem]{Corollary}

\usepackage[nobysame]{amsrefs}
\BibSpec{article}{%
+{}{\PrintAuthors} {author}
+{,}{ \textrm} {title}
+{.}{ \textit} {journal}
+{,}{ \textbf} {volume}
+{}{ \parenthesize} {date}
+{,}{ } {pages}
+{.}{ arXiv:} {eprint}
+{.}{} {transition}
}
\BibSpec{book}{%
+{}{\PrintAuthors} {author}
+{,}{ \textit} {title}
+{.}{ \textrm} {series} 
+{,}{ Vol.} {volume} 
+{.}{ } {publisher}
+{,}{ } {date}
+{.}{} {transition}
}

\def\var{{\mathrm{var}}}

\def\eps{\varepsilon}
\def\e{\mathrm{e}}

\def\<{{\langle}}
\def\>{{\rangle}}
\def\({{\Big(}}
\def\){{\Big)}}

\def\bx{{\mathbf{x}}}
\def\tr{\mathrm {tr}}

\def\dif{{\mathord{{\rm d}}}}

\def\no{\nonumber}
\def\={&\!\!=\!\!&}

\def\cE{{\mathcal E}}

\def\cG{{\mathcal G}}
\def\cH{{\mathcal H}}

\def\cN{{\mathcal N}}

\def\cS{{\mathcal S}}

\def\cU{{\mathcal U}}

\def\mE{{\mathbb E}}

\def\mI{{\mathbb I}}

\def\mN{{\mathbb N}}

\def\mP{{\mathbb P}}

\def\mR{{\mathbb R}}

\def\1{{\mathbf{1}}}

\def\sI{{\mathscr I}}

\def\sN{{\mathscr N}}

\def\geq{\geqslant}
\def\leq{\leqslant}

\def\var{{\mathrm{var}}}

\def\eps{\varepsilon}

\def\e{\mathrm{e}}

\def\<{{\langle}}
\def\>{{\rangle}}
\def\({{\Big(}}
\def\){{\Big)}}

\def\bx{{\mathbf{x}}}
\def\tr{\mathrm {tr}}

\def\dif{{\mathord{{\rm d}}}}

\def\no{\nonumber}
\def\={&\!\!=\!\!&}
\def\bt{\begin{theorem}}
\def\et{\end{theorem}}
\def\bl{\begin{lemma}}
\def\el{\end{lemma}}
\def\br{\begin{remark}}
\def\er{\end{remark}}
\def\bx{\begin{Examples}}
\def\ex{\end{Examples}}
\def\bd{\begin{definition}}
\def\ed{\end{definition}}
\def\bp{\begin{proposition}}
\def\ep{\end{proposition}}
\def\bc{\begin{corollary}}
\def\ec{\end{corollary}}

\def\geq{\geqslant}
\def\leq{\leqslant}

\def\<{\langle} \def\>{\rangle}

\def\bpf{\begin{proof}}
\def\epf{\end{proof}}

\def\balpha{\boldsymbol{\alpha}}
\def\bxi{{\boldsymbol{\xi}}}

\allowdisplaybreaks

\begin{document}

\title[Zeroth-Order Optimization via Adaptive Sampling]
{Sampling-Based Zero-Order Optimization Algorithms}
\author{Xicheng Zhang}

\address{Xicheng Zhang:
School of Mathematics and Statistics, Beijing Institute of Technology, Beijing 100081, China\\
Faculty of Computational Mathematics and Cybernetics, Shenzhen MSU-BIT University, 518172 Shenzhen, China\\
		Email: xczhang.math@bit.edu.cn
 }


\thanks{
This work is supported by National Key R\&D program of China (No. 2023YFA1010103) and NNSFC grant of China (No. 12131019)  and the DFG through the CRC 1283 ``Taming uncertainty and profiting from randomness and low regularity in analysis, stochastics and their applications''. }

\begin{abstract}
We propose a novel zeroth-order optimization algorithm based on an efficient sampling strategy. Under mild global regularity conditions on the objective function, we establish non-asymptotic convergence rates for the proposed method. Comprehensive numerical experiments demonstrate the algorithm's effectiveness, highlighting three key attributes: 
(i) Scalability: consistent performance in high-dimensional settings (exceeding 100 dimensions); 
(ii) Versatility: robust convergence across a diverse suite of benchmark functions, including Schwefel, Rosenbrock, Ackley, Griewank, Lévy, Rastrigin, and Weierstrass; and 
(iii) Robustness to discontinuities: reliable performance on non-smooth and discontinuous landscapes. 
These results illustrate the method's strong potential for black-box optimization in complex, real-world scenarios.
\end{abstract}

\maketitle
\setcounter{tocdepth}{2}


\section{Introduction}

Optimization constitutes a fundamental mathematical and computational discipline that arises across diverse fields including engineering, economics, machine learning, and operations research. Traditional optimization methodologies primarily focus on identifying the minimum or maximum of a given objective function, often subject to constraints. In its broadest sense, optimization refers to the process of systematically adjusting inputs to a system or model to achieve optimal performance as measured by a predefined criterion. This criterion is typically expressed through an \textit{objective function}, which may be either maximized (e.g., profit, efficiency) or minimized (e.g., cost, error, time) \cite{nocedal2006numerical}.

Formally, a standard optimization problem in $\mathbb{R}^d$ can be expressed as:
\begin{align}
x^* = \arg\max_{x \in \mathbb{R}^d} f(x) \quad \text{or} \quad x_* = \arg\min_{x \in \mathbb{R}^d} f(x)
\end{align}
subject to
\begin{align} \label{Cons}
g_i(x) \leq 0, \quad i = 1, \dots, m \quad \text{or} \quad h_j(x) = 0, \quad j = 1, \dots, p,
\end{align}
where $f(x)$ denotes the objective function, $g_i(x)$ represent inequality constraints, and $h_j(x)$ equality constraints.

Optimization problems can be classified according to various criteria:
\begin{itemize}
    \item \textbf{Linear vs. Nonlinear:} If both the objective function and all constraints are linear, the problem constitutes a linear program (LP); otherwise, it is a nonlinear program (NLP) \cite{boyd2004convex}.
    \item \textbf{Convex vs. Non-convex:} A problem is convex if the objective function is convex and the feasible region forms a convex set. Convex problems generally admit efficient solutions with guaranteed global optimality \cite{boyd2004convex}.
    \item \textbf{Unconstrained vs. Constrained:} Unconstrained optimization imposes no restrictions on variables, while constrained optimization requires solutions to satisfy additional equality or inequality constraints as specified in \eqref{Cons}.
    \item \textbf{Deterministic vs. Stochastic:} Deterministic optimization assumes perfectly known parameters, whereas stochastic optimization incorporates randomness or uncertainty in model data \cite{shapiro2009lectures}.
    \item \textbf{Discrete vs. Continuous:} Discrete optimization restricts variables to discrete values (e.g., integers), while continuous optimization permits variables to assume any real values within specified ranges.
\end{itemize}

The selection of an appropriate optimization methodology is primarily determined by the underlying problem structure:
\begin{itemize}
    \item \textbf{Gradient-Based Methods:} These techniques utilize derivative information to guide the optimization process. Notable examples include:
    \begin{itemize}
        \item Gradient descent: Iteratively moves in the direction of steepest descent
        \begin{equation}
            x_{k+1} = x_k - \alpha \nabla f(x_k),
        \end{equation}
        where $\alpha$ denotes the learning rate.
        
        \item Newton's method: Incorporates second-order derivative information for accelerated convergence
        \begin{equation}
            x_{k+1} = x_k - [\nabla^2 f(x_k)]^{-1} \nabla f(x_k).
        \end{equation}
        
        \item Stochastic gradient descent (SGD): Introduces controlled noise for improved generalization
        \begin{equation}
            x_{k+1} = x_k - \alpha \nabla f(x_k) + \sqrt{\alpha \varepsilon}\xi_k,
        \end{equation}
        where $\xi_k \sim \mathcal{N}(0,I_d)$ are independent Gaussian random variables.
    \end{itemize}
    Comprehensive treatments can be found in \cite{nocedal2006numerical}.
    
    \item \textbf{Derivative-Free Methods:} These approaches are particularly valuable when gradient information is unavailable or computationally prohibitive. Prominent methods include the Nelder-Mead simplex algorithm, pattern search techniques, and Powell's conjugate direction method \cite{wright1999direct}.
    
    \item \textbf{Metaheuristic Algorithms:} These high-level strategies excel in complex optimization landscapes characterized by non-convexity, non-differentiability, or high dimensionality. Representative techniques include evolutionary algorithms (e.g., genetic algorithms), swarm intelligence methods (e.g., particle swarm optimization), and physics-inspired algorithms (e.g., simulated annealing). See \cite{eiben2003introduction} for a comprehensive overview.
    
    \item \textbf{Convex Optimization Techniques:} For problems exhibiting convexity, specialized methods such as interior-point algorithms and barrier methods provide theoretically guaranteed polynomial-time convergence \cite{boyd2004convex}.
\end{itemize}

Despite their widespread utility, optimization methods face several fundamental challenges including computational complexity in high-dimensional spaces, convergence guarantees in non-convex landscapes, and sensitivity to noise and numerical instabilities. Contemporary research focuses on promising directions such as integration with machine learning paradigms, exploitation of problem structure (e.g., sparsity, separability), and development of parallel and distributed optimization frameworks \cite{bottou2018optimization}.

The increasing computational capabilities and algorithmic innovations have significantly expanded the application domains of optimization techniques. Today, they form the computational backbone of numerous fields including large-scale machine learning model training \cite{goodfellow2016deep}, complex systems design and control, operations research and logistics planning, and quantitative finance and portfolio optimization. The continued development of efficient, robust, and scalable optimization algorithms remains a critical research frontier as we confront increasingly complex problems across science and engineering disciplines.

This paper develops a novel global zero-order optimization algorithm that builds upon recent advances in sampling techniques \cite{Zh25}. The proposed method offers three significant advantages over conventional optimization approaches:
\begin{enumerate}[label=(\roman*)]
    \item \textbf{Minimal regularity requirements:} The algorithm imposes only weak assumptions on the target distribution, making it applicable to a broad class of optimization problems.
    \item \textbf{High-dimensional non-convex optimization:} The method maintains its effectiveness even for challenging, non-convex landscapes in high-dimensional spaces.
    \item \textbf{Initialization robustness:} Unlike many existing methods, our approach demonstrates low sensitivity to initial conditions, ensuring more reliable convergence.
\end{enumerate}

\subsection{Intuitive Foundation}
Let $f:\mathbb{R}^d \rightarrow [0,\infty)$ be a probability density function. Given $N$ independent samples $\{X_i\}_{i=1}^N$ drawn from $f$, the maximum likelihood principle ensures that for any $\varepsilon > 0$ and sufficiently large $N$, the following holds with high probability:
\[
\exists X_k \in \{X_i\}_{i=1}^N \text{ such that } \|X_k - x^*\| < \varepsilon,
\]
where $x^* = \arg\max_{x\in\mathbb{R}^d} f(x)$ denotes the global maximizer, and $\|\cdot\|$ represents the Euclidean norm in $\mathbb{R}^d$. Defining the empirical maximizer:
\begin{equation}
    \xi^*_N = \arg\max_{1\leq i\leq N} f(X_i),
\end{equation}
we establish the following convergence result (see Theorem~\ref{Th24}):
\begin{equation}\label{eq:convergence}
    \lim_{N\to\infty} \left(N^{2/d} \|f(\xi^*_N) - f(x^*)\|_{L^p(\Omega)}\right) < \infty,
\end{equation}
provided $\|\nabla^2 f\|_\infty < \infty$. However, this intuitive estimate becomes computationally infeasible for high-dimensional problems due to the excessively large sample size $N$ required.

Instead, we consider the minimization problem of a potential function via exponential tilting. Let $U:\mathbb{R}^d\to[0,\infty]$ be a measurable function satisfying the integrability condition:
\[
\int_{\mathbb{R}^d}\e^{-U(x)} \dif x < \infty.
\]
For any $\theta\geq 1$, we define the probability measure:
\[
\mu_\theta(\dif x) \propto\e^{-\theta U(x)} \dif x.
\]
A fundamental fact is that:
\[
\lim_{\theta\to\infty}\left(\int_{\mathbb{R}^d}\e^{-\theta U(x)} \dif x\right)^{1/\theta} =\e^{-U_*},
\]
where $U_* = \inf_{x\in\mathbb{R}^d}U(x)$. Moreover, as established in \cite{Hw80}, we have the weak convergence:
\begin{align}\label{DA21}
\mu_\theta \xrightarrow{w} \mu \quad \text{as} \quad \theta\to\infty,
\end{align}
where $\mu$ is supported on the minimizing manifold:
\[
\sN := \{x\in\mathbb{R}^d: U(x) = U_*\}.
\]
The characterization of the measure $\mu$ on $\sN$ was originally established in \cite{Hw80} through an analysis of the Hessian matrix $\nabla^2 U$. We note that while the existence result in \eqref{DA21} guarantees convergence, it does not yield explicit convergence rate estimates. Nevertheless, under the polynomial growth condition: for some $\kappa_0,\kappa_1,m>0$,
\begin{align}\label{UU96}
\kappa_0|x|^{m} \leq U(x) \leq \kappa_1|x|^{m},
\end{align}
where $0$ is the unique minimizer (without regularity assumptions on $U$), we derive quantitative bounds for the convergence of the objective function $U(x)$ to its minimum value $U_*$ along the sampling points, and the convergence of the empirical minimizers to the true minimizer $x_*$. More precisely, for i.i.d. samples $(X^\theta_n)_{n=1,\cdots, N} \sim \mu_\theta$ and any $p,\theta \geq 1$ (see Corollary \ref{Cor29}):
\begin{align}
\left\|1 -\e^{-\min_{1\leq n\leq N} U(X^\theta_n)}\right\|_{L^p(\Omega)} \leq \frac{d}{\theta m} +\e^{-N(\kappa_0/\kappa_1)^{d/m}/(2p)},
\end{align}
and for $\alpha \in (0,\frac{1}{m})$, there exists $C = C(p,d,\alpha,\kappa_0,\kappa_1,m) > 0$ such that (see Theorem \ref{Th29})
\begin{align}
\|X^\theta_n - x_*\|_{L^p(\Omega)} \leq C\theta^{-\alpha}, \quad \theta\geq 1.
\end{align}

In practical implementation, given a sufficiently large inverse temperature $\theta \gg 1$, we generate $N$ independent samples $\mathcal{S}_N = \{X_i\}_{i=1}^N$ from the Gibbs distribution $\mu_\theta(\dif x)$; compute the objective values $\{U(X_i)\}_{i=1}^N$; and finally extract the empirical minimizer $\xi_* = \arg\min_{1\leq i\leq N} U(X_i)$ as an approximation to the global minimizer $x_*$. Since single sampling may not guarantee convergence, we develop iterative schemes (detailed in subsequent sections) to refine the minimizer approximation.

\subsection{Sampling via Stochastic Differential Equations}
The key computational challenge in our algorithm is sampling from the distribution $\mu_\theta$. The most common approach utilizes Langevin diffusion:
\[
   \dif X_t = -\nabla U(X_t)\dif t + \sqrt{2\theta^{-1}}\dif W_t, 
\]
where $W_t$ is a $d$-dimensional standard Brownian motion. This SDE has $\mu_\theta$ as its unique stationary distribution, satisfying the Fokker-Planck equation:
\[
\partial_t \mu_\theta = \nabla \cdot \left(\mu_\theta\nabla U + \theta^{-1}\nabla\mu_\theta\right) = 0.
\]
Under dissipative conditions $\langle \nabla U(x), x\rangle \geq m\|x\|^2 - b$, exponential convergence is guaranteed \cite{EGZ19}:
\[
\mathcal{W}_2(\mu_t, \mu_\theta) \leq C\e^{-\lambda t},
\]
where $\mathcal{W}_2$ denotes the 2-Wasserstein distance. However, this approach has limitations:
\begin{itemize}
    \item Requires $U$ to be differentiable (or subdifferentiable).
    \item Convergence rates depend strongly on $U$'s convexity properties.
\end{itemize}
This approach corresponds to continuous-time stochastic gradient descent.

To overcome these limitations, we employ a transport-based method introduced in \cite{Zh25}. Consider any target distribution $\mu$ on $\mathbb{R}^d$ (not necessarily of Gibbs form). Let $\sigma_t,\beta_t:[0,1]\to[0,1]$ be $C^1$-functions satisfying:
\begin{align}\label{Sig0}
\beta_t = 1 - \sigma_t, \quad \sigma_0 = 1, \quad \sigma_1 = 0, \quad \sigma'_t < 0.
\end{align}
For $\varepsilon, \gamma > 0$, consider the following SDE:
\begin{align}\label{SDE-0}
\dif X_t = b_t(X_t)\dif t + \sqrt{\varepsilon\beta_t'}\dif W_t, \quad X_0 = \xi \sim \mathcal{N}(0,\gamma I_d),
\end{align}
where the drift term $b_t(x)$ is defined for $t \in (0,1)$ and $\eta \sim \mu$ as:
\begin{align}\label{Ex5}
b_t(x) = \frac{\sigma_t'\mathbb{E}[(x - \eta)\e^{-\|x - \beta_t\eta\|^2/(2\ell_t\sigma_t)}]}{\sigma_t\mathbb{E}\e^{-\|x - \beta_t\eta\|^2/(2\ell_t \sigma_t)}} \quad \text{with} \quad \ell_t := \varepsilon\beta_t + \gamma\sigma_t.
\end{align}
The following result is proven in \cite{Zh25}.
\begin{theorem}\label{Th01}
Suppose that $\eta$ is bounded by $K$. Then, for any $\varepsilon,\gamma > 0$ and $\xi \sim \mathcal{N}(0,\gamma I_d)$ independent of $\eta$, the SDE \eqref{SDE-0} admits a unique solution $(X_t)_{t \in [0, 1]}$ such that $X_1 \sim \eta$ and
\[
\|X_t - X_1\|_{L^2(\Omega)} \leq \sigma_t\|\xi - \eta\|_{L^2(\Omega)} + 2\sqrt{d\varepsilon\sigma_t}.
\]
Moreover, for $t \in [0,1)$, the distribution density of $X_t$ is given by
\begin{align}\label{Ex6}
\phi_t(x) := (2\pi\ell_t\sigma_t)^{-d/2}\mathbb{E}e^{-\|x - \beta_t\eta\|^2/(2\ell_t\sigma_t)}.
\end{align}
\end{theorem}

\textbf{Flexibility of the Sampling Scheme:}
Theorem~\ref{Th01} imposes no regularity conditions on the target measure $\mu$, which yields two significant advantages:
\begin{itemize}
    \item It can handle \textit{non-smooth} potential functions that violate standard differentiability requirements.
    \item It naturally admits extensions to manifold-constrained optimization problems, which we will explore in future work.
\end{itemize}

\textbf{Computational Consideration:} 
The drift term computation involves high-dimensional integration, a challenge we resolve in Section~\ref{sec:particle} through a provably convergent Monte Carlo approximation scheme. Specifically, our approach:
\begin{itemize}
    \item Maintains the theorem's generality while ensuring computational tractability.
    \item Provides explicit convergence rates under measurable conditions.
\end{itemize}

\subsection{Adaptive Zooming Algorithm}
Our main algorithm employs an adaptive coordinate-wise scaling strategy. Let $\boldsymbol{\alpha} = (\alpha_1,...,\alpha_d) \in (0,1)^d$ be a multi-index controlling coordinate scaling. Define the operations:
\[
    \boldsymbol{\alpha} \odot x := (\alpha_1x_1,...,\alpha_dx_d).
\]
Given initial point $x_0 \in \mathbb{R}^d$, scaling parameter $\theta \gg 1$, and multi-index $\boldsymbol{\alpha} \in (0,1)^d$, our algorithm primarily involves sampling from the following probability density function via SDE \eqref{SDE-0}:
\[
f(x) = \frac{e^{-\theta U(\boldsymbol{\alpha} \odot x + x_0)}}{\int_{\mathbb{R}^d}\e^{-\theta U(\boldsymbol{\alpha} \odot x + x_0)} \dif x}.
\]
Below is the pseudocode of our algorithm.
\begin{algorithm}[H]
\caption{Adaptive Zooming Optimization}
\begin{algorithmic}[1]
\State \textbf{Input:} Initial point $x_0$, scaling $\theta \gg 1$, multi-index $\boldsymbol{\alpha}_0$
\State \textbf{Output:} Optimizer estimate $x^*$
\For{$k = 0,1,2,...$ until convergence}
    \State Define density: $f_k(x) \propto \exp(-\theta U(\boldsymbol{\alpha}_k \odot x + x_k))$
    \State Sample $X_k \sim f_k$ via SDE \eqref{SDE-0}
    \State Update: $x_{k+1} = x_k + \boldsymbol{\alpha}_k \odot X_k$
    \State Adapt scaling: Update $\boldsymbol{\alpha}_{k+1}$ (see Section 4)
\EndFor
\end{algorithmic}
\end{algorithm}
The algorithm iteratively:
\begin{itemize}
    \item Samples from increasingly localized distributions via \eqref{SDE-0}.
    \item Adapts the search scope through coordinate-wise scaling.
    \item Combines stochastic exploration with deterministic refinement.
\end{itemize}
Section 4 provides complete implementation details, including the adaptive $\boldsymbol{\alpha}_k$ update rule.

\subsection{Organization and Notations}
The remainder of this paper is organized as follows: 
\begin{itemize}
    \item Section~\ref{sec:convergence} establishes non-asymptotic convergence rates for optimization via sampling, providing theoretical guarantees for our approach.
    \item Section~\ref{sec:particle} analyzes the approximation error of our particle-based approximation for SDE~\eqref{SDE-0}, with explicit convergence rates.
    \item Section~\ref{sec:algorithm} presents the complete optimization algorithm, including implementation details and computational complexity analysis.
    \item Section~\ref{sec:numerics} demonstrates the effectiveness of our method through comprehensive numerical experiments on synthetic and real-world benchmark functions.
\end{itemize}

Before proceeding, we introduce the following key notations:
\begin{itemize}
\item The gradient of $U:\mathbb{R}^d\to\mathbb{R}$ is denoted by $\nabla U = (\partial_{x_1}U,\cdots,\partial_{x_d} U)$.
\item For $K>0$, let $B_K := \{x: \|x\| < K\}$ be the open ball in $\mathbb{R}^d$.
\item For $p\geq 1$, let $\|\cdot\|_{L^p(\Omega)}$ denote the $L^p$-norm in probability space $(\Omega,\mathcal{F},\mathbb{P})$.
\item For $\alpha>0$, the Gamma function and lower incomplete Gamma function are defined respectively as
\[
\Gamma(\alpha) := \int^\infty_0 t^{\alpha-1}e^{-t}\dif t, \quad
\gamma(\alpha,x) := \int^x_0 t^{\alpha-1}e^{-t}\dif t, \quad x>0,
\]
satisfying the fundamental relations:
\begin{align}\label{SE1}
\alpha\Gamma(\alpha) = \Gamma(\alpha+1), \quad
\inf_{\alpha\geq 1}\frac{\gamma(\alpha,\alpha)}{\Gamma(\alpha)} \geq \frac{1}{2}.
\end{align}
\item The total variation norm of a signed measure $\mu$ is given by
\[
\|\mu\|_{\mathrm{TV}} := \sup_{A\in\mathcal{B}(\mathbb{R}^d)}|\mu(A)|.
\]
\item For probability measures $\mu,\nu$, the relative entropy is
\[
\mathcal{H}(\mu|\nu) :=
\begin{cases}
\mu(\log f), & \text{if } \mu = f\nu, \\
\infty, & \text{otherwise}.
\end{cases}
\]
\end{itemize}

\section
{Convergence rate estimate of optimization by sampling}\label{sec:convergence}

In this section we present some quantitative convergence rate estimate for finding the minimum of a function by sampling.
We first show the following simple result about seeking the maximum of a function by sampling,
and provide the basic  convergence  rate estimate.
\bl\label{Th41}
Let $N\in\mN$ and $X_1,\cdots, X_N$ be i.i.d. random variables in $\mR^d$ with distribution $\mu$. 
Let $f:\mR^d\to[0,\infty)$ be a measurable function.
Suppose that
$$
f^*:=\sup_{x\in\mR^d}f(x)\in[0,\infty).
$$
Then for any $p>0$, it holds that
\begin{align}\label{RA1}
\mE\left|\sup_{n=1,\cdots N}f(X_n)-f^*\right|^p\leq p\int^{f^*}_0 \e^{-N\mu\{x: f(x)>f^*-r\}}r^{p-1}\dif r.
\end{align}
In particular, if $f:\mR^d\to[0,\infty)$ is continuous, and satisfies
$$
\frak{m}_f:=\int_{\mR^d}f(x)\dif x<\infty, \ \ \ \lim_{|x|\to\infty} f(x)=0,
$$ 
then for $\mu(\dif x)=f(x)\dif x/\frak{m}_f$,
\begin{align}\label{Lim1}
\lim_{N\to\infty}\mE\left|\sup_{n=1,\cdots N}f(X_n)-f^*\right|^p=0.
\end{align}
\el
\begin{proof}
Note that
$$
\left|\sup_{n=1,\cdots N}f(X_n)-f^*\right|=\inf_{n=1,\cdots,N}(f^*-f(X_n))=:Y_N.
$$
Since $Y_N\in[0, f^*]$, we have
\begin{align}\label{SQ1}
\mE\left|\sup_{n=1,\cdots N}f(X_n)-f^*\right|^p=p\int^{f^*}_0 \mP(Y_N\geq r) r^{p-1}\dif r.
\end{align}
For $r>0$, by the independence of $(X_n)_{n=1}^N$, we clearly have
\begin{align*}
\mP(Y_N\geq r)&=\prod_{n=1}^N\mP(f^*-f(X_n)\geq r)=\mP(f^*-f(X_1)\geq r)^N\\
&=\mu\{x: f^*-f(x)\geq r\}^N=(1-\mu\{x: f(x)>f^*-r\})^N,
\end{align*}
which together with \eqref{SQ1} implies  \eqref{RA1} by noting that for $s\in(0,1)$,
\begin{align}\label{RS1}
(1-s)^N\leq \e^{-Ns}.
\end{align}
Moreover, for any $r\in(0,f^*)$, since $f$ is continuous, we have
$$
\mu\{x: f(x)>f^*-r\}=\int_{\{x: f(x)>f^*-r\}}\frac{f(x)}{\frak{m}_f}\dif x\geq\frac{f^*-r}{\frak{m}_f}{\rm vol}\{x: f(x)>f^*-r\}>0.
$$
The limit \eqref{Lim1} now follows by \eqref{RA1} and the dominated convergence theorem.
\end{proof}

\br 
A more explicit and rough estimate for \eqref{RA1} is as follows: for any $\eps\in[0,f^*]$,
\begin{align*}
\mE\left|\sup_{n=1,\cdots N}f(X_n)-f^*\right|^p
&\leq p\left(\int^\eps_0+\int^{f^*}_\eps\right) \e^{-N\mu\{x: f(x)>f^*-r\}}r^{p-1}\dif r\\
&\leq p\int^\eps_0r^{p-1}\dif r+  \left(p\int^{f^*}_\eps r^{p-1}\dif r\right) \e^{-N\mu\{x: f(x)>f^*-\eps\}}\\
&=\eps^p+((f^*)^p-\eps^p)\e^{-N\mu\{x: f(x)>f^*-\eps\}}.
\end{align*}
Below we shall consider concrete cases for obtaining quantitative convergence rate.
\er

When $\mu(\dif x)=f(x)\dif x$, we have the following quantitative convergence rate estimate.
\bt\label{Th24}
Let $f:\mR^d\to[0,\infty)$ be a probability density function with $\|\nabla^2 f\|_\infty<\infty$.
Let $X_1,\cdots, X_N$ be i.i.d. random variables with distribution $\mu(\dif x)=f(x)\dif x$.
Suppose that there is a maximum point $x^*\in\mR^d$ of $f$ so that
$f^*:=f(x^*)=\sup_{x\in\mR^d}f(x)<\infty.$
Then for any $p\geq 1$, 
$$
\varlimsup_{N\to\infty}\left(N^{2/d}\left\|\sup_{n=1,\cdots N}f(X_n)-f^*\right\|_{L^p(\Omega)}\right)\leq 
\frac{\Gamma(\tfrac{2p}d+1)^{1/p}\Gamma(\tfrac d2+1)^{2/d}\kappa_f}{2\pi (f^*)^{2/d}},
$$
where for $\alpha>0$, $\Gamma(\alpha):=\int^\infty_0\e^{-s}s^{\alpha-1}\dif s$ is the Gamma function, and
\begin{align}\label{Gf}
\kappa_f:=\sup_{|\omega|=1}\sup_{x\in\mR^d}[-\tr(\omega\otimes \omega\cdot\nabla^2 f(x))]\geq 0.
\end{align}
In particular, for $d\geq 2$ and $p=d/2$, we have
$$
\varlimsup_{N\to\infty}\left(N\mE\left|\sup_{n=1,\cdots N}f(X_n)-f^*\right|^{d/2}\right)\leq 
\left(\frac{\kappa_f}{2\pi}\right)^{d/2}\frac{\Gamma(\tfrac d2+1)}{f^*}.
$$
\et
\begin{proof}
By Lemma \ref{Th41}, the key point is to make an estimate for
$$
\mu\{x: f(x)>f^*-r\}=\int_{\{x: f(x)>f^*-r\}}f(x)\dif x\geq (f^*-r){\rm vol}\{x: f^*-f(x)<r\}.
$$
Without loss of generality, we may assume $x^*=0$. Otherwise, we may consider the function $\tilde f(x):=f(x+x^*)$.
Since $f\in C^2$, we have $\nabla f(0)=0$.
By Taylor's expansion, we have
\begin{align}\label{SS1}
f^*-f(x)&=f(0)-f(x)+x\cdot\nabla f(0)=-\int^1_0s\int^1_0\tr(x\otimes x\cdot\nabla^2 f(ss'x))\dif s'\dif s\leq 
\frac{\kappa_f}2|x|^2,
\end{align}
where $\kappa_f$ is defined by \eqref{Gf}.
Hence, for $r\in[0,f^*]$,
$$
{\rm vol}\{x: f^*-f(x)<r\}\geq {\rm vol}\left\{x: \kappa_f|x|^2/2<r\right\}
=\left(\frac{2 r}{\kappa_f}\right)^{d/2}\frac{\pi^{d/2}}{\Gamma(\tfrac d2+1)}=K_f r^{d/2},
$$
where $K_f:=(2\pi/\kappa_f)^{d/2}/\Gamma\big(\tfrac d2+1\big).$ Thus,
\begin{align}
\mu\{x: f(x)>f^*-r\}\geq(f^*-r)K_f r^{d/2}.\label{NN1}
\end{align}
For given $\eps\in(0,1)$, we split the right hand side of \eqref{RA1} denoted by $\sI_{p,N}$  into two parts: 
$$
\sI_{p,N}=p\left(\int^{f^*}_{\eps f^*}+\int^{\eps f^*}_0\right) \e^{-N\mu\{x: f(x)>f^*-r\}}r^{p-1}\dif r
=:\sI^{(1)}_{p,N}(\eps)+\sI^{(2)}_{p,N}(\eps).
$$
For $\sI^{(1)}_{p,N}(\eps)$, by \eqref{NN1} with $r=\eps f^*$, we have
$$
\sI^{(1)}_{p,N}(\eps)\leq \left(p\int^{f^*}_{\eps f^*} r^{p-1}\dif r\right)\e^{-N\mu\{x: f(x)>(1-\eps) f^*\}} 
\leq (f^*)^p\e^{-N (1-\eps) f^*K_f(\eps f^*)^{d/2}}.
$$
For $\sI^{(2)}_{p,N}(\eps)$, by \eqref{NN1} with $r\in[0,\eps f^*]$, we also have
\begin{align*}
\sI^{(2)}_{p,N}(\eps)&\leq p\int^{\eps f^*}_0 \e^{-N(f^*-r)K_f r^{d/2}}r^{p-1}\dif r\\
&\leq p\int^{\eps f^*}_0 \e^{-N(1-\eps) f^*K_f r^{d/2}}r^{p-1}\dif r\\
&\leq p\int^\infty_0 \e^{-N(1-\eps)f^*K_f r^{d/2}}r^{p-1}\dif r.
\end{align*}
Using the change of variable $N(1-\eps)f^*K_fr^{d/2}=t$, we obtain
\begin{align*}
\sI^{(2)}_{p,N}(\eps)
&\leq (N(1-\eps) f^*K_f)^{-2p/d}\frac{2p}{d}\int^\infty_0 \e^{-t}t^{2p/d-1}\dif t
=(N(1-\eps) f^*K_f)^{-2p/d}\Gamma\Big(\tfrac{2p}d+1\Big).
\end{align*}
Combining the above calculations, we get for any $\eps\in(0,1)$,
$$
\sI_{p,N}\leq (f^*)^p\e^{-N (1-\eps) f^*K_f(\eps f^*)^{d/2}} +
(N(1-\eps) f^*K_f)^{-2p/d}\Gamma\Big(\tfrac{2p}d+1\Big).
$$
Hence, by \eqref{RA1},
\begin{align*}
N^{2/d}\left\|\sup_{n=1,\cdots N}f(X_n)-f^*\right\|_{L^p(\Omega)}
&\leq (f^*N)^{2/d}\e^{-N (1-\eps) f^*K_f(\eps f^*)^{d/2}/p}\\
&+\Gamma\Big(\tfrac{2p}d+1\Big)^{1/p}/ (K_f f_*(1-\eps))^{2/d}.
\end{align*}
By firstly letting $N\to\infty$ and then $\eps\to 0$, we obtain the desired upper limit bound.
\end{proof}
\br
Estimate \eqref{SS1} implies that $\nabla^2 f(x^*)$ must be negative definite.
Moreover, the $C^2$-assumption of $f$ is only used in \eqref{SS1}. Of course, if we directly assume that
$$
f(x^*)-f(x)\leq\kappa_f |x-x^*|^2/2,\ \ \forall x\in\mR^d,
$$ 
then we have the same conclusion.
If in addition, $\frak{m}_f:=\int_{\mR^d}f(x)\dif x<\infty$, then by considering 
$\tilde f(x)=f(x)/\frak{m}_f$, we have the following convergence rate
$$
\varlimsup_{N\to\infty}\left(N^{2/d}\left\|\sup_{n=1,\cdots N}f(X_n)-f^*\right\|_{L^p(\Omega)}\right)\leq 
\frac{\Gamma(\tfrac{2p}d+1)^{1/p}\Gamma(\tfrac d2+1)^{2/d}}{2\pi}\left(\frac{\frak{m}_f}{f^*}\right)^{2/d}\kappa_f.
$$
\er
The following corollary is a direct application.
\bc
Let $U:\mR^d\to[0,\infty)$ be a $C^2$-function. Suppose that
$$
\lim_{|x|\to\infty} U(x)=\infty,\ \ \frak{m}_U:=\int_{\mR^d}\e^{-U(x)}\dif x<\infty,
$$ 
and 
$$
\kappa_{U}:=
\sup_{|\omega|=1}\sup_{x\in\mR^d}\left[\big(\nabla_\omega\nabla_\omega U(x)-|\nabla_\omega U(x)|^2\big) \e^{-U(x)}\right]<\infty,
$$
where $\nabla_\omega U(x):=\<\omega,\nabla U(x)\>.$
Let $X_1,\cdots, X_N$ be i.i.d. random variables with common 
distribution $\mu(\dif x)\propto\e^{-U(x)}\dif x$. Then for any $p\geq 1$, we have
$$
\varlimsup_{N\to\infty}\left(N^{2/d}\left\|\sup_{n=1,\cdots N}\e^{-U(X_n)}-\e^{-U_*}\right\|_{L^p(\Omega)}\right)\leq 
\frac{\Gamma(\tfrac{2p}d+1)^{1/p}\Gamma(\tfrac d2+1)^{2/d}}{2\pi (\e^{-U_*}/\frak{m}_U)^{2/d}} \kappa_{U},
$$
where $U_*:=\inf_{x\in\mR^d} U(x)$.
\ec
\begin{proof}
We now apply Theorem \ref{Th24} to $f(x)=\e^{-U(x)}/\frak{m}_U$.  By the chain rule, we have
$$
\nabla^2 f(x)=[\nabla U(x)\otimes\nabla U(x)-\nabla^2 U(x)]\e^{-U(x)}/\frak{m}_U.
$$ 
Thus, we clearly have
$$
\kappa_f=\sup_{|\omega|=1}\sup_{x\in\mR^d}[-\tr(\omega\otimes \omega\cdot\nabla^2 f(x))]=\kappa_{U}/\frak{m}_U.
$$
The desired estimate now follows.
\end{proof}
\br
Let $x_*$ be the minimum point of $U$. Then $\nabla U(x_*)=0$ and
$$
\e^{-U(x_*)}\sup_{|\omega|=1}\nabla_\omega\nabla_\omega U(x_*)
\leq \kappa_{U}\leq \e^{-U(x_*)}\sup_{x\in\mR^d}
\sup_{|\omega|=1}\nabla_\omega\nabla_\omega U(x).
$$
\er
\br
If we take $U(x)=|x|^2/2$, then $\frak{m}_U=(2\pi)^{d/2}$ and $\kappa_U=1$.
In this case, we have
$$
\lim_{N\to\infty}\left(N^{2/d}\left\|\sup_{n=1,\cdots N}\e^{-|X_n|^2/2}-1\right\|_{L^2(\Omega)}\right)
\leq\Gamma(\tfrac 4d+1)^{1/2} \Gamma(\tfrac d2+1)^{2/d}.
$$
Since $X_n$ are i.i.d. Gaussian, by tedious calculation, the above inequality indeed takes the equality. 
\er

A natural question arises regarding the optimality of the convergence rate $N^{-2/d}$. Indeed, this rate critically depends on the local behavior of the potential function $U$ near its minimizer $x_*$. In our subsequent analysis, we specifically examine the case where the measure $\mu$ takes the Gibbs form $\mu(\dif x) \propto \e^{-\theta U(x)}\dif x$, with assumptions imposed directly on the potential $U$ rather than its exponential transformation.

From a computational perspective, the dimension-dependent rate $N^{-2/d}$ becomes particularly problematic in high-dimensional settings ($d \gg 1$). In such cases, achieving sufficient precision in locating the minimizer $x_*$ would require prohibitively large sample sizes $N$. To address this fundamental limitation, we introduce a temperature parameter $\theta$ (analogous to its role in Langevin Monte Carlo methods), which enables the following crucial non-asymptotic estimate:

\bt\label{Th28}
Let $U: \mR^d\to[0,\infty]$ be a 
measurable function so that for fixed $\theta\geq 1$, 
$$
\frak{m}_{\theta,U}:=\int_{\mR^d}\e^{-\theta U(x)}\dif x<\infty.
$$ 
Suppose that $U_*:=\inf_{x} U(x)\in(-\infty,\infty)$. For $r\geq 0$, define
\begin{align}\label{Low1}
\ell_\theta(r):=\int_{\{x: 0\leq U(x)-U_*<r\}}\frac{\e^{-\theta U(x)}}{\frak{m}_{\theta,U}}\dif x.
\end{align}
Let 
$X^\theta_1,\cdots, X^\theta_N$ be  i.i.d. random variables with common distribution 
$\mu_\theta(\dif x)\varpropto\e^{-\theta U(x)}\dif x.$ Then for fixed $N\in\mN$ and any $\eps\in(0,1)$ and $p\geq 1$,
\begin{align}\label{DP1}
\Big\|1-\e^{U_*-\inf_{n=1,\cdots N}U(X^\theta_n)}\Big\|_{L^p(\Omega)}\leq
\eps+\e^{-N \ell_\theta(\eps)/p}.
\end{align}
Moreover, if  for some $\alpha,\kappa>0$, 
\begin{align}\label{As1}
\ell_\theta(r/\theta)\geq \kappa r^\alpha,\ \  r\in[0,1],
\end{align}
then  for  any $p\geq 1$,
\begin{align}\label{DP2}
\varlimsup_{N\to\infty}\left(N^{1/\alpha}\Big\|1-\e^{U_*-\inf_{n=1,\cdots N}U(X^\theta_n)}\Big\|_{L^p(\Omega)}\right)\leq 
\frac{\Gamma(\frac p\alpha+1)}{\kappa^{1/\alpha}\theta}.
\end{align}
\et
\begin{proof}
Without loss of generality, we assume $U_*=0$. Otherwise, we may consider
$$
\widetilde U(x):=U(x)-U_*.
$$
Applying  \eqref{RA1} to $f(x)=\e^{-U(x)}$ and $\mu_\theta(\dif x)\propto\e^{-\theta U(x)}\dif x/\frak{m}_{\theta,U}$, we have
$$
\mE\left(1-\e^{-\inf_{n=1,\cdots N}U(X_n)}\right)^p=
\mE\left|1-\sup_{n=1,\cdots N}\e^{-U(X_n)}\right|^p\leq p\int^1_0\e^{-N\delta_{\theta, r}} r^{p-1}\dif r,
$$
where 
$$
\delta_{\theta,r}=\mu_\theta\Big\{x: \e^{-U(x)}>1-r\Big\}=\int_{\e^{-U(x)}>1-r}\frac{\e^{-\theta U(x)}}{\frak{m}_{\theta,U}}\dif x.
$$
Noting that $\e^{-s}\geq 1-s$ for $s>0$,  by the definition \eqref{Low1} of $\ell_\theta$, we have
$$
\delta_{\theta,r}\geq \int_{\{x: U(x)<r\}}\frac{\e^{-\theta U(x)}}{\frak{m}_{\theta,U}}\dif x=\ell_\theta(r),\ \ \forall r\in[0,1].
$$
Hence, for any $\theta\geq 1$ and $\eps\in(0,1)$,
\begin{align}
&\mE\left(1-\e^{-\inf_{n=1,\cdots N}U(X_n)}\right)^p\leq p\int^1_0\e^{-N \ell_\theta(r)} r^{p-1}\dif r\no\\
&\qquad=p\left(\int^{\eps}_0+\int^1_{\eps}\right)\e^{-N \ell_\theta(r)} r^{p-1}\dif r=:I_1+I_2.\label{As2}
\end{align}
We first look at \eqref{DP1}. For $I_1$, we clearly have
$$
I_1\leq p\int^{\eps}_0r^{p-1}\dif r=\eps^p.
$$
For $I_2$, since $\ell_\theta(r)$ is increasing in $r$, we have
\begin{align}
I_2\leq p\left(\int^1_{\eps}r^{p-1}\dif r\right)\e^{-N \ell_\theta(\eps)}\leq \e^{-N \ell_\theta(\eps)}.
\label{As3}
\end{align}
Hence,
$$
\mE\left(1-\e^{-\inf_{n=1,\cdots N}U(X_n)}\right)^p
\leq\eps^p+\e^{-N \ell_\theta(\eps)}.
$$
Thus we obtain \eqref{DP1}. Next we look at \eqref{DP2}. Let $\eps\in(0,1/\theta)$.
For $I_1$, by \eqref{As1} and the change of variable $N\kappa(r\theta)^\alpha=t$, we have
\begin{align}
I_1\leq p\int^\eps_0\e^{-N\kappa (r\theta)^\alpha}r^{p-1}\dif r
&=\frac {p/\alpha}{(N\kappa)^{\frac p\alpha}\theta^p}
\int^{N\kappa(\eps\theta)^\alpha}_0\e^{-t} t^{\frac p\alpha-1}\dif t\leq
\frac{\Gamma(\tfrac p\alpha +1)}{(N\kappa)^{\frac p\alpha}\theta^p}.\label{As4}
\end{align}
Combining \eqref{As2}, \eqref{As3} and \eqref{As4}, we obtain
$$
\mE\left(1-\e^{-\inf_{n=1,\cdots N}U(X^\theta_n)}\right)^p\leq
\frac{\Gamma(\tfrac p\alpha+1)}{(N\kappa)^{\frac p\alpha}\theta^p}+\e^{-N \ell_\theta(\eps)}
$$
and
\begin{align}\label{DS91}
\left\|1-\e^{-\inf_{n=1,\cdots N}U(X^\theta_n)}\right\|_{L^p(\Omega)}\leq
\frac{\Gamma(\frac p\alpha+1)^{1/p}}{(N\kappa)^{\frac 1\alpha}\theta}+\e^{-N \ell_\theta(\eps)/p}.
\end{align}
The proof is complete.
\end{proof}
\noindent {\bf Example 1:} 
Suppose that $U(x) = \kappa |x|^m$ for some $\kappa, m > 0$. Then we have
\begin{align*}
\ell_\theta(r) &= \displaystyle\int_{\{x : \kappa|x|^m < r\}} \e^{-\theta \kappa |x|^m}  \dif x/\displaystyle\int_{\mathbb{R}^d} \e^{-\theta \kappa |x|^m}  \dif x \\
&=\displaystyle\int_{\{x : |x|^m < r\}} \e^{-\theta |x|^m}  \dif x/
\displaystyle\int_{\mathbb{R}^d} \e^{-\theta |x|^m}  \dif x \\
&= \displaystyle\int_0^{r^{1/m}} \e^{-\theta s^m} s^{d-1}  \dif s/\displaystyle\int_0^\infty \e^{-\theta s^m} s^{d-1}  \dif s \\
&= \displaystyle\int_0^r \e^{-t} t^{d/m - 1}  \dif t/\displaystyle\int_0^\infty \e^{-t} t^{d/m - 1}  \dif t \\
&= \gamma(d/m, r\theta)/\Gamma(d/m),
\end{align*}
where the second equality follows from a change of variable, the third uses spherical coordinates, and the fourth is obtained via the substitution $t = \theta s^m$.
Note that $\frac{\gamma(\alpha, \alpha)}{\Gamma(\alpha)} > \frac{1}{2}$ (see \eqref{SE1}). Therefore, for any $\theta \geq d/(m\eps)$,
\[
\ell_\theta(\eps) = \frac{\gamma(d/m, \eps\theta)}{\Gamma(d/m)} \geq \frac{1}{2},
\]
and hence by \eqref{DP1},
\[
\left\|1 - \e^{-\inf_{n=1,\dots,N} U(X^\theta_n)}\right\|_{L^p(\Omega)} \leq \frac{d}{m\theta} + \e^{-N/(2p)}.
\]
Moreover, for $r \in [0,1]$, it is straightforward to show that
\[
\ell_\theta(r) \geq \frac{m r^{d/m}}{d \, \Gamma(d/m)}.
\]
This example illustrates that the curse of dimensionality can be mitigated by choosing a sufficiently large $\theta$.

The following corollary is more general and provides useful information.
\bc\label{Cor29}
Let $0\leq\kappa_0<\kappa_1$ and $m>0$.
Suppose that for some $R\in(0,\infty]$ and $x_*\in\mR^d$, 
\begin{align}\label{UU09}
U(x)=+\infty,\ \ |x-x_*|>R,
\end{align}
and
\begin{align}\label{UU9}
\kappa_0|x-x_*|^{m}\leq U(x)-U(x_*)\leq\kappa_1|x-x_*|^{m},\ \ |x-x_*|<R.
\end{align}
Then for all $p\geq 1$ and $N\in\mN$, it holds that for $\theta\geq d/(m\kappa_1 R^m)$,
$$
\left\|1-\e^{U(x_*)-\inf_{n=1,\cdots N}U(X^\theta_n)}\right\|_{L^p(\Omega)}\leq\frac{d}{\theta m}+\e^{-N(\kappa_0/\kappa_1)^{d/m}/(2p)}.
$$
In particular,
$$
\lim_{\theta\to\infty}\left\|1-\e^{-\inf_{n=1,\cdots N}U(X^\theta_n)}\right\|_{L^p(\Omega)}\leq\e^{-N(\kappa_0/\kappa_1)^{d/m}/p},
$$
and
$$
\lim_{N\to\infty}\left\|1-\e^{-\inf_{n=1,\cdots N}U(X^\theta_n)}\right\|_{L^p(\Omega)}\leq\frac{d}{\theta m}.
$$
\ec
\begin{proof}
Without loss of generality, we assume $x_*=0$ and $U(x_*)=0$.
Let $\eps=\frac{d}{\theta m}$.
Notice that by the change of variable, 
\begin{align*}
\ell_\theta(\eps)&=\int_{\{x: U(x)-U(0)<\eps\}}\e^{-\theta U(x)}\dif x/\int_{\mR^d}\e^{-\theta U(x)}\dif x\\
&\geq\int_{\{x: \kappa_1|x|^{m}<\eps,|x|\leq R\}}\e^{-\theta \kappa_1|x|^{m}}\dif x/\int_{|x|\leq R}\e^{-\theta \kappa_0|x|^{m}}\dif x\\
&=\int^{(\eps/\kappa_1)^{1/m}\wedge R}_0\e^{-\theta \kappa_1s^{m}}s^{d-1}\dif s/\int^R_0\e^{-\theta \kappa_0s^{m}}s^{d-1}\dif s\\
&=\frac{\kappa_0^{d/m}}{\kappa_1^{d/m}}
\int^{(\eps\theta)\wedge(\theta\kappa_1R^{m})}_0\e^{-t}t^{d/m-1}\dif s/\int^{\theta\kappa_0 R^{m}}_0\e^{-t}t^{d/m-1}\dif s\\
&\geq\frac{\kappa_0^{d/m}}{\kappa_1^{d/m}}\cdot
\frac{\gamma(\frac dm, \frac dm\wedge(\theta\kappa_1R^{m})}{\Gamma(\frac dm)}.
\end{align*}
For $\theta$ being large enough so that $\theta \kappa_1 R^{m}>d/m$,
by \eqref{DP1} and \eqref{SE1}, we obtain the desired estimate.
\end{proof}
\br
Conditions \eqref{UU09} and \eqref{UU9} imply that $x_*$ is the unique minimizer of $U$.
Below is a picture of the domain $\{(x,y): |x|^2\leq y \leq 8|x|^2\}$ where the function $U$ stays.
\begin{center}
\begin{tikzpicture}
    \begin{axis}[
        domain=-1.5:1.5, 
        samples=100,
        xlabel={$x$},
        ylabel={$y$},
        axis lines=middle,
        legend pos=outer north east,
        grid=both,
        ymin=0, ymax=8,
        xmin=-1.5, xmax=1.5,
        width=12cm,
        height=6cm
    ]

    \addplot [
        name path=D,
        domain=-1.5:1.5,
        samples=200,
        thick,
        red
    ] {8*abs(x^2)};

    \addplot [
        name path=E,
        domain=-1.5:1.5,
        samples=200,
        draw=none
    ] {abs(x^2)};

    \addplot [
        fill=gray, 
        fill opacity=0.5
    ] fill between[
        of=E and D,
        soft clip={domain=-1.5:1.5}
    ];

    \end{axis}
\end{tikzpicture}
\end{center}
\er

The following lemma shows that under certain assumptions on $ U $, 
the sampling point $ X_\theta $ converges to the true minimizer as $ \theta \to \infty $.
\bt\label{Th29}
Let $U:\mR^d\to[0,\infty)$ be a nonnegative measurable function so that for any $\theta\geq 1$,
$$
\frak{m}_{\theta,U}:=\int_{\mR^d}\e^{-\theta U(x)}\dif x<\infty.
$$
Suppose that there are $\delta\in(0,1)$ and $\kappa_0, \kappa_1, \kappa_2>0$ such that
for $|x-x_*|\leq\delta$,
\begin{align}\label{RA41}
\kappa_0|x-x_*|^m\leq U(x)-U(x_*)\leq \kappa_1|x-x_*|^m,
\end{align}
and for $|x-x_*|>\delta$,
\begin{align}\label{RA42}
U(x)-U(x_*)\geq \kappa_2.
\end{align}
Let $\mu_\theta(\dif x)=\e^{-\theta U(x)}\dif x/\frak{m}_{\theta, U}$ and $X_\theta\sim \mu_\theta$. 
Then for any $p>0$ and $\alpha\in(0,\tfrac 1m)$, 
there is a constant $C=C(d,m,p,\kappa_0,\kappa_1,\kappa_2,\alpha)>0$ such that for all $\theta\geq 1$,
\begin{align}\label{RA01}
\mE|X_\theta-x_*|^p\leq C\theta^{-\alpha p}.
\end{align}
\et
\begin{proof}
Without loss of generality we assume $x_*=0$ and $U_*=0$. Otherwise, we may consider 
$$
\widetilde U(x):=U(x+x_*)-U(x_*).
$$
Fix $p>0$. Note that by definition,
$$
\mE|X_\theta|^p=p\int^\infty_0\mu_\theta\{x: |x|>\lambda\}\lambda^{p-1}\dif\lambda
=\frac{p}{\frak{m}_{\theta, U}}(I_1+I_2),
$$
where 
\begin{align*}
I_1:=\int^\infty_\delta\!\!\int_{\{x: |x|>\lambda\}}\e^{-\theta U(x)}\dif x\lambda^{p-1}\dif\lambda,\ \ \
I_2:=\int^\delta_0\!\!\int_{\{x: |x|>\lambda\}}\e^{-\theta U(x)}\dif x\lambda^{p-1}\dif\lambda.
\end{align*}
For $I_1$, by \eqref{RA42}, we have
\begin{align*}
I_1
&\leq \e^{-(\theta-1)\kappa_2}\int^\infty_\delta\!\!\int_{\{x:|x|>\lambda\}}\e^{-U(x)}\dif x\lambda^{p-1}\dif\lambda\\
&\leq \e^{-(\theta-1)\kappa_2}\int^\infty_0\!\!\int_{\{x:|x|>\lambda\}}\e^{-U(x)}\dif x\lambda^{p-1}\dif\lambda\\
&=\e^{-(\theta-1)\kappa_2}\mE|X_1|^p \frak{m}_{1,U}/p.
\end{align*}
For $I_2$, by \eqref{RA41} and \eqref{RA42}, we have
\begin{align*}
I_2
&=\frac{\delta^p}{p}\int_{\{x:|x|>\delta\}}\e^{-\theta U(x)}\dif x
+\int^\delta_0\!\!\int_{\{x:\lambda<|x|\leq \delta\}}\e^{-\theta U(x)}\dif x\lambda^{p-1}\dif\lambda\\
&\leq \frac{\delta^p\frak{m}_{1,U}}{p\e^{(\theta-1)\kappa_2}}
+\frac1p\int_{\{x:|x|\leq \delta\}}\e^{-\theta\kappa_0|x|^m}|x|^p\dif x.
\end{align*}
Fix $\eps\in(0,\delta)$. 
By the change of variable, we have
\begin{align*}
\int_{\{x:|x|\leq \delta\}}\e^{-\theta\kappa_0|x|^m}|x|^p\dif x
&\leq\int_{\{x:|x|\leq \eps\}}\e^{-\theta\kappa_0|x|^m}|x|^p\dif x+\int_{\{x:\eps<|x|\leq \delta\}}\e^{-\theta\kappa_0|x|^m}|x|^p\dif x\\
&\leq\eps^p\int_{\{x:|x|\leq \eps\}}\e^{-\theta\kappa_0|x|^m}\dif x
+\e^{-\theta\kappa_0\eps^m}\int_{\{x:\eps<|x|\leq \delta\}}|x|^p\dif x\\
&=\frac{2\pi^{d/2}}{\Gamma(d/2)}\left(\eps^p
\int^\eps_0\e^{-\theta\kappa_0 s^m} s^{d-1}\dif s
+\e^{-\theta\kappa_0\eps^m}\int_{\eps}^{\delta}s^{p+d-1}\dif s\right)\\
&\leq\frac{2\pi^{d/2}}{\Gamma(d/2)}\left(\frac{\eps^p\Gamma(d/m)}{(\theta\kappa_0)^{d/m}m}
+\e^{-\theta\kappa_0\eps^m}\frac{\delta^{d+p}}{d+p}\right).
\end{align*}
Combining the above calculations,  we obtain
\begin{align}\label{RA44}
\mE|X_\theta|^p\leq  \frac{\frak{m}_{1,U}(\mE|X_1|^p +\delta^p)}{\frak{m}_{\theta, U}\e^{(\theta-1)\kappa_2} }
+\frac{2\pi^{d/2}}{\frak{m}_{\theta, U}\Gamma(d/2)}\left(\frac{\eps^p\Gamma(d/m)}{(\theta\kappa_0)^{d/m}m}
+\e^{-\theta\kappa_0\eps^m}\frac{\delta^{d+p}}{d+p}\right).
\end{align}
On the other hand, by \eqref{RA41} and the change of variable, we also have
\begin{align*}
\frak{m}_{\theta, U}
&\geq \int_{\{x:|x|<\delta\}}\e^{-\theta U(x)}\dif x\geq\int_{\{x:|x|<\delta\}}\e^{-\theta\kappa_1 |x|^m}\dif x\\
&=\frac{2\pi^{d/2}}{\Gamma(d/2)}\int^\delta_0\e^{-\theta\kappa_1 s^m} s^{d-1}\dif s
=\frac{2\pi^{d/2}\gamma(\tfrac dm,\theta\kappa_1 \delta^m)}{m\Gamma(d/2)(\theta\kappa_1)^{d/m}}.
\end{align*}
If $\theta\kappa_1 \delta^m>d/m$, then by \eqref{SE1}, we further have
$$
\frak{m}_{\theta, U}\geq \frac{\pi^{d/2}\Gamma(\tfrac dm)}{m\Gamma(d/2)(\theta\kappa_1)^{d/m}}.
$$
Substituting this into \eqref{RA44}, we obtain that for $\theta>d/(m\kappa_1 \delta^m)$,
\begin{align*}
\mE|X_\theta|^p
&\leq   \frac{\frak{m}_{1,U}(\mE|X_1|^p +\delta^p)}{\frak{m}_{\theta, U}\e^{(\theta-1)\kappa_2} }
+\frac{2m(\theta\kappa_1)^{d/m}}{\Gamma(\frac dm)}\left(\frac{\eps^p\Gamma(d/m)}{(\theta\kappa_0)^{d/m}m}
+\e^{-\theta\kappa_0\eps^m}\frac{\delta^{d+p}}{d+p}\right)\\
&\leq  \frac{\frak{m}_{1,U}m\Gamma(\frac d2)(\theta\kappa_1)^{d/m}(\mE|X_1|^p +\delta^p)}{\pi^{d/2}\Gamma(\tfrac dm)\e^{(\theta-1)\kappa_2} }
+2\eps^p\left(\tfrac{\kappa_1}{\kappa_0}\right)^{d/m}+
\frac{2m(\theta\kappa_1)^{d/m}\delta^{d+p}}{\e^{\theta\kappa_0\eps^m}\Gamma(\frac dm)(d+p)}.
\end{align*}
From this, by choosing $\eps=\theta^{-\alpha}$ with $\alpha\in(0,\frac1m)$, we derive the desired estimate.
\end{proof}
\br
If $U(x)=|x|^2$, it is easy to see that
$$
\mE|X_\theta|^p\leq C\theta^{-p/2}.
$$
Moreover, conditions \eqref{RA41} and \eqref{RA42} are satisfied when the potential $U$ admits a unique minimizer $x_* \in \mathbb{R}^d$, is $C^2$-smooth in a neighborhood of $x_*$, and has a positive definite Hessian $\nabla^2 U(x_*)$.
\er

\section{Sampling any distribution via SDE}\label{sec:particle}

From a theoretical standpoint, the SDE \eqref{SDE-0} provides a universal framework for sampling arbitrary distributions. However, the drift term $b_t(x)$ depends nonlinearly on the target distribution $\mu$, making its exact computation prohibitively expensive in practice. 
When $\mu$ admits a density $\mu(\dif x) = f(x)\dif x$, we can employ Monte Carlo approximation for $b_t(x)$. Specifically, through the change of variable technique, we obtain:
\bl
Let  $\tau_t:(0,1)\to(0,\infty)$ be any function. The drift $b$ can be written as 
\begin{align}
b_t(x)
=\tau_t\sigma_t'\mathbb{E} \left[\xi\cH_t(x,\xi)\right]/\mathbb{E} \cH_t(x,\xi),\label{Ex1}
\end{align}
where $\xi\sim \cN(0,\mI_d)$ and
\begin{align}
\cH_t(x,y):=\e^{-\sigma_t|x+\tau_t\beta_t y|^2/(2\ell_t)+|y|^2 / 2} f(x-\tau_t\sigma_ty).
\label{Ex21}
\end{align}
Moreover, we also have
\begin{align}
\phi_t(x):=(2\pi\ell_t\sigma_t)^{-d/2}\mE\e^{-|x - \beta_t\eta|^2/(2\ell_t\sigma_t)}=\left(\tau_t\sqrt{\sigma_t/\ell_t}\right)^{d}\mathbb{E} \cH_t(x,\xi).
\label{Ex2}
\end{align}
\el
\begin{proof}
By the change of variable $y=x -\tau_t\sigma_tz$ and using $\sigma_t+\beta_t=1$, we have
\begin{align*}
\mE[(x-\eta)\e^{-|x - \beta_t\eta|^2/(2\ell_t\sigma_t)}]
&=\int_{\mR^d}(x-y) \e^{-|x - \beta_ty|^2/(2\ell_t\sigma_t)}f(y)\dif y\\
&=(\tau_t\sigma_t)^{1+d}\int_{\mR^d}z \e^{-\sigma_t|x+\tau_t \beta_tz|^2/(2\ell_t)}
f(x -\tau_t\sigma_tz)\dif z\\
&=(\tau_t\sigma_t)^{1+d}(2\pi)^{d/2}\mE\left[\xi\cH_t(x,\xi)\right],
\end{align*}
where $\xi\sim \cN(0,\mI_d)$ and $\cH_t(x,y)$ is defined by \eqref{Ex21}.
Similarly, we have
\begin{align}\label{DQ1}
\mE\e^{-|x - \beta_t\eta|^2/(2\ell_t\sigma_t)}
=(\tau_t\sigma_t)^{d}(2\pi)^{d/2}\mE\cH_t(x,\xi).
\end{align}
By the definition \eqref{Ex5} of $b_t(x)$, we complete the proof.
\end{proof}

To evaluate the drift term computationally, we employ a Monte Carlo approximation of the integral. Let $\{\xi_i\}_{i=1}^N \stackrel{\text{i.i.d.}}{\sim} \cN(0,\mI_d)$ be independent Gaussian random vectors. The Monte Carlo approximation $b_t^N(x)$ of the drift $b_t(x)$ is given by:
$$
b_t^N(x) = \tau_t\sigma_t' \sum_{i=1}^N \xi_i \mathcal{H}_t(x,\xi_i)/\sum_{i=1}^N \mathcal{H}_t(x,\xi_i).
$$
For notational convenience, we define the auxiliary function:
\begin{equation}
\mathcal{G}_t(x,y) := (x - \tau_t\sigma_t y)\mathcal{H}_t(x,y)
\end{equation}
which yields the equivalent representation:
\begin{equation}\label{SD1}
b_t^N(x) = \frac{\sigma_t'}{\sigma_t}\left[x - \frac{\sum_{i=1}^N \mathcal{G}_t(x,\xi_i)}{\sum_{i=1}^N \mathcal{H}_t(x,\xi_i)}\right].
\end{equation}
Consider the approximate SDE driven by this drift:
\begin{equation}\label{SDE9}
\dif X_t^N = b_t^N(X_t^N)\dif t + \sqrt{\epsilon\beta_t'}\dif W_t, \quad X_0^N \sim \mathcal{N}(0,\gamma \mI_d).
\end{equation}
Below we make the following assumptions:
\begin{itemize}
    \item $f$ has compact support in the ball $B_K$,
    \item $\inf_{t\in[0,1]} \beta_t' > 0$.
\end{itemize}
In particular, we have the uniform bound:
\begin{equation}
|b_t^N(x)| \leq \frac{|\sigma_t'|}{\sigma_t}(|x| + K).
\end{equation}
By standard SDE theory, there exists a unique strong solution to \eqref{SDE9}. Moreover, we establish the following convergence rate estimate:

\bt\label{Th45}
Let $d\geq 3$ and $\eps\in(0,\gamma]$, $\lambda\geq 9$, where $\ell_t=\eps\beta_t+\gamma\sigma_t$. Let $\tau_t=\sqrt{\lambda\ell_t/\sigma_t}$.
Suppose that $|\eta| \leq K$. 
Let $\mu^N_t$ and $\mu_t$ be the laws of $X^N_t$ and $X_t$, respectively.
Then for all $N \in \mathbb{N}$ and $t\in(0,1)$,
$$
\|\mu^N_t - \mu_t\|_{\var}\leq \frac{\sqrt{8}K\big(9\pi\lambda\gamma/4\big)^{d/4}\e^{3K^2/(\eps\sigma_t)}}{\sqrt{\eps(d-2)N}}
\left(\int_{\mR^d}f^2(z)\dif z\right)^{1/2}.
$$
\et
\begin{proof}
By the classical Csiszár–Kullback–Pinsker (CKP) inequality (see \cite[(22.25)]{Villani2009})
and the entropy formula (see \cite{La21}), we have
\begin{align*}
    \|\mu^N_t - \mu_t\|_{\var}^2 &\leq 2 \cH(\mu_t | \mu^N_t)
    = \int^t_0 \frac{\mathbb{E}|b^N_s(X_s) - b_s(X_s)|^2}{\eps \beta'_s} \dif s.
\end{align*}
Note that by \eqref{SD1} and $\sigma'_s=-\beta_s'\leq0$,
\begin{align*}
&b^N_t(x)-b_t(x)=\frac{\sigma'_t}{\sigma_t}\left[\frac{\mE \cG_t(x,\xi)}{\mE \cH_t(x,\xi)}-\frac{\sum_{i=1}^N\cG_t(x,\xi_i)}{\sum_{i=1}^N\cH_t(x,\xi_i)}\right]
=:-\frac{\beta'_t}{\sigma_t}\cE_t(x,\bxi).
\end{align*}
Thus, by the independence of $X_\cdot$ and $\bxi$ as well as $X_s\sim\phi_s(x)$
(see Theorem \ref{Th01}), we further have
\begin{align}
\|\mu^N_t - \mu_t\|_{\var}^2 \leq \int^t_0 \frac{\beta'_s\mathbb{E}|\cE_s(X_s,\bxi)|^2}{\eps \sigma^2_s} \dif s
=\int^t_0 \frac{\beta'_s}{\eps \sigma^2_s}\left(\int_{\mR^d}\mathbb{E}|\cE_s(x,\bxi)|^2\phi_s(x)\dif x\right) \dif s.
\label{Ex4}
\end{align}
Now, we devote to make an estimate for $\mE|\cE_t(x,\bxi)|^2$. Observe that
\begin{align*}
&\left|\sum_{i=1}^N\cG_t(x,\xi_i)\mE\cH_t(x,\xi)-\sum_{i=1}^N\cH_t(x,\xi_i)\mE\cG_t(x,\xi)\right|\\
&\quad\leq\sum_{i=1}^N|\cG_t(x,\xi_i)|\left|\sum_{i=1}^N\cH_t(x,\xi_i)-\mE\cH_t(x,\xi)\right|\\
&\quad+\sum_{i=1}^N\cH_t(x,\xi_i)\left|\sum_{i=1}^N\cG_t(x,\xi_i)-\mE\cG_t(x,\xi)\right|.
\end{align*}
Since $f$ has support in $B_K$, by the definition of $\cE_t(x,\bxi)$, we have
$$
|\cE_t(x,\bxi)|\leq
\frac{K\left|\frac1N\sum_{i=1}^N\cH_t(x,\xi_i)-\mE\cH_t(x,\xi)\right|+\left|\frac1N\sum_{i=1}^N\cG_t(x,\xi_i)-\mE\cG_t(x,\xi)\right|}
{\mE\cH_t(x,\xi)}.
$$
Since $\cG_t(x,\xi_i)$ and $\cH_t(x,\xi_i)$ are i.i.d. random variable, we further have
\begin{align*}
\|\cE_t(x,\bxi)\|_{L^2(\Omega)}
\leq\frac{K\|\cH_t(x,\xi)\|_{L^2(\Omega)}+\|\cG_t(x,\xi)\|_{L^2(\Omega)}}{\sqrt{N}\mE\cH_t(x,\xi)}
\leq\frac{2K\|\cH_t(x,\xi)\|_{L^2(\Omega)}}{\sqrt{N}\mE\cH_t(x,\xi)}.
\end{align*}
By  \eqref{Ex6}, \eqref{Ex4}, \eqref{Ex2} and $\tau_t=\sqrt{\lambda\ell_t/\sigma_t}$, we obtain
\begin{align*}
\|\mu^N_t - \mu_t\|_{\var}^2 \leq
\int^t_0 \frac{4K^2\beta'_s\lambda^{d/2}}{\eps \sigma^2_s N}\left(\int_{\mR^d}\frac{\mE|\cH_s(x,\xi)|^2}{\mE\cH_s(x,\xi)}\dif x\right) \dif s.
\end{align*}
By \eqref{DQ1} and Jensen's inequality, we have
\begin{align*}
(\tau_t\sigma_t)^{d}(2\pi)^{d/2}\mE\cH_t(x,\xi)\geq 
\e^{-\mE|x - \beta_t\eta|^2/(2\ell_t \sigma_t)},
\end{align*}
and by the change of variable $x-\tau_t\sigma_ty=z$ and $\tau_t=\sqrt{\lambda\ell_t/\sigma_t}$,
\begin{align*}
(\tau_t\sigma_t)^{d}(2\pi)^{d/2}\mE|\cH_t(x,\xi)|^2
&=(\tau_t\sigma_t)^{d}\int_{\mR^d} \e^{-\sigma_t|x +\tau_t\beta_ty|^2/\ell_t+|y|^2/2}f^2(x-\tau_t\sigma_ty)\dif y\\
&=\int_{\mR^d} \e^{-|x -\beta_tz|^2/(\ell_t\sigma_t)+|x-z|^2/(2\lambda\ell_t\sigma_t)}f^2(z)\dif z.
\end{align*}
Thus, by Fubini's theorem, we get
\begin{align*}
\int_{\mR^d}\frac{\mE|\cH_t(x,\xi)|^2}{\mE\cH_t(x,\xi)}\dif x
\leq\int_{\mR^d}\left(\int_{\mR^d} \e^{g_t(x,z)/(2\ell_t\sigma_t)}\dif x\right)f^2(z)\dif z,
\end{align*}
where
$$
g_t(x,z):=\mE|x -\beta_t\eta|^2-2|x -\beta_tz|^2+\tfrac{|x-z|^2}\lambda.
$$
Note that
\begin{align*}
g_t(x,z)
&=(\tfrac1\lambda-1)|x|^2+2\<x,\beta_t(2z-\mE\eta)-\tfrac z\lambda\>+\beta_t^2\mE|\eta|^2+(\tfrac1\lambda-2\beta_t^2)|z|^2\\
&=(\tfrac1\lambda-1)|x+a_\lambda(z)/(\tfrac1\lambda-1)|^2
+\beta_t^2\mE|\eta|^2+(\tfrac1\lambda-2\beta_t^2)|z|^2+\tfrac{|a_\lambda(z)|^2}{1-1/\lambda},
\end{align*}
where 
$$
a_\lambda(z):=\beta_t(2z-\mE\eta)-\tfrac z\lambda.
$$
Since $|\eta|\leq K$, it is easy to see that for $\lambda\geq 9$ and $|z|\leq K$, 
\begin{align*}
g_t(x,z)&\leq
(\tfrac1\lambda-1)|x+a_\lambda(z)/(\tfrac1\lambda-1)|^2
+(1+\tfrac1\lambda+\tfrac{(3+1/\lambda)^2}{1-1/\lambda})K^2\\
&\leq
(\tfrac19-1)|x+a_\lambda(z)/(\tfrac1\lambda-1)|^2
+(1+\tfrac19+\tfrac{(3+1/9)^2}{1-1/9})K^2\\
&=-8|x+a_\lambda(z)/(\tfrac1\lambda-1)|^2/9+12K^2.
\end{align*}
Hence,
\begin{align*}
\int_{\mR^d} \e^{g_t(x,z)/(2\ell_t\sigma_t)}\dif x 
&\leq\int_{\mR^d}\e^{-4|x|^2/(9\ell_t\sigma_t)+6K^2/(\ell_t\sigma_t)}\dif x
=\big(9\pi\ell_t\sigma_t/4\big)^{d/2}\e^{6K^2/(\ell_t\sigma_t)}.
\end{align*}
Combining the above calculations, we obtain
$$
\|\mu^N_t - \mu_t\|_{\var}^2
\leq\int^t_0 \frac{4K^2\beta'_s\lambda^{d/2}}{\eps \sigma^2_s N}
\big(9\pi\ell_s\sigma_s/4\big)^{d/2}\e^{6K^2/(\ell_s\sigma_s)}\dif s\int_{\mR^d}f^2(z)\dif z.
$$
Noting that for $\eps\in(0,\gamma]$, by $\sigma_s+\beta_s=1$,
$$
\gamma\geq \ell_s=\eps+(\gamma-\eps)\sigma_s\geq\eps,
$$
we have
\begin{align*}
\|\mu^N_t - \mu_t\|_{\var}^2
&\leq\frac{4K^2\big(9\pi\lambda \gamma/4\big)^{d/2}}{\eps N}\e^{6K^2/(\eps\sigma_t)}\int^t_0 
\beta_s'\sigma_s^{d/2-2}\dif s\int_{\mR^d}f^2(z)\dif z\\
&\leq\frac{8K^2\big(9\pi\lambda \gamma/4\big)^{d/2}}{\eps N (d-2)}\e^{6K^2/(\eps\sigma_t)}
\int_{\mR^d}f^2(z)\dif z.
\end{align*}
The proof is complete.
\end{proof}

Using the above theorem and applying Euler's discretization scheme to the time variable, 
we propose the following algorithm for sampling from the probability density function $f$:

\begin{algorithm}[H]
\caption{ Monte Carlo Sampling via SDEs with Time-Dependent Drifts}
\begin{algorithmic}[1]
\State \textbf{Input:} 
    \begin{itemize}
        \item Target distribution: $f: \mathbb{R}^d \to \mathbb{R}_+$ (integrable non-negative measurable)
        \item $N = 1000$: Monte Carlo samples per iteration
        \item $M = 10$: Time steps
        \item $\lambda > 0$: Importance sampling intensity
        \item $\gamma > 0$: Initial noise amplitude
        \item $\epsilon > 0$: Brownian motion coefficient
    \end{itemize}
\State  \textbf{Initialize}:
    \begin{itemize}
        \item $\Delta t \gets \frac{1}{M}$: Time step
        \item $\sigma_t \gets 1 - t$, $\beta_t \gets t$: Time-dependent parameters
        \item $\{\xi_i\}_{i=1}^N \stackrel{iid}{\sim} \mathcal{N}(0,\mI_d)$: Standard normal vectors
        \item $X_0 \sim \mathcal{N}(0, \gamma\mI_d)$: Initial state
    \end{itemize}

\State\textbf{Iterate} for $k = 0$ to $M-1$:
    \begin{enumerate}[label=(\roman*)]
        \item Update time parameters:
       $$
            t_k \gets k \Delta t, \ \
            \ell_k \gets \epsilon \beta_{t_k} + \gamma \sigma_{t_k}, \ \
            \tau_k \gets \sqrt{\lambda \ell_k / \sigma_{t_k}}
$$
 
        \item Compute drift term:
        $b_k \gets \tau_k \sigma_{t_k}' \cdot \sum\limits_{i=1}^N \xi_i w_i(X_k)/\sum\limits_{i=1}^N w_i(X_k)$, where
      $$
      w_i(x) \gets \exp\left(-\frac{|x + \tau_k \beta_{t_k} \xi_i|^2}{2\ell_k} + \frac{|\xi_i|^2}{2}\right) f(x - \tau_k \sigma_{t_k} \xi_i)
      $$  
        \item Generate noise: $\zeta_k \sim \mathcal{N}(0,\mI_d)$
        
        \item Update state:
        \begin{align*}
            X_{k+1} &\gets X_k + b_k \Delta t + \sqrt{\epsilon \beta_{t_k}' \Delta t} \cdot \zeta_k
        \end{align*}
    \end{enumerate}
 \State    \textbf{Output}: $X_M$ (final state)
\end{algorithmic}
\end{algorithm}

\section{Sampling-based zero-order optimization algorithm}\label{sec:algorithm}

In this section, based on Corollary  \ref{Cor29}, Theorems~\ref{Th29} and \ref{Th45}, we develop a zero-order optimization algorithm. To this end, we first introduce necessary notation.

Let $\boldsymbol{\alpha} = (\alpha_1, \cdots, \alpha_d) \in (0,1)^d$ be a multi-index that controls the scaling of the target function along different coordinates. For any $x = (x_1, \cdots, x_d) \in \mathbb{R}^d$, we define the following operations:
\begin{align*}
\boldsymbol{\alpha} \odot x := (\alpha_1 x_1, \cdots, \alpha_d x_d), \ \
|\boldsymbol{\alpha}| := \alpha_1 \cdots \alpha_d.
\end{align*}
Let $U: \mathbb{R}^d \to [0, \infty]$ denote the objective function whose minimum we seek. Throughout this section, we maintain the fundamental assumption:
$$
\int_{\mathbb{R}^d} \e^{-U(x)} \, \dif x < \infty.
$$

Let $ x_0 \in \mathbb{R}^d $, $ \theta \gg 1 $ and $ \boldsymbol{\alpha} \in (0,1)^d $. 
Our algorithm mainly consists of sampling from the following probability density function:
\begin{align}\label{FF1}
f(x) := \frac{ \exp \big( -\theta \, U(\boldsymbol{\alpha} \odot x + x_0) \big) }{ \int_{\mathbb{R}^d} \exp \big( -\theta \, U(\boldsymbol{\alpha} \odot y + x_0) \big) \, \mathrm{d}y }.
\end{align}

\begin{figure}[h]
\centering
\includegraphics[width=4in, height=2in]{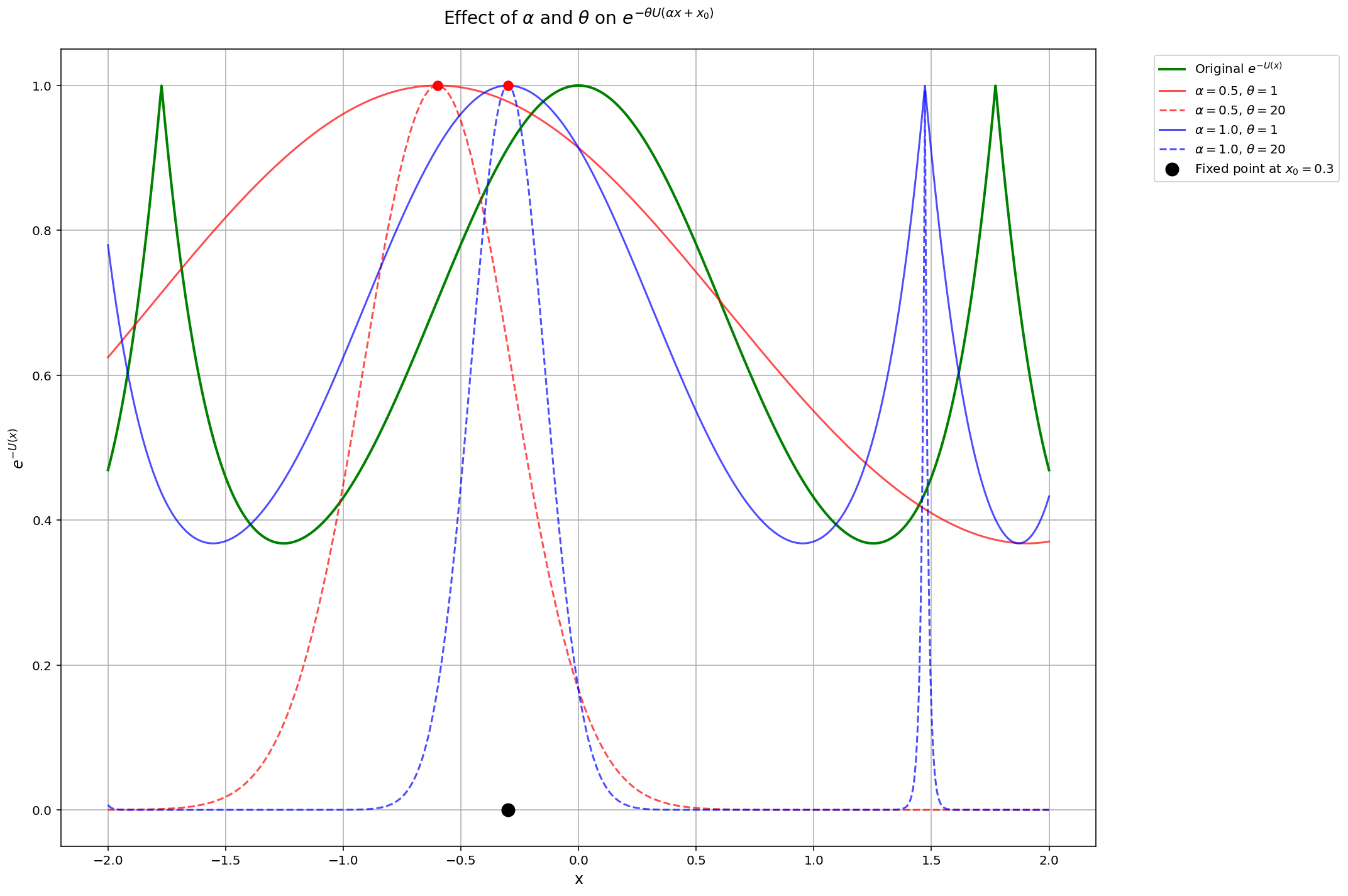}
\caption{Illustration of parameter effects. Larger $\theta$ values yield steeper density functions $f(x)$, while smaller $\boldsymbol{\alpha}$ components produce flatter densities.}
\label{fig:params}
\end{figure}

\bl
Let $f$ be given in \eqref{FF1}. Under \eqref{UU9}, for $\theta\geq d/(m\kappa_1 R^{m})$, we have
\begin{align}\label{DA1}
\int_{\mR^d}f^2(z)\dif z\leq4m|\balpha|(\theta\kappa^2_1/(2\kappa_0))^{d/m}/\Gamma(\tfrac dm).
\end{align}
\el
\begin{proof}
By the definition and  the change of variable, we have
\begin{align*}
\int_{\mR^d}f^2(z)\dif z&=\int_{\mR^d}\e^{-2\theta U(\balpha  \odot x+x_0)}\dif x/\left(\int_{\mR^d}\e^{-\theta U(\balpha  \odot x+x_0)}\dif x\right)^2\\
&=|\balpha|\int_{\mR^d}\e^{-2\theta U(x)}\dif x/\left(\int_{\mR^d}\e^{-\theta U(x)}\dif x\right)^2\\
&=|\balpha|\int_{\mR^d}\e^{-2\theta U(x)}\dif x/\left(\int_{\mR^d}\e^{-\theta U(x)}\dif x\right)^2\\
&\leq|\balpha|\int_{|x|\leq R}\e^{-2\theta\kappa_0|x|^m}\dif x/\left(\int_{|x|\leq R}\e^{-\theta\kappa_1|x|^m}\dif x\right)^2\\
&=|\balpha|\int^{R}_0\e^{-2\theta \kappa_0s^{m}}s^{d-1}\dif s
/\left(\int^{R}_0\e^{-\theta \kappa_1s^{m}}s^{d-1}\dif s\right)^2\\
&=\frac{|\balpha|m(\theta\kappa_1)^{2d/m}}{(2\theta\kappa_0)^{d/m}}
\int^{2\theta\kappa_0 R^{m}}_0\e^{-t}t^{d/m-1}\dif t
/\left(\int^{\theta\kappa_1 R^{m}}_0\e^{-t}t^{d/m-1}\dif t\right)^2\\
&\leq|\balpha|(\theta\kappa^2_1/(2\kappa_0))^{d/m}
\Gamma(\tfrac dm)/\gamma(\tfrac dm, \theta\kappa_1 R^{m})^2.
\end{align*}
In particular, for $\theta\kappa_1 R^{m}\geq\tfrac dm$, by \eqref{SE1}, we obtain the desired estimate.
\end{proof}
\br
The convergence rate in Theorem \ref{Th45} depends on the $L^2$-norm of the density. From \eqref{DA1}, it can be observed that smaller values of $|\balpha|$ lead to faster convergence in the particle approximation.
\er

The following lemma is direct by the change of variable.
\bl\label{Le32}
Let $X\sim \e^{-\theta U(\balpha  \odot x+x_0)}/\int_{\mR^d}\e^{-\theta U(\balpha  \odot x+x_0)}\dif x.$
Then 
$$
\boldsymbol{\alpha}\odot X+x_0\sim \e^{-\theta U(x)}\dif x/\int_{\mR^d}\e^{-\theta U(x)}\dif x.
$$
\el
Our proposed algorithm implements an adaptive zooming strategy through the iterative update rule:
\begin{equation}\label{eq:update_rule}
    x_{k+1} = x_k + \boldsymbol{\alpha}_k \odot X_k,
\end{equation}
where
\begin{itemize}
    \item $\boldsymbol{\alpha}_k \in \mathbb{R}^d$ denotes the adaptive zooming parameter vector that dynamically controls the exploration scope.
    \item $X_k$ is sampled from the probability distribution:
    \begin{equation}\label{eq:sampling_dist}
        f_k(x) \propto \exp\left(-\theta U(\boldsymbol{\alpha}_k \odot x + x_k)\right).
    \end{equation}
\end{itemize}

From Lemma~\ref{Le32}, we derive the key distributional property:
\begin{equation}\label{eq:target_dist}
    x_{k+1} \sim \frac{\exp\left(-\theta U(x)\right)}{\int_{\mathbb{R}^d} \exp\left(-\theta U(x)\right) \dif x}.
\end{equation}

The algorithm's convergence is governed by two fundamental mechanisms:
\begin{enumerate}
    \item \textbf{Progressive zooming}: The parameters $\boldsymbol{\alpha}_k$ decay asymptotically ($\boldsymbol{\alpha}_k \to \mathbf{0}$ as $k \to \infty$), enabling increasingly refined local search.
    
    \item \textbf{Temperature-controlled concentration}: For large $\theta \gg 1$, the distribution exhibits exponential concentration around the global minimizer $x_*$, ensuring (see Theorem \ref{Th29}):
    \begin{itemize}
        \item $x_k \to x_*$ (convergence to optimum).
        \item $\|X_k\| \to 0$ (vanishing exploration steps).
    \end{itemize}
\end{enumerate}

\subsection{Zooming Parameter Update Strategies}
\label{sec:zooming}
In the above algorithm,
the adaptive control of the zooming parameter $\balpha_k \in \mathbb{R}^d$ governs the algorithm's exploration-exploitation trade-off. We formalize two update methodologies: 
\begin{itemize}
\item[(i)] Exponential Decay Update (EDU):
\begin{equation}
\label{eq:exp_update}
\balpha_{k+1} = \balpha^{1+k}, \quad \text{where} \quad \balpha \sim \mathcal{U}\left([\alpha_{\text{min}}, \alpha_{\text{max}}]^d\right),
\end{equation}
where $\alpha_{\text{min}} > 0$ prevents degenerate collapse, and $\alpha_{\text{max}} < 1$ guarantees contraction, $\cU$ is the uniform distribution.

\item[(ii)] Sample-adaptive Variance Update (SVU): for given samples $\cS_N = \{Y_n\}_{n=1}^N$ at iteration $k$,
we compute component-wise variances:
\[
\sigma_{ki} := \var(\{Y_{ni}\}_{n=1}^N), \quad \forall i \in \{1,\dots,d\}.
\]
Then we update $\balpha_k$ by
\begin{equation}
\label{eq:SVU}
\balpha_{k+1} = \balpha_k \otimes \bsigma_k,
\end{equation}
where $\bsigma_k$ is the normalized variance vector:
\[
\bsigma_k := \left(\frac{\sigma_{k1}}{\|\bsigma_k\|}, \dots, \frac{\sigma_{kd}}{\|\bsigma_k\|}\right), \quad \|\bsigma_k\| := \sqrt{\sum_{i=1}^d \sigma_{ki}^2}.
\]
\end{itemize} 
In the above SVU, $\sigma_{ki}$ reflects the concentration of each component.

\subsection{Pseudocode of algorithm} Below is the concrete algorithm.

\begin{algorithm}[H]
\caption{ Global Optimization via Density Sampling}
\begin{algorithmic}[1]
\State \textbf{Input}: 
    \begin{itemize}
        \item Objective function: $U: \mathbb{R}^d \to \mathbb{R}_+$ (non-negative measurable)
        \item $N=10$: Samples per iteration
        \item $M = 200$: Maximum iterations
        \item $\theta = 10^{100}$: Sampling intensity
        \item $\alpha_{\text{min}} = 0.1$, $\alpha_{\text{max}} = 1$: Scaling bounds
    \end{itemize}
\State \textbf{Initialize}:
    \begin{itemize}
        \item $U_{\min} \gets +\infty$: Current minimum value
        \item $x_0 \gets \mathbf{0} \in \mathbb{R}^d$: Initial estimate
        \item $\boldsymbol{\alpha}_0 \gets \mathbf{1} \in \mathbb{R}^d$: Initial scaling vector
    \end{itemize}

\State \textbf{Iterate} for $k = 0$ to $M-1$:
    \begin{enumerate}[label=(\roman*)]
        \item Define scaled function: $U_k(x) \gets U(\boldsymbol{\alpha}_k \odot x + x_k)$
        
        \item Sample points $\mathcal{S}_N = \{Y_1,\ldots,Y_N\}$ using \textbf{Algorithm 2} with:
        \begin{align*}
            f_k(x) \propto \exp(-\theta U_k(x))
        \end{align*}
        
        \item Find best sample: $m \gets \operatorname{argmin}_{1\leq n\leq N} U_k(Y_n)$
        
        \item \textbf{If} $U_k(Y_m ) < U_{\min}$ \textbf{then}:
        \begin{itemize}
            \item $U_{\min} \gets U_k(Y_m )$
            \item $x_{k+1} \gets \boldsymbol{\alpha}_k \odot Y_m + x_k$
        \end{itemize}
        
        \item Update zooming parameters $\balpha_k$ by EDU \eqref{eq:exp_update} or  SVU \eqref{eq:SVU}.
    \end{enumerate}

\State \textbf{Output}: $(x_M, U(x_M))$ (estimated minimizer and value)
\end{algorithmic}
\end{algorithm}


\section{Numerical Experiments}\label{sec:numerics}

In this section, we evaluate the proposed algorithm on a suite of standard benchmark functions commonly employed in the global optimization literature (see \cite{JY13}). These test problems are carefully selected to assess algorithmic performance across several critical dimensions of optimization difficulty:

\begin{itemize}
    \item \textbf{Modality:} This refers to the number of local minima and maxima in the objective function. Multimodal functions pose significant challenges due to the increased likelihood of algorithms becoming trapped in suboptimal local minima.

    \item \textbf{Separability:} A function $U(x_1,\dots,x_d)$ is said to be separable if it can be expressed as
    $$
    U(x_1,\dots,x_d) = \sum_{i=1}^d u_i(x_i).
    $$
    In such cases, each variable can be optimized independently, greatly reducing the problem's complexity.

    \item \textbf{Differentiability:} The differentiability of a function determines whether gradient-based optimization methods, such as Gradient Descent (GD) or Stochastic Gradient Descent (SGD), are applicable.

    \item \textbf{Discontinuity:} Functions exhibiting discontinuities typically preclude the use of gradient-based techniques, necessitating derivative-free approaches such as genetic algorithms, simulated annealing, or Bayesian optimization.

    \item \textbf{Scalability:} This dimension evaluates the algorithm's capacity to maintain performance as the problem dimensionality or structural complexity increases. Scalability is influenced by both the intrinsic properties of the objective function and the design of the optimization method.
\end{itemize}

\subsection{Optimization of high dimensional benchmark functions}
In this subsection we provide nine examples for testing our algorithm.
For each benchmark function, we present three visualizations to illustrate performance:

\begin{itemize}
    \item \textbf{Left Panel:} Contour plot of the two-dimensional function to reveal the shape of the optimization landscape.

    \item \textbf{Middle Panel:} Convergence plot showing the objective function value versus iteration number, with key algorithmic parameters annotated.

    \item \textbf{Right Panel:} Optimization trajectory projected onto the first two coordinates, providing insights into the search dynamics.
\end{itemize}

\bx\rm
\texttt{(Schwefel Function)}: differentiable, separable, scalable, multimodal.

The Schwefel function is notorious for its numerous local minima and a global minimum situated in a narrow basin, posing significant challenges for optimization algorithms to converge to the optimal solution.
For $d \geq 2$ with $x= (x_1, \dots, x_d)$, the $d$-dimensional Schwefel function is defined as:
$$
U(x) = 418.9829 \cdot d - \sum_{i=1}^{d} x_i \sin\left(\sqrt{|x_i|}\right), \quad x\in [-500, 500]^d,
$$
where the global minimum $U(x_*) = 0$ is achieved at $x_* = (420.9687, \dots, 420.9687)$, a point distant from the origin. Our numerical experiments yield the following results:
\begin{center}
\begin{minipage}{\textwidth}
\centering
\includegraphics[width=1.9in, height=1.4in]{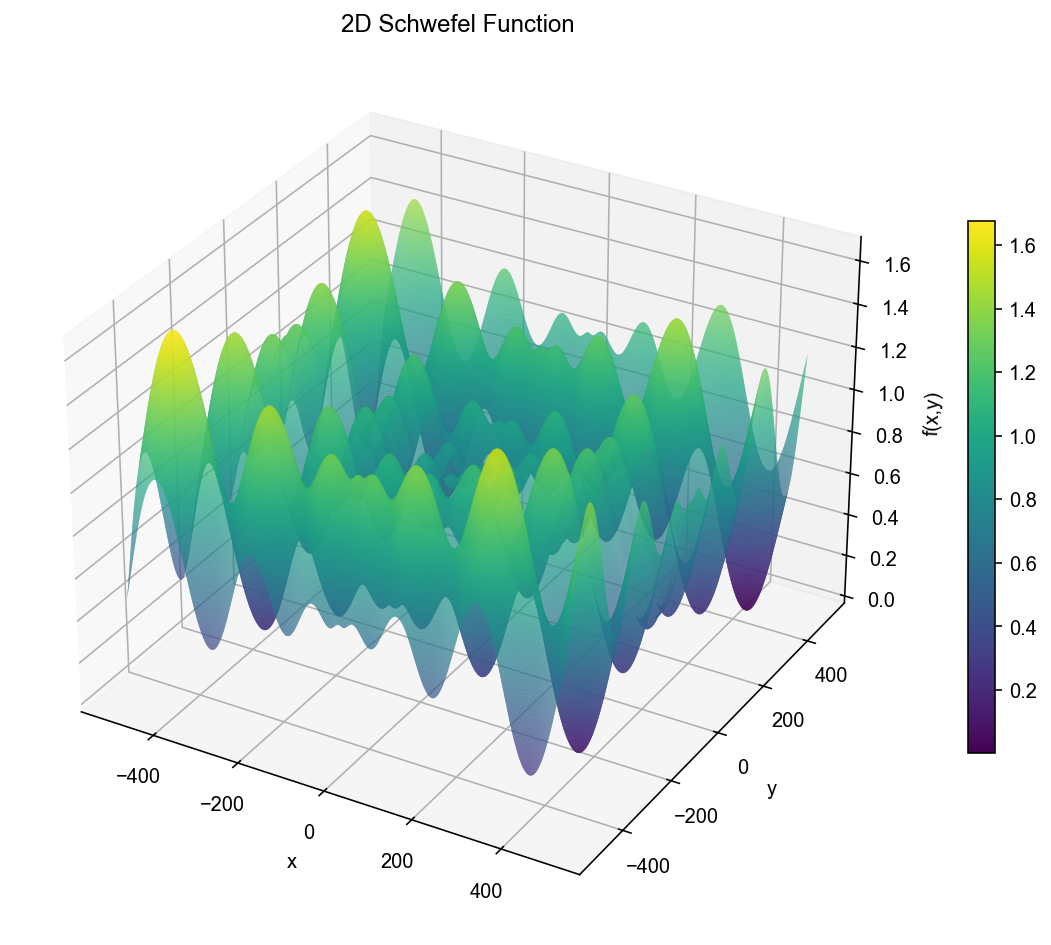}
\includegraphics[width=1.9in, height=1.4in]{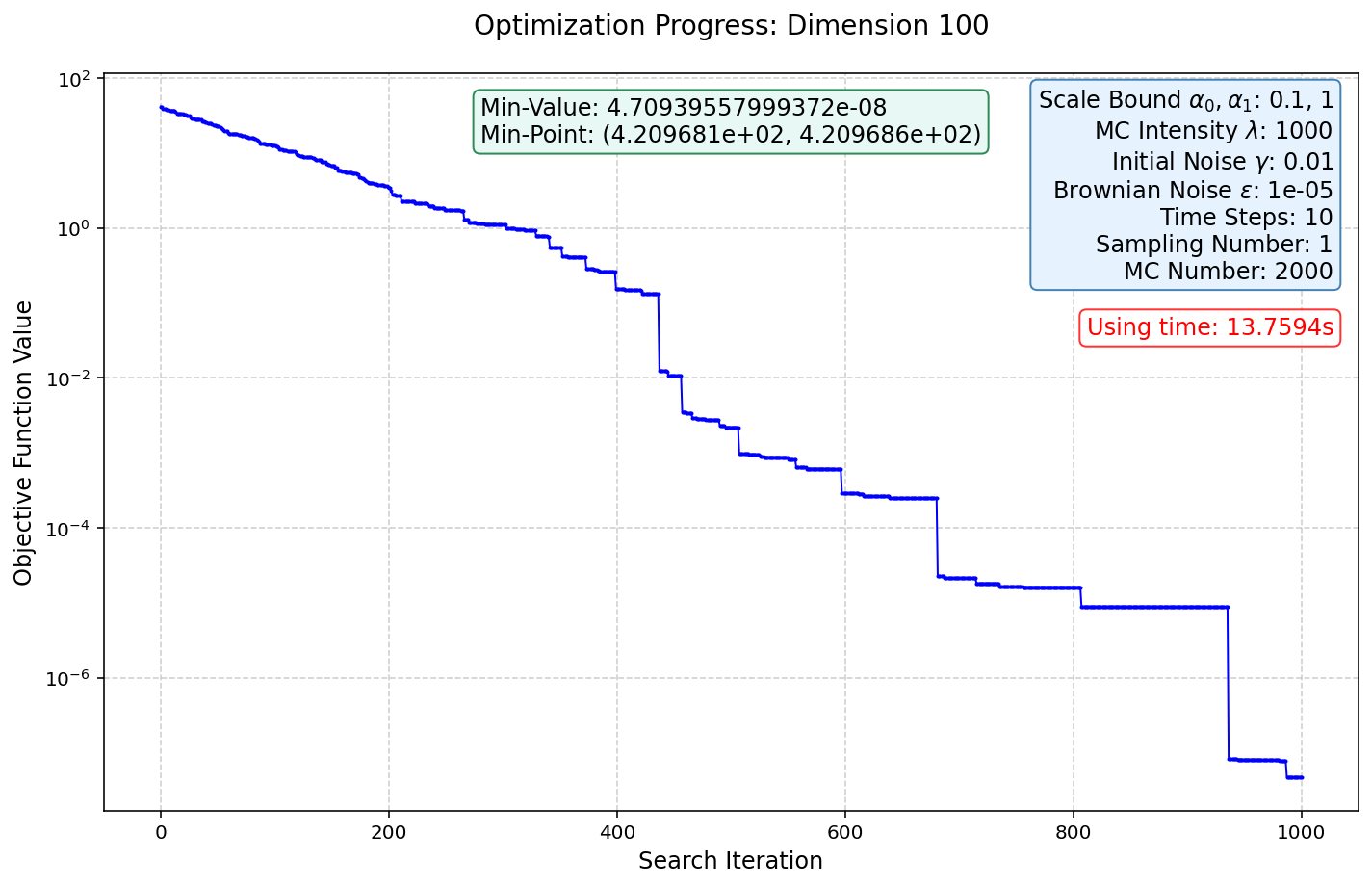}
\includegraphics[width=1.9in, height=1.4in]{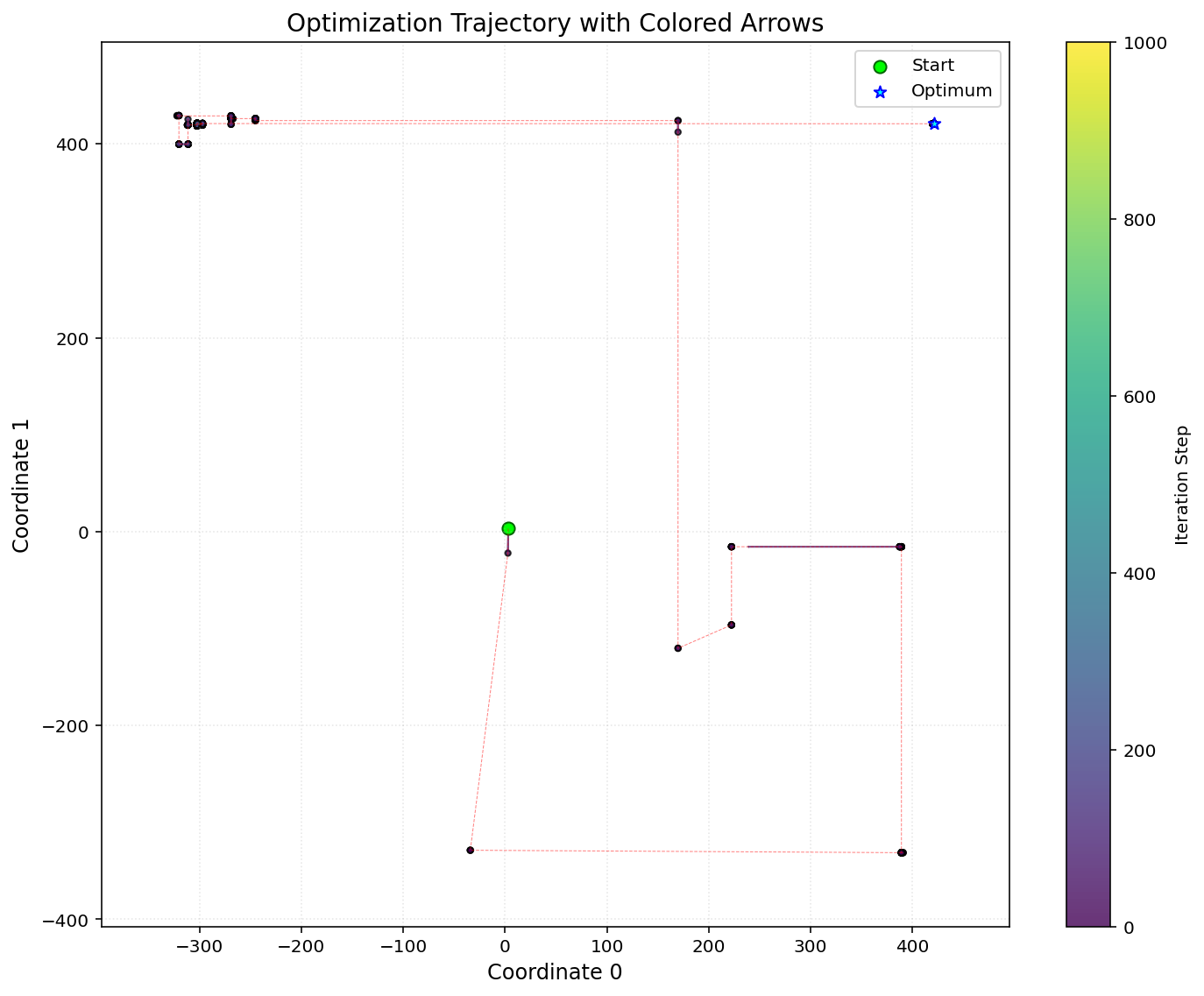}
\captionof{figure}{Schwefel Function}
\end{minipage}
\end{center}
\ex

\bx\rm
\texttt{Rosenbrock Function}: differentiable, non-separable, scalable, unimodal.

This classic benchmark features a narrow curved valley, which challenges gradient-based methods due to poor conditioning and slow convergence. The $d$-dimensional Rosenbrock function is defined as:
$$
U(x) = \sum_{i=1}^{d-1} \left[ 100(x_{i+1} - x_i^2)^2 + (1 - x_i)^2 \right],
$$
where the global minimum $U(x_*) = 0$ is achieved at $x_* = (1, \dots, 1)$.
The numerical experimental results are presented below.
\begin{center}
\begin{minipage}{\textwidth}
\centering
\includegraphics[width=1.9in, height=1.4in]{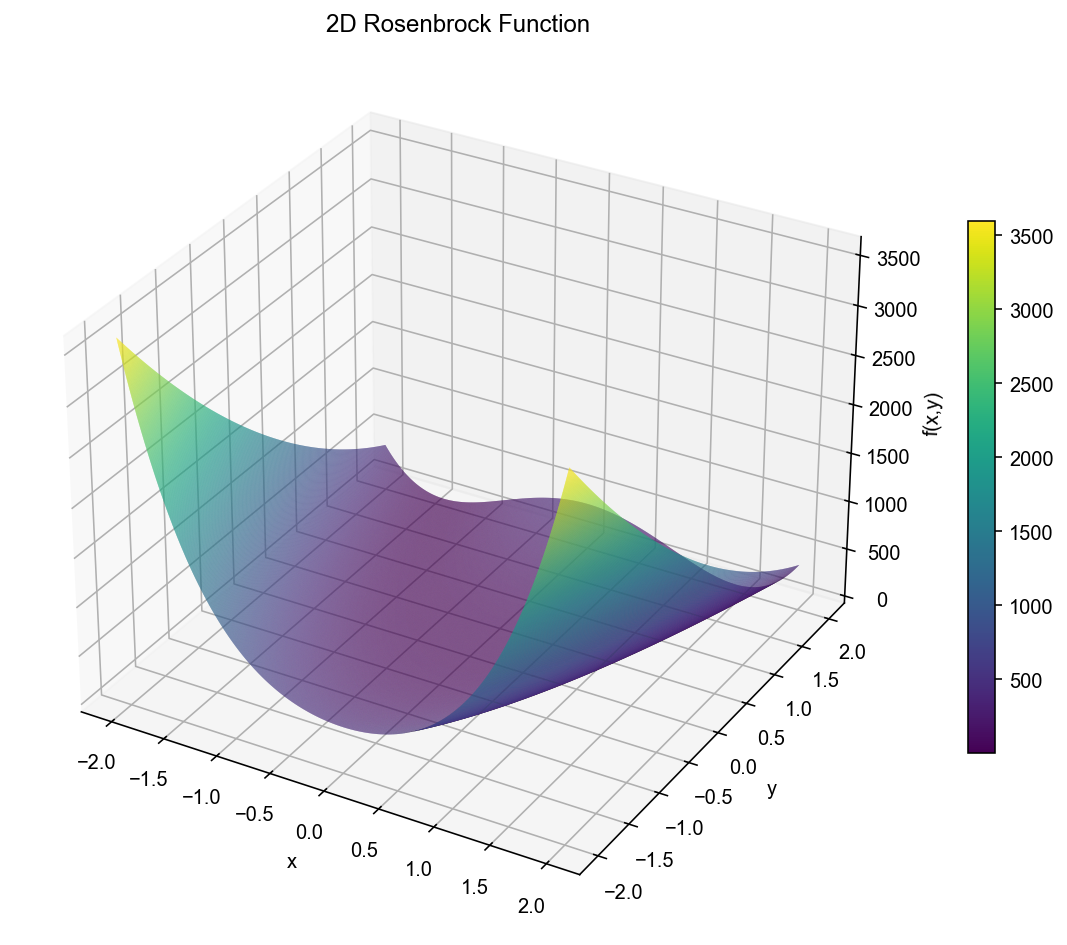}
\includegraphics[width=1.9in, height=1.4in]{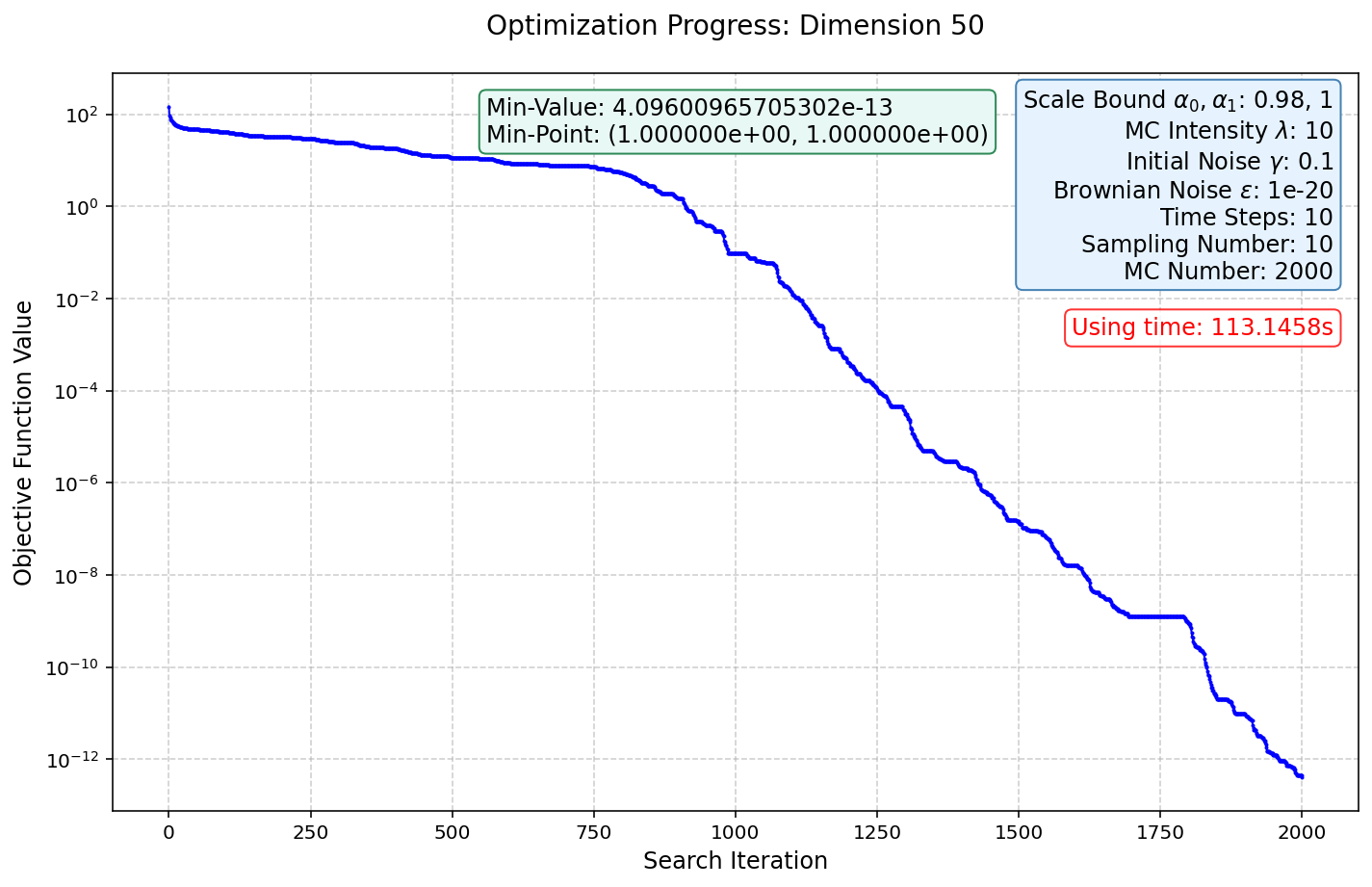}
\includegraphics[width=1.9in, height=1.4in]{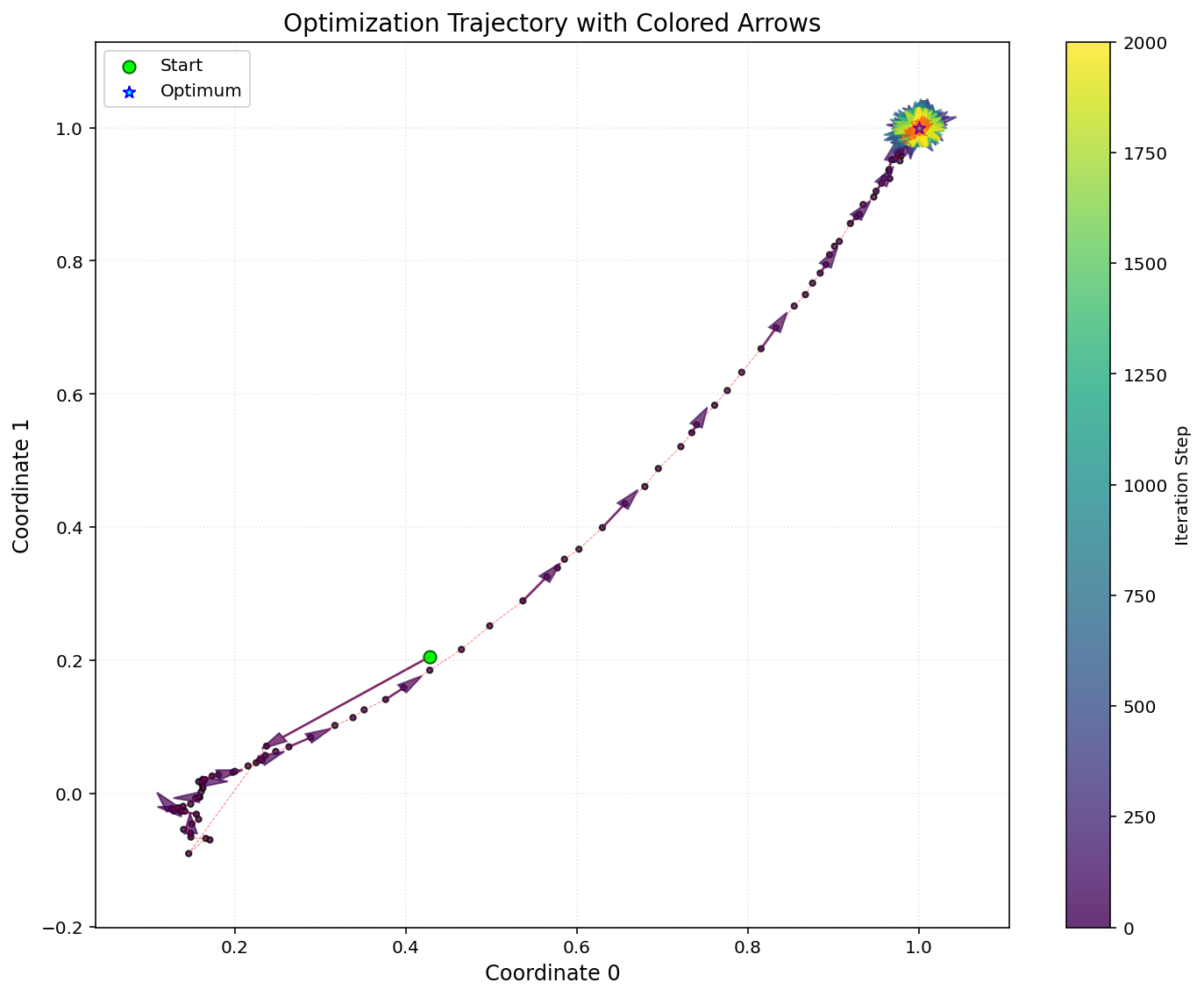}
\captionof{figure}{Rosenbrock Function}
\end{minipage}
\end{center}
\ex

\bx\rm
\texttt{Ackley Function}: differentiable, non-separable, scalable, multimodal.

The Ackley function contains a nearly flat outer region and a sharp global minimum in the center, making it ideal for testing convergence properties of global optimizers.
The $d$-dimensional Ackley function is defined as:
$$
U(x) = -20 \exp \left(-0.2 \sqrt{\frac{1}{d} \sum_{i=1}^{d} x_i^2} \right) - \exp \left(\frac{1}{d} \sum_{i=1}^{d} \cos(2\pi x_i) \right) + 20 + \e,
$$
where the global minimum $U(x_*) = 0$ is achieved at $x_* = (0, \dots, 0)$.
Below are the results of our numerical experiments. 
\begin{center}
\begin{minipage}{\textwidth}
\centering
\includegraphics[width=1.9in, height=1.4in]{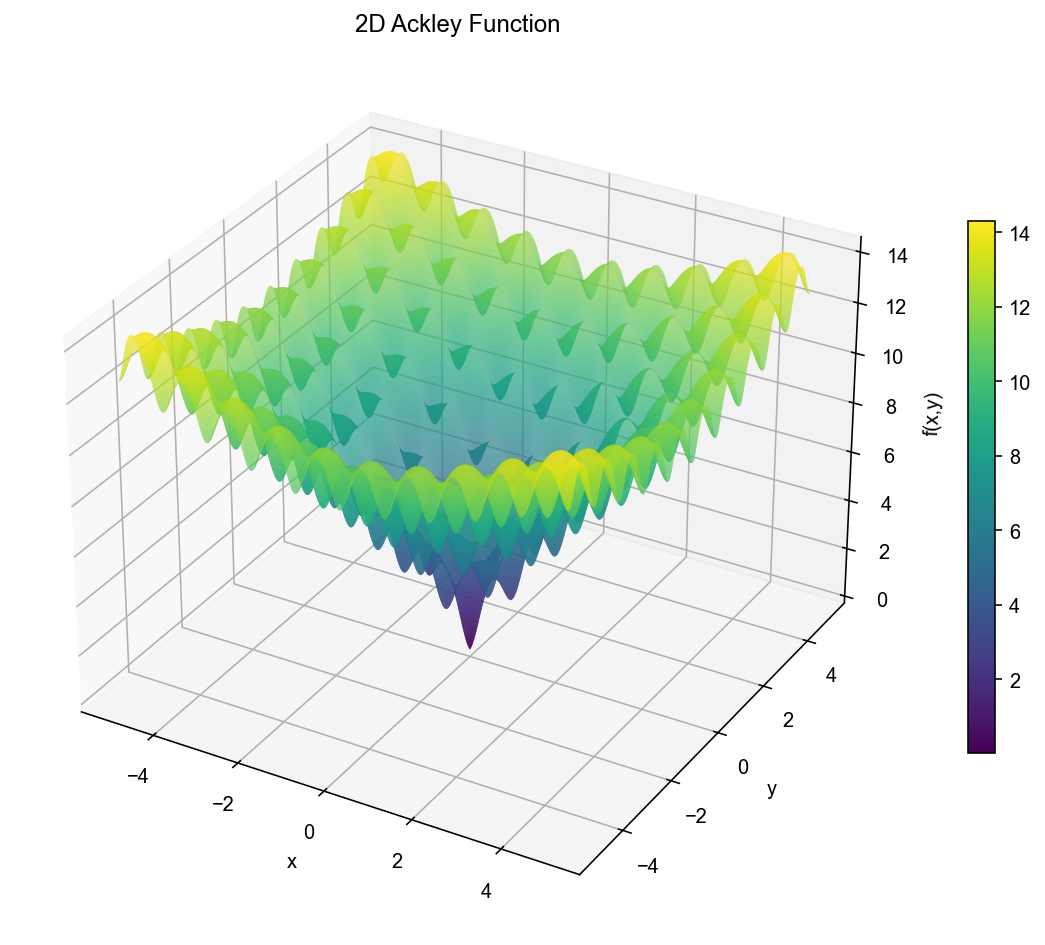}
\includegraphics[width=1.9in, height=1.4in]{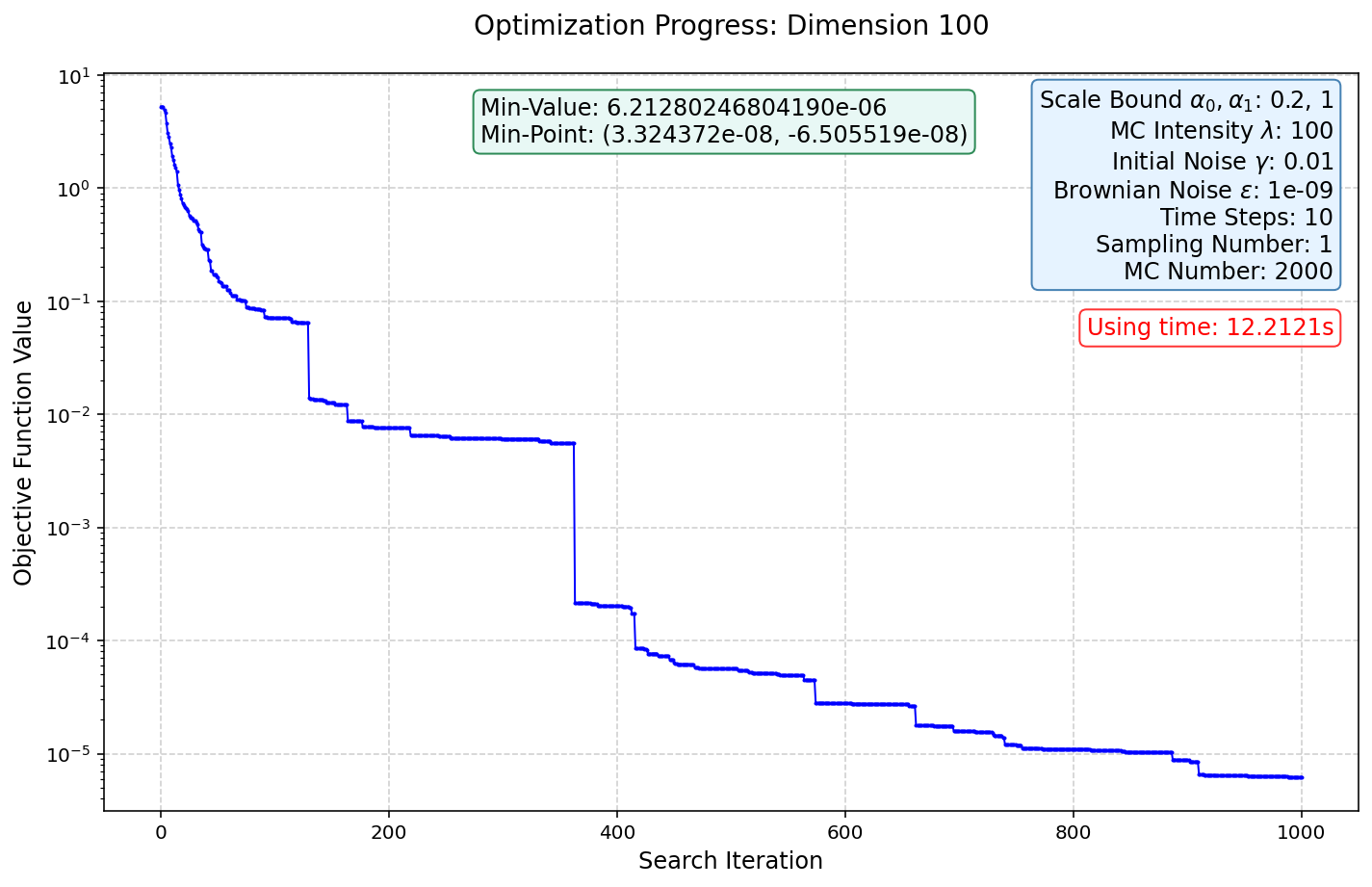}
\includegraphics[width=1.9in, height=1.4in]{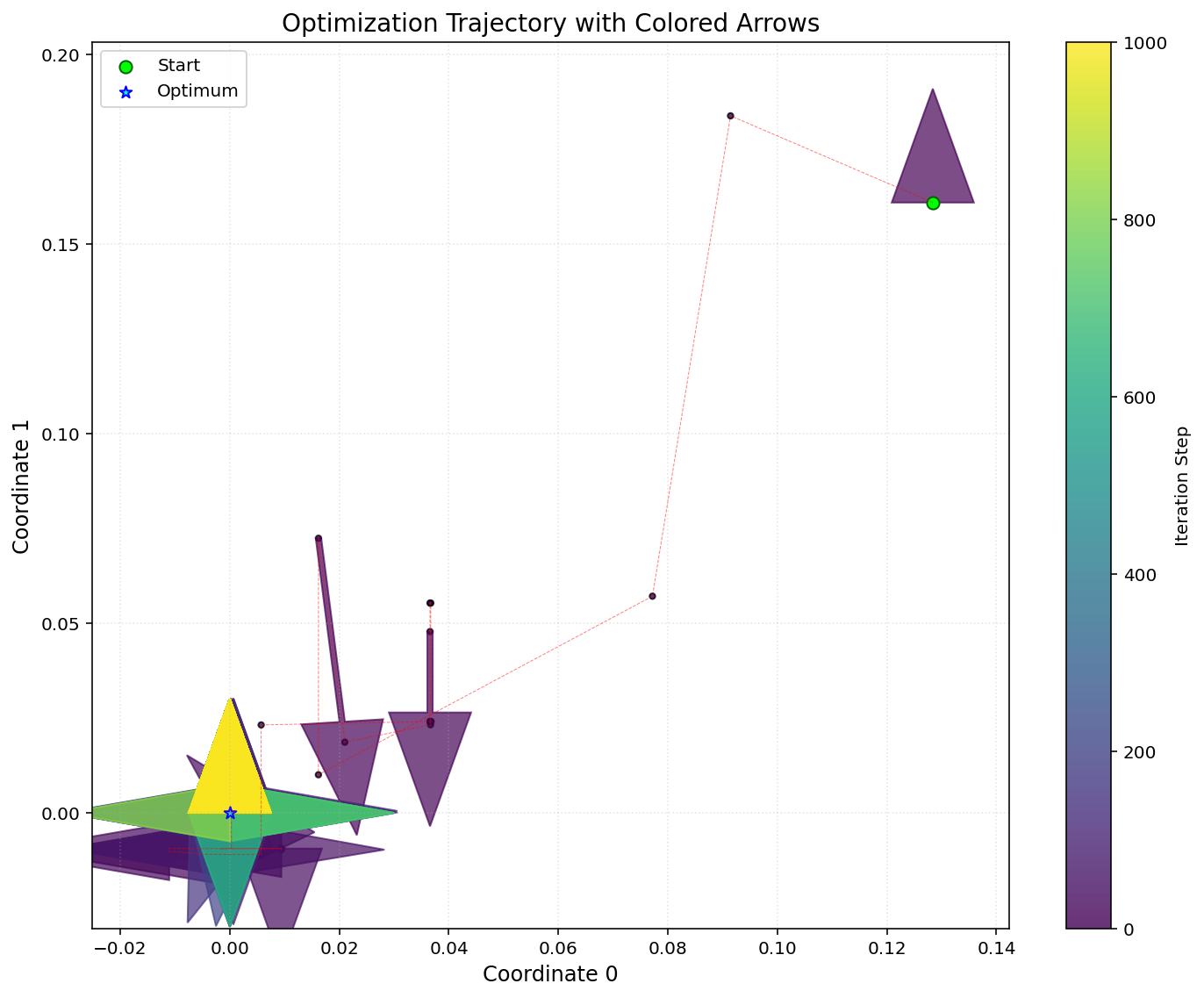}
\captionof{figure}{Ackley Function}
\end{minipage}
\end{center}
\ex

\bx\rm
\texttt{Griewank Function}: differentiable, non-separable, scalable, multimodal.

This function has many regularly spaced local minima, providing a useful test for evaluating global exploration.
The $d$-dimensional Griewank function is defined as:
$$
U(x) = 1 + \frac{1}{4000} \sum_{i=1}^{d} x_i^2 - \prod_{i=1}^{d} \cos\left(\frac{x_i}{\sqrt{i}}\right),
$$
where  the global minimum $U(x_*) = 0$ is achieved at $x_* = (0, \dots, 0)$.
Numerical results are shown below.
\begin{center}
\begin{minipage}{\textwidth}
\centering
\includegraphics[width=1.9in, height=1.4in]{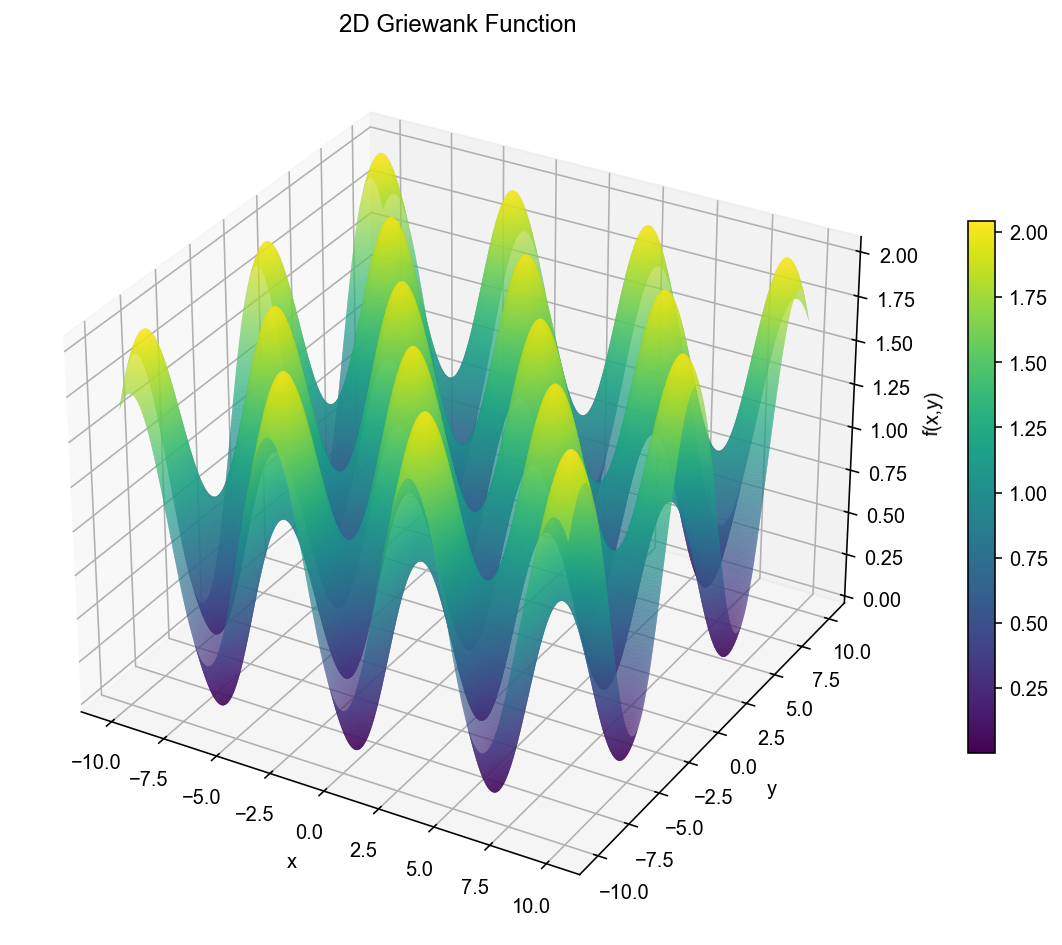}
\includegraphics[width=1.9in, height=1.4in]{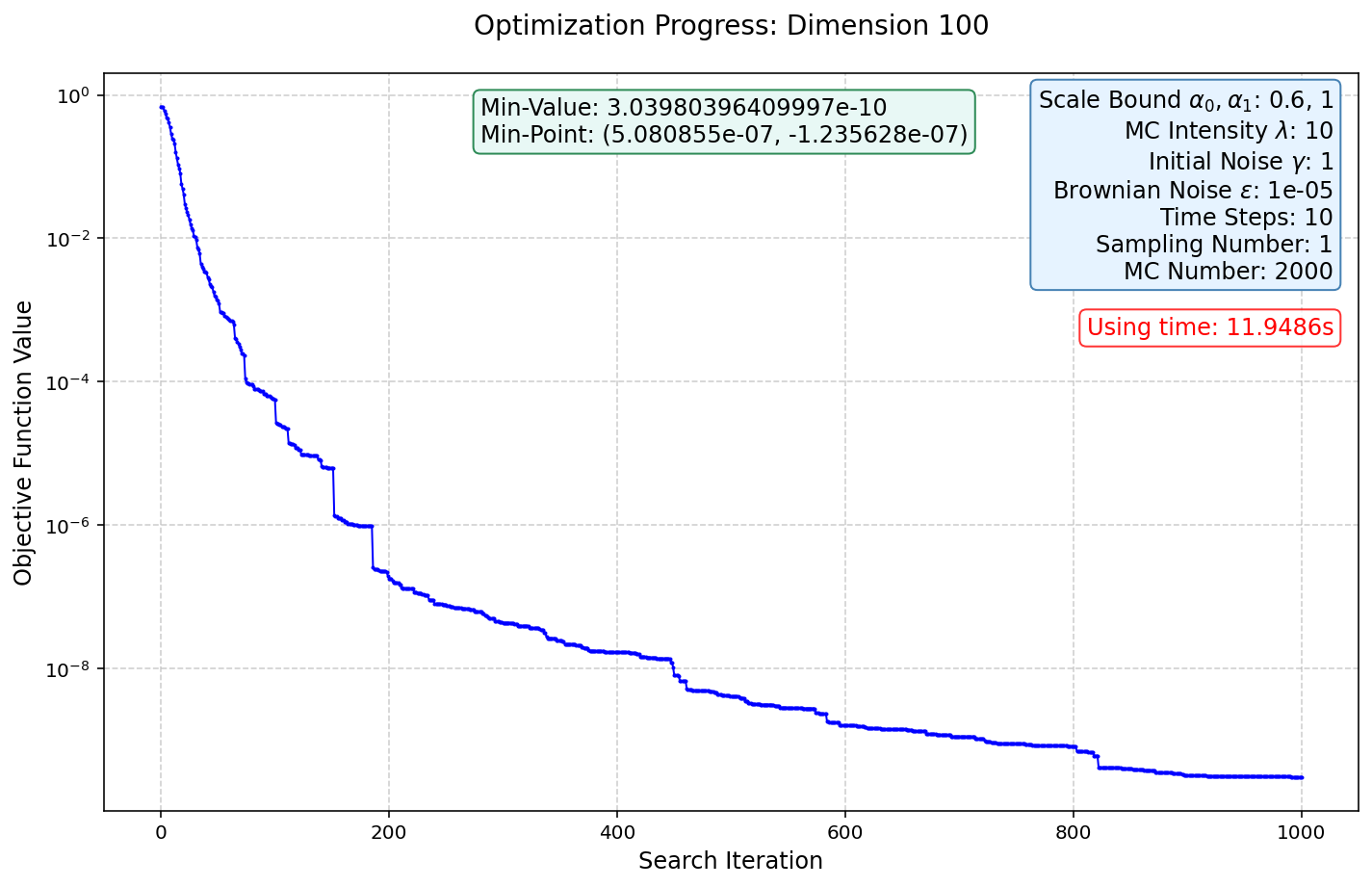}
\includegraphics[width=1.9in, height=1.4in]{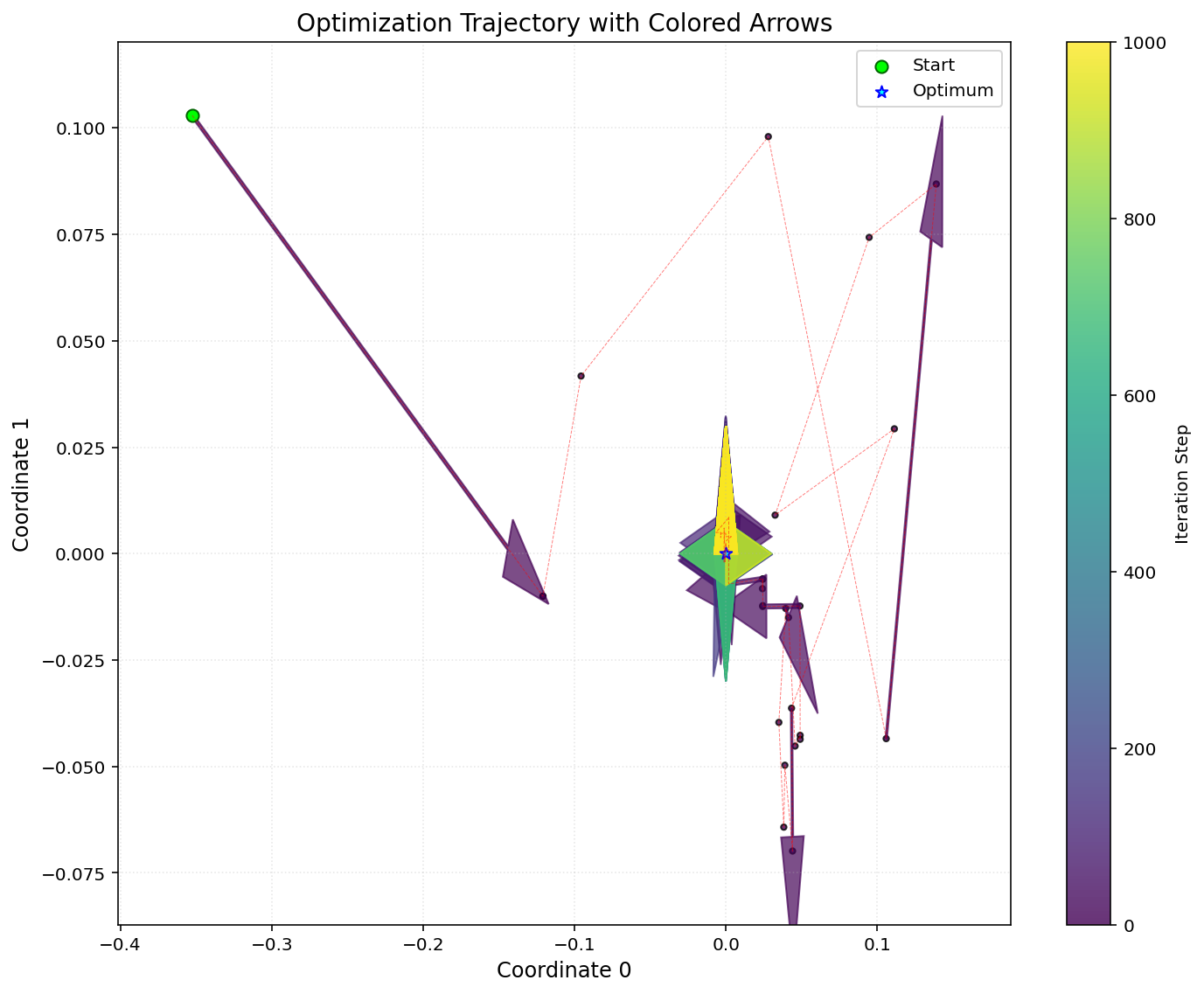}
\captionof{figure}{Griewank Function}
\end{minipage}
\end{center}
\ex

\bx\rm
\texttt{Rastrigin Function}: differentiable, separable, scalable, multimodal.

The Rastrigin function is composed of a quadratic term and cosine modulation, yielding a landscape full of regular local minima.
The $d$-dimensional Rastrigin function is defined as:
$$
U(x) = 10*d + \sum_{i=1}^d \left[ x_i^2 - 10\cos(2\pi x_i) \right], \quad x\in [-5.12, 5.12]^d,
$$
where the global minimum $U(x_*) = 0$ is achieved at $x_* = (0, \dots, 0)$.
Below are the results of our numerical experiments. 
\begin{center}
\begin{minipage}{\textwidth}
\centering
\includegraphics[width=1.9in, height=1.4in]{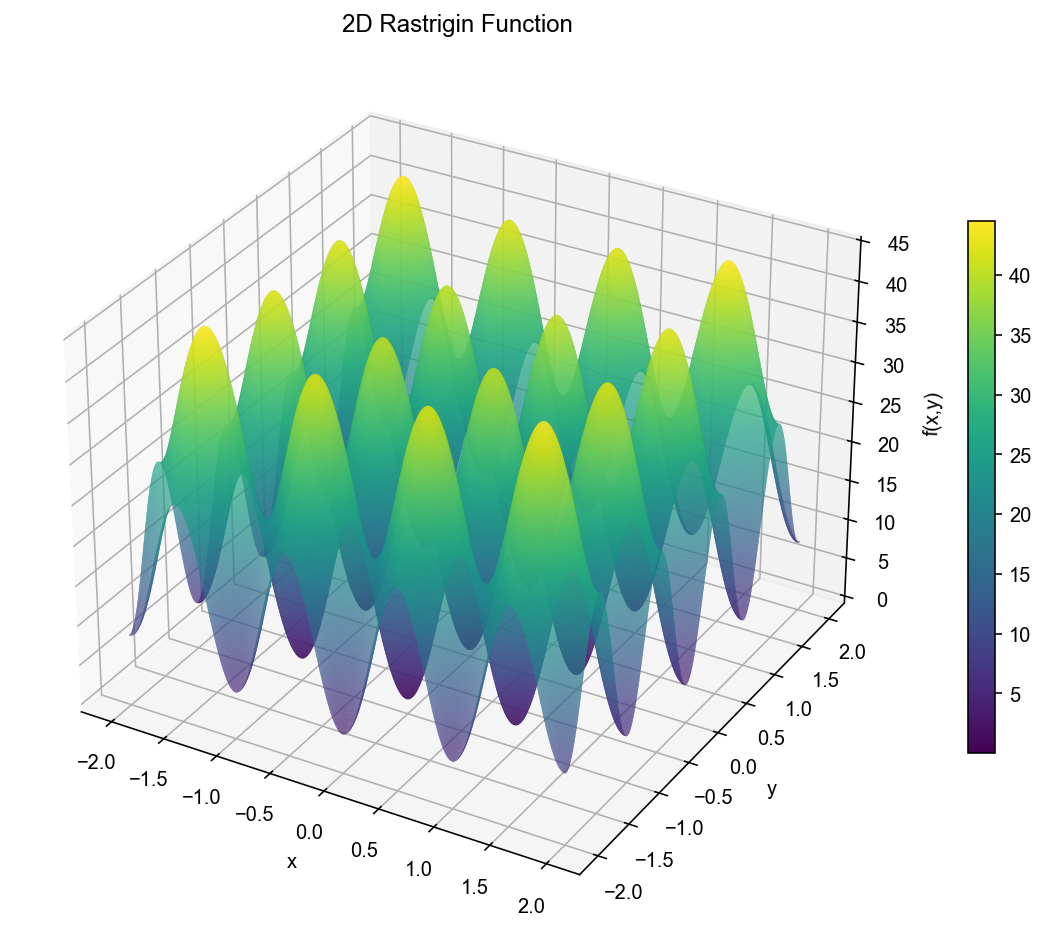}
\includegraphics[width=1.9in, height=1.4in]{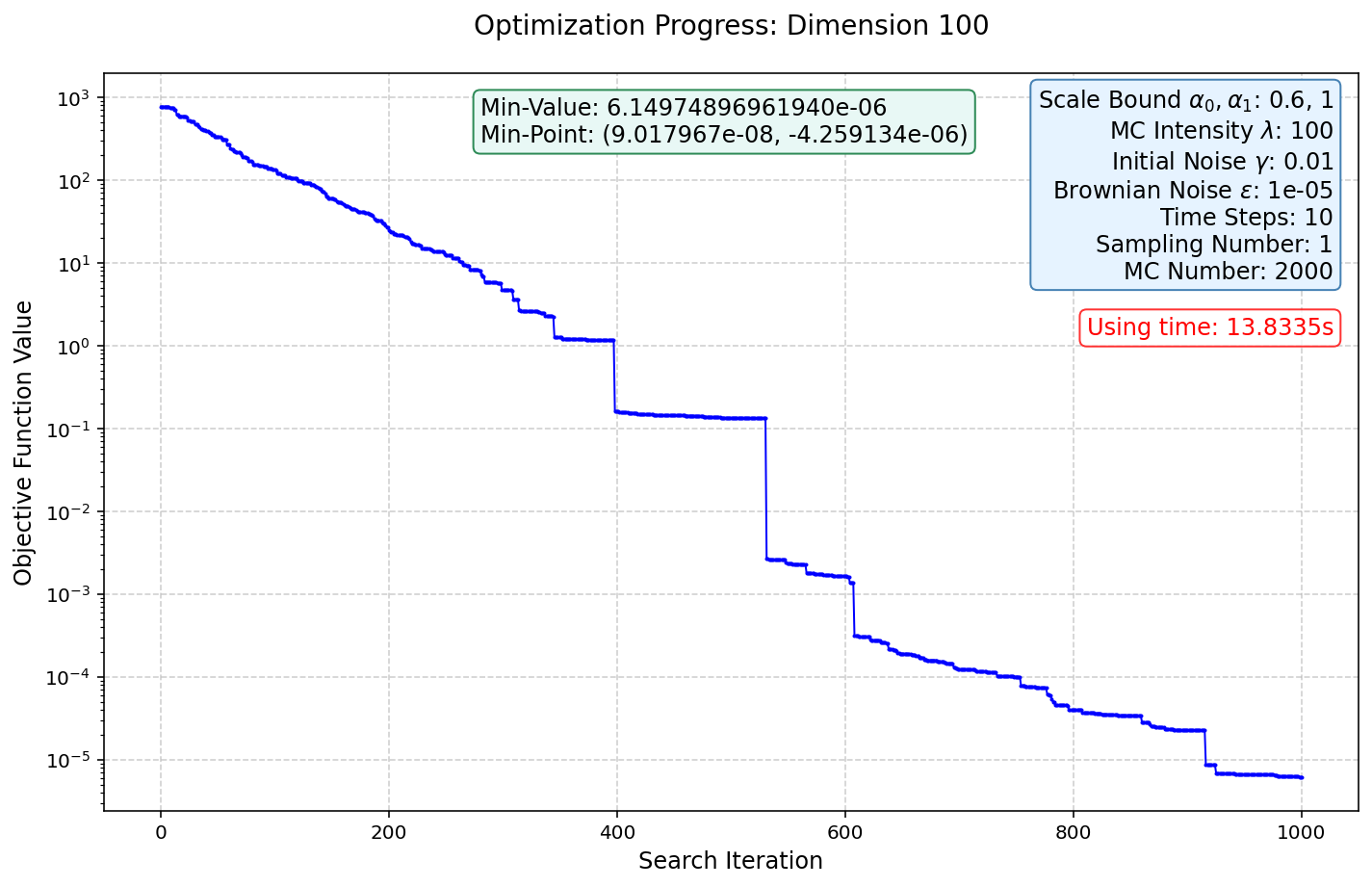}
\includegraphics[width=1.9in, height=1.4in]{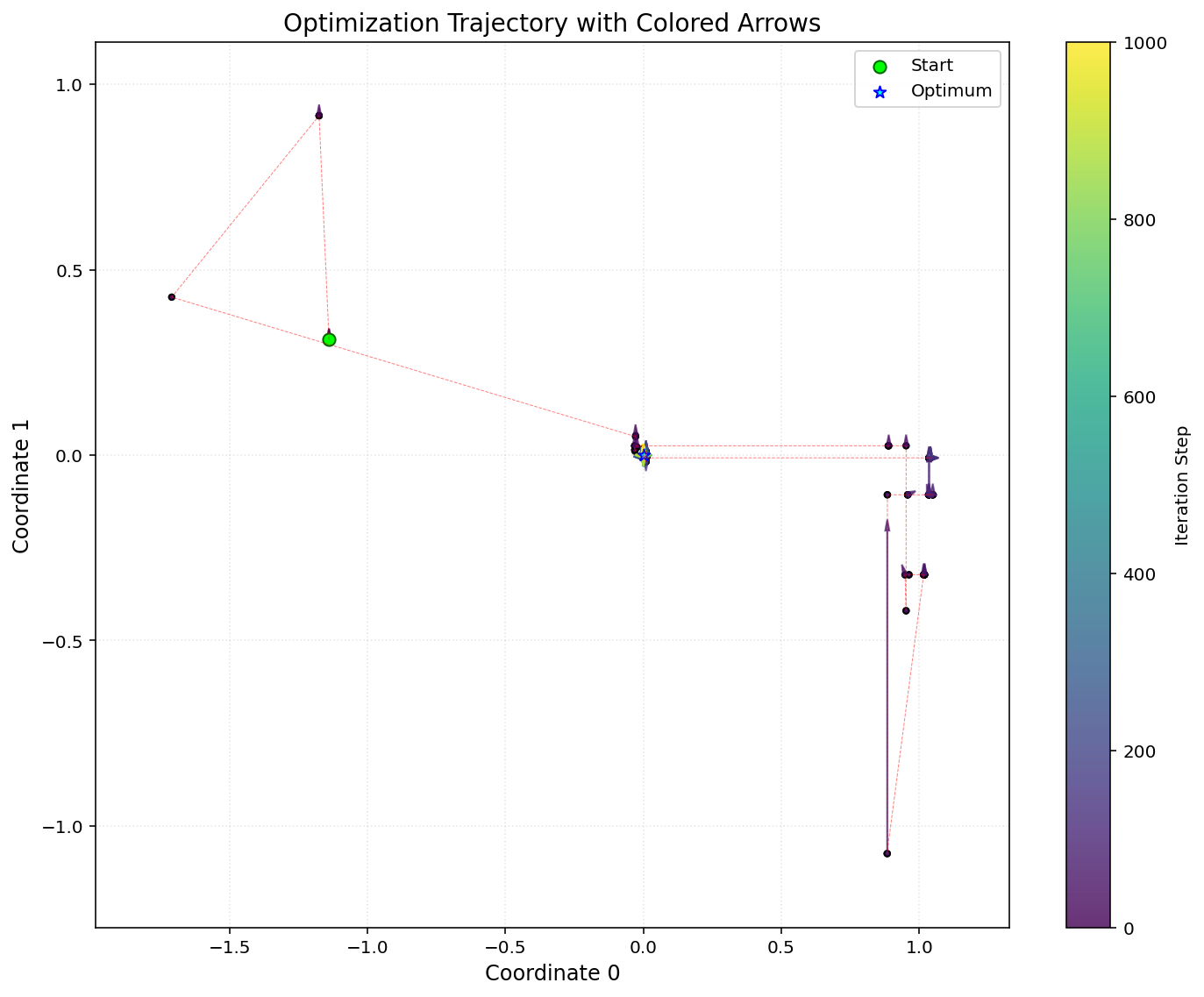}
\captionof{figure}{Rastrigin Function}
\end{minipage}
\end{center}
\ex

\bx\rm
\texttt{L\'evy Function}: differentiable, separable, scalable, multimodal.

L\'evy function is characterized by a rugged, complex landscape with numerous sharp local minima.
The $d$-dimensional L\'evy function is defined as:
$$
U(x) = \sin^2(\pi w_1) + \sum_{i=1}^{d-1} (w_i - 1)^2 \left[1 + 10\sin^2(\pi w_i + 1)\right] + (w_d - 1)^2 \left[1 + \sin^2(2\pi w_d)\right],
$$
where
$$
w_i = 1 + \frac{x_i - 1}{4}, \quad x\in [-10, 10]^d.
$$
The global minimum $U(x_*) = 0$ is achieved at $x_* = (0, \dots, 0)$.
Below are the results of our numerical experiments. 
\begin{center}
\begin{minipage}{\textwidth}
\centering
\includegraphics[width=1.9in, height=1.4in]{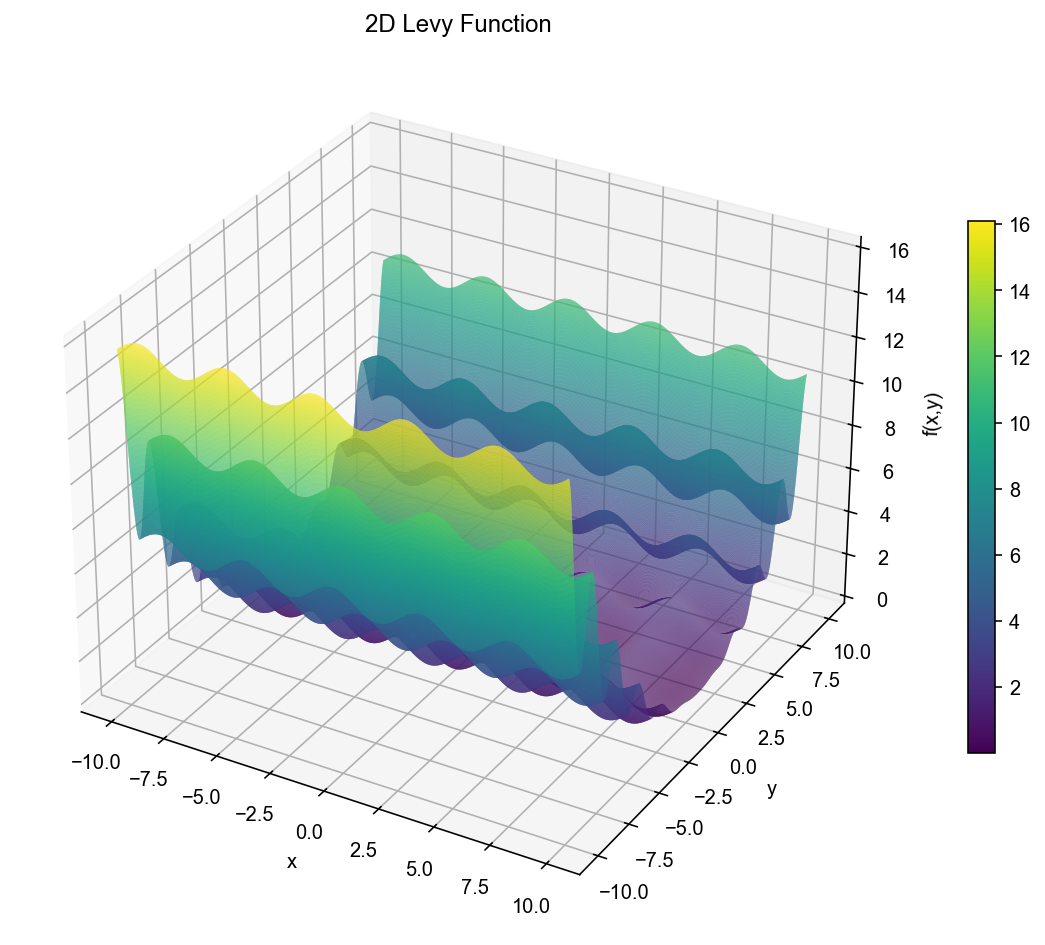}
\includegraphics[width=1.9in, height=1.4in]{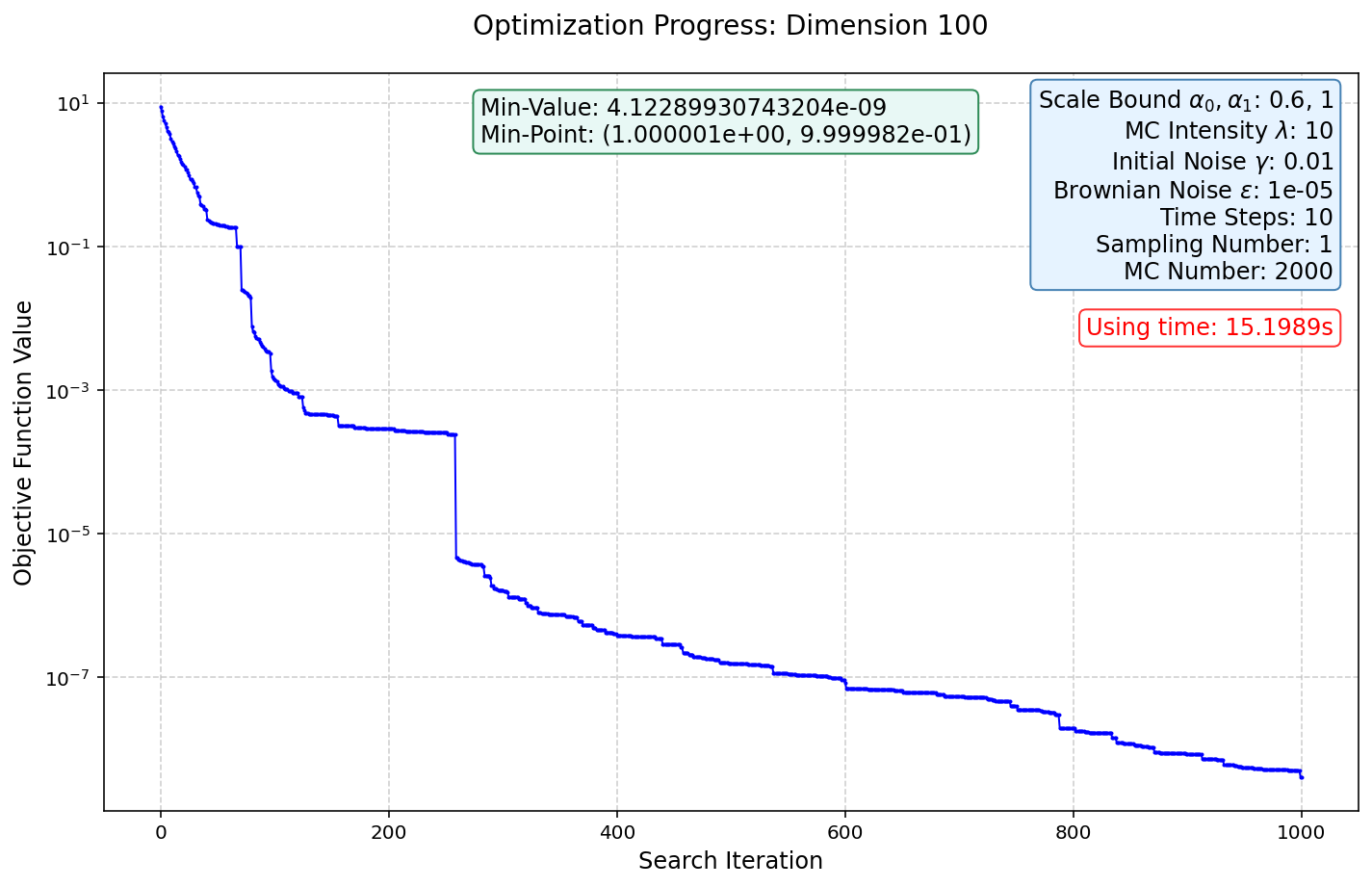}
\includegraphics[width=1.9in, height=1.4in]{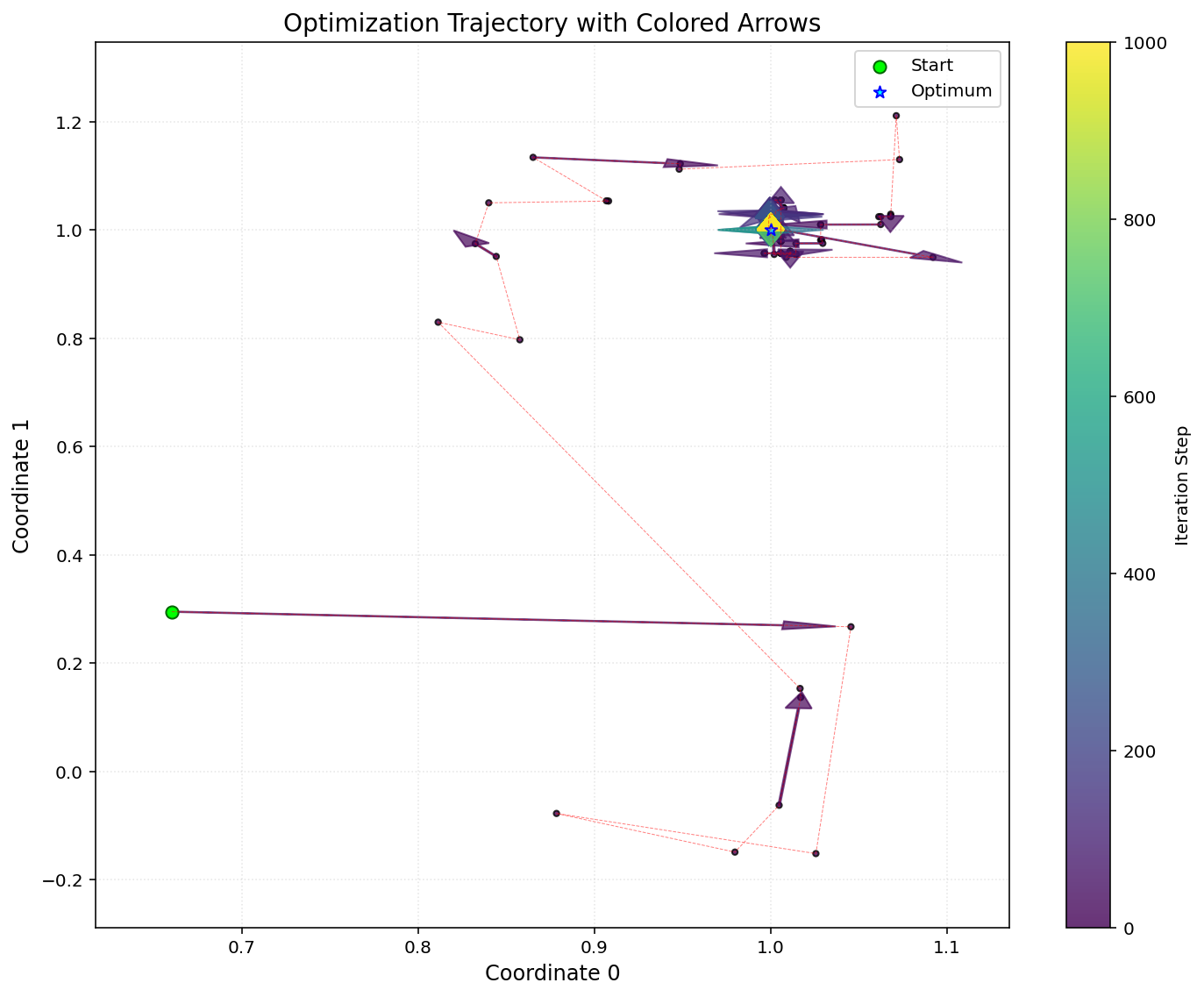}
\captionof{figure}{L\'evy Function}
\end{minipage}
\end{center}
\ex

\bx\rm
\texttt{Weierstrass Function}: non-differentiable, separable, scalable, multimodal.

This function is continuous everywhere but differentiable nowhere, representing a hard test case for classical optimizers, which is defined by
$$
U(x) = \sum_{i=1}^d\sum_{k=0}^{k_{\rm max}} a^k\cos(b^k\pi x_i) + \frac{d}{1-a},
$$
where $a \in (0,1)$, $b$ is an odd integer, and $ab > 1 + \frac{3\pi}{2}$. 

In our setup, $a = 0.5$, $b = 13$. In this case, its global minimum is approximately:
$$
\min(U) \approx d \cdot \frac{a^{k_{\rm max}+1}}{1 - a}.
$$
Below are the results of our numerical experiments. 
\begin{center}
\begin{minipage}{\textwidth}
\centering
\includegraphics[width=1.9in, height=1.4in]{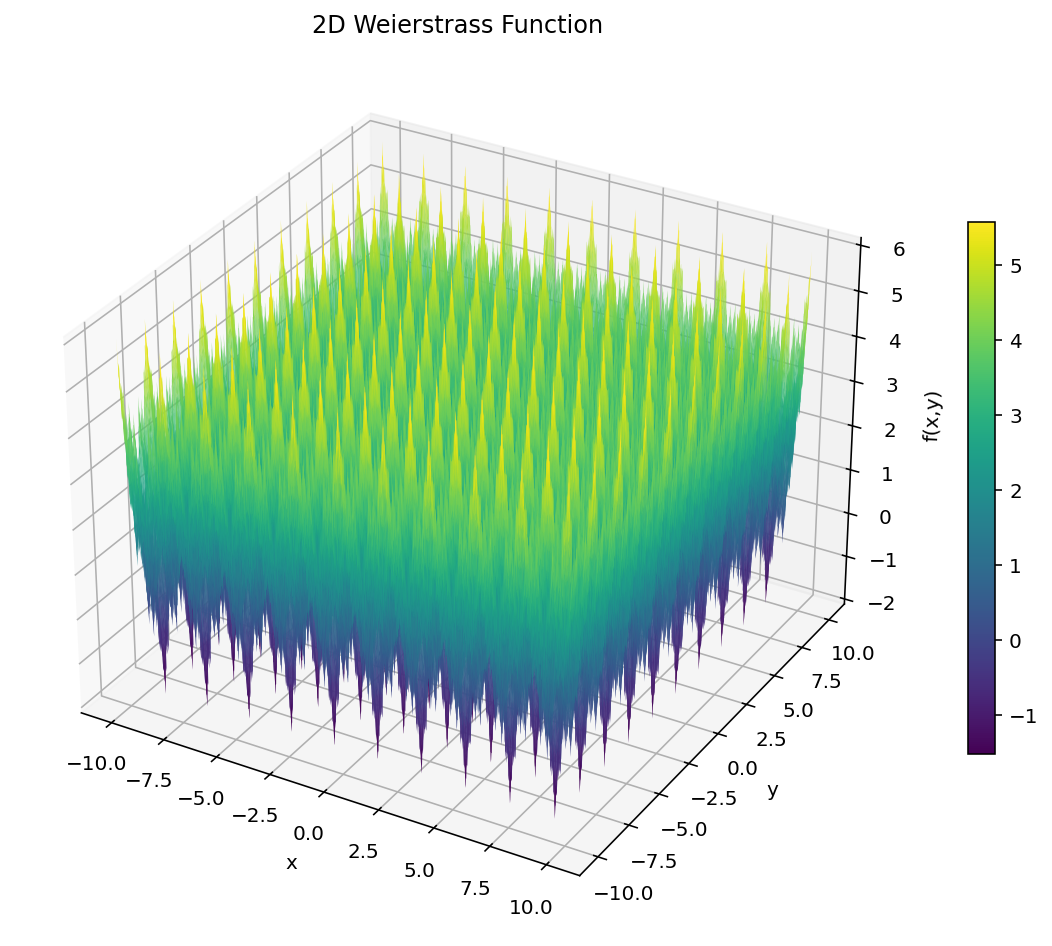}
\includegraphics[width=1.9in, height=1.4in]{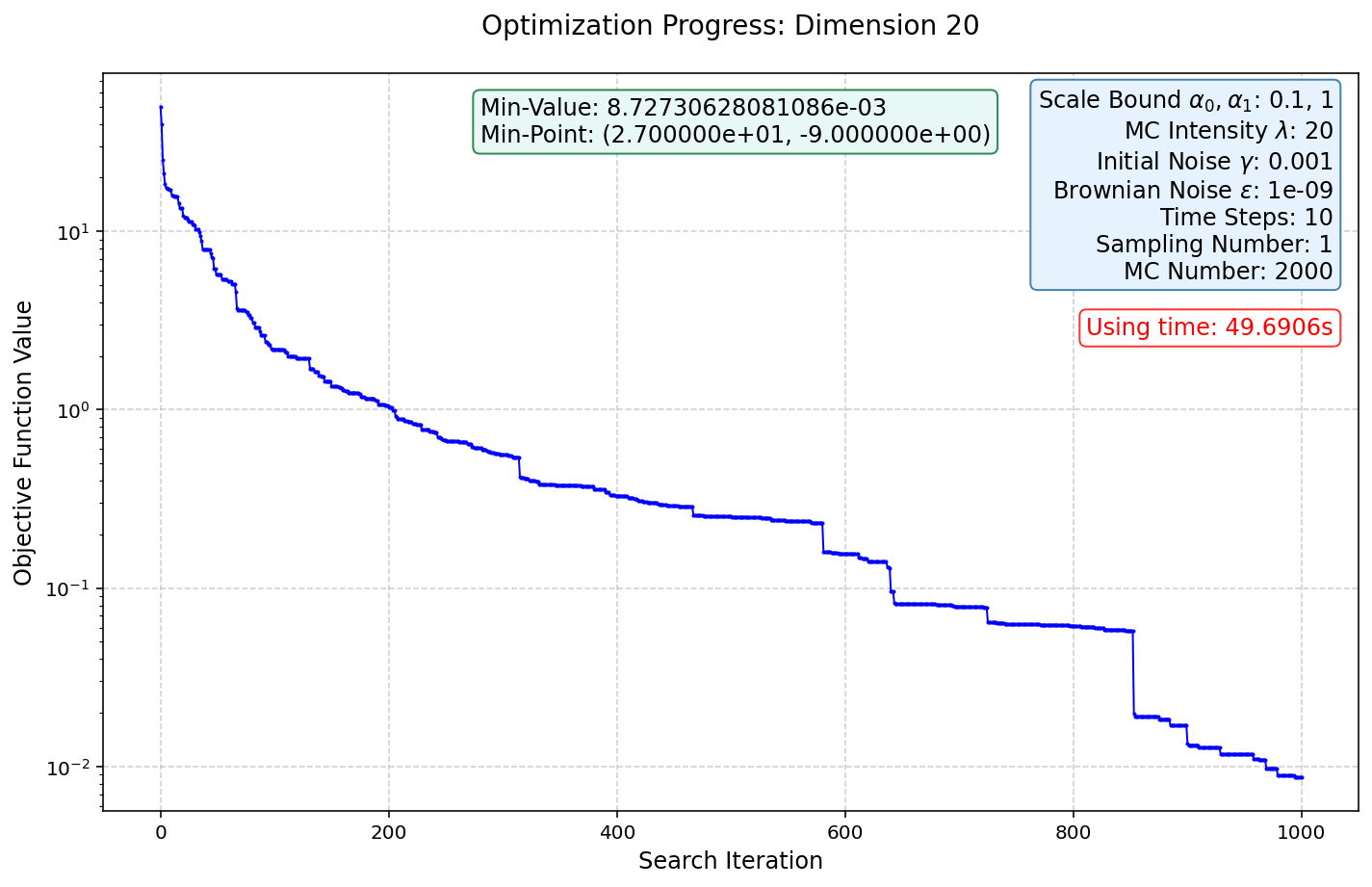}
\includegraphics[width=1.9in, height=1.4in]{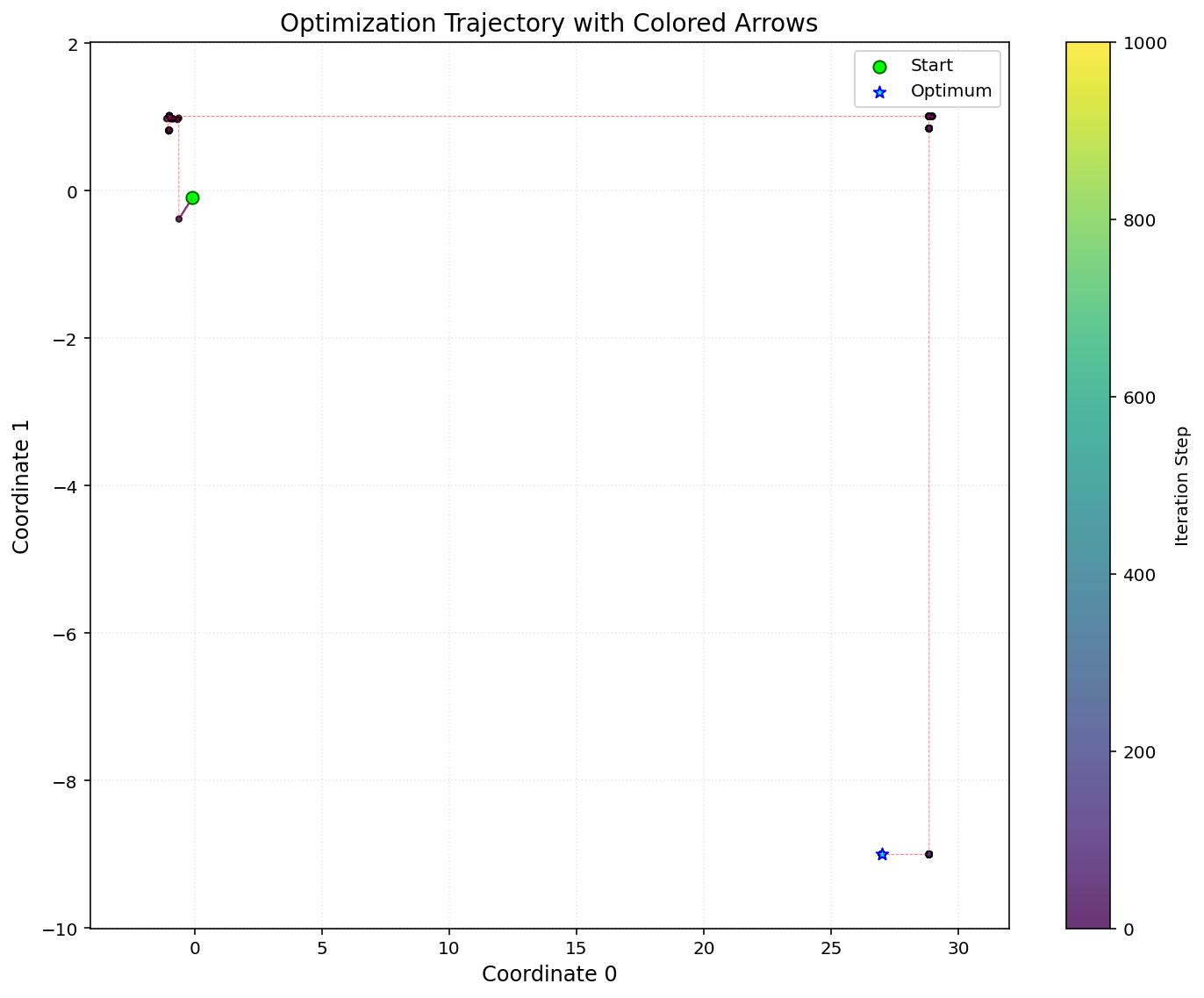}
\captionof{figure}{Weierstrass Function}
\end{minipage}
\end{center}
\ex

\bx\rm
\texttt{Discontinuous Function}: discontinuous, separable, scalable, multimodal.

Consider the following discontinuous function
$$
U(x) = \sum_{i=1}^d\lfloor x_i + 0.5 \rfloor^2,
$$
with global minimum $\min(U) = 0.$
This function features integer-valued discontinuities and abrupt jumps, making it particularly suitable for evaluating non-gradient-based optimization algorithms. The landscape contains numerous local minima. The numerical results are presented below.
\begin{center}
\begin{minipage}{\textwidth}
\centering
\includegraphics[width=1.9in, height=1.4in]{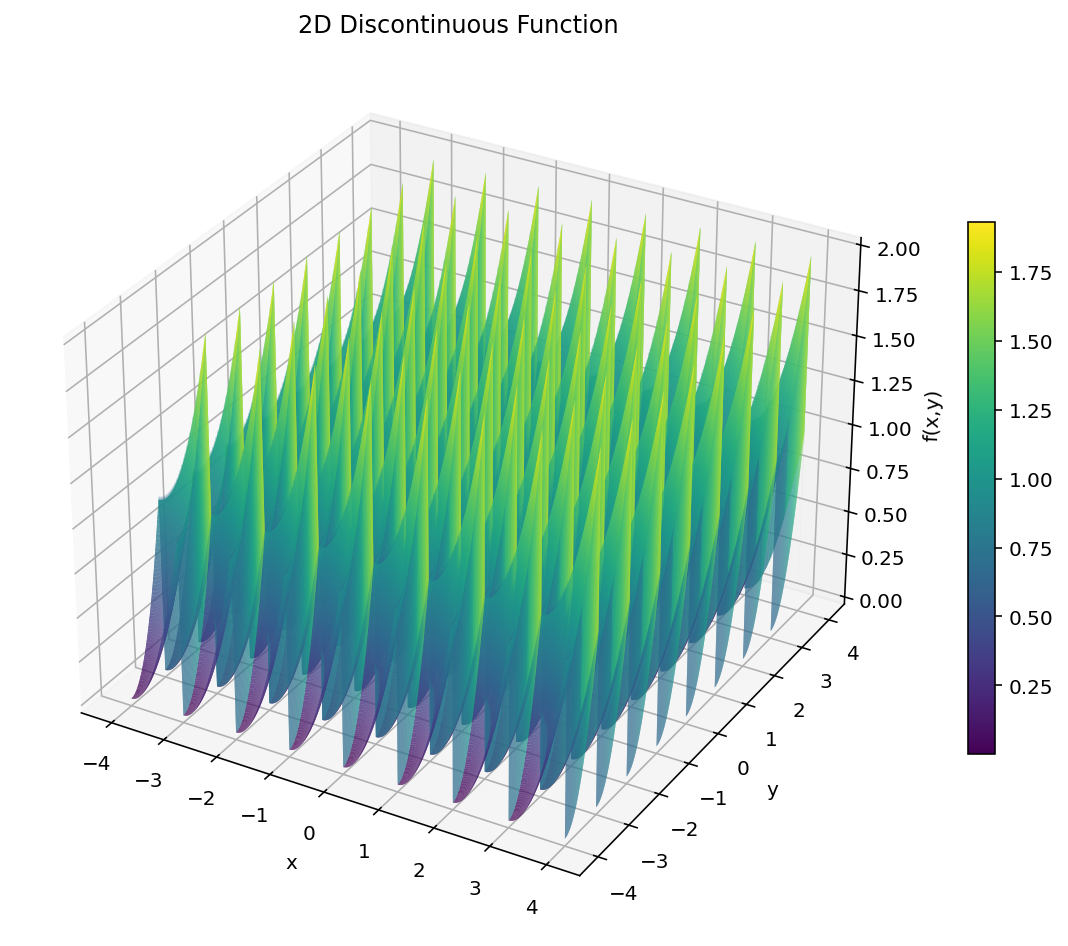}
\includegraphics[width=1.9in, height=1.4in]{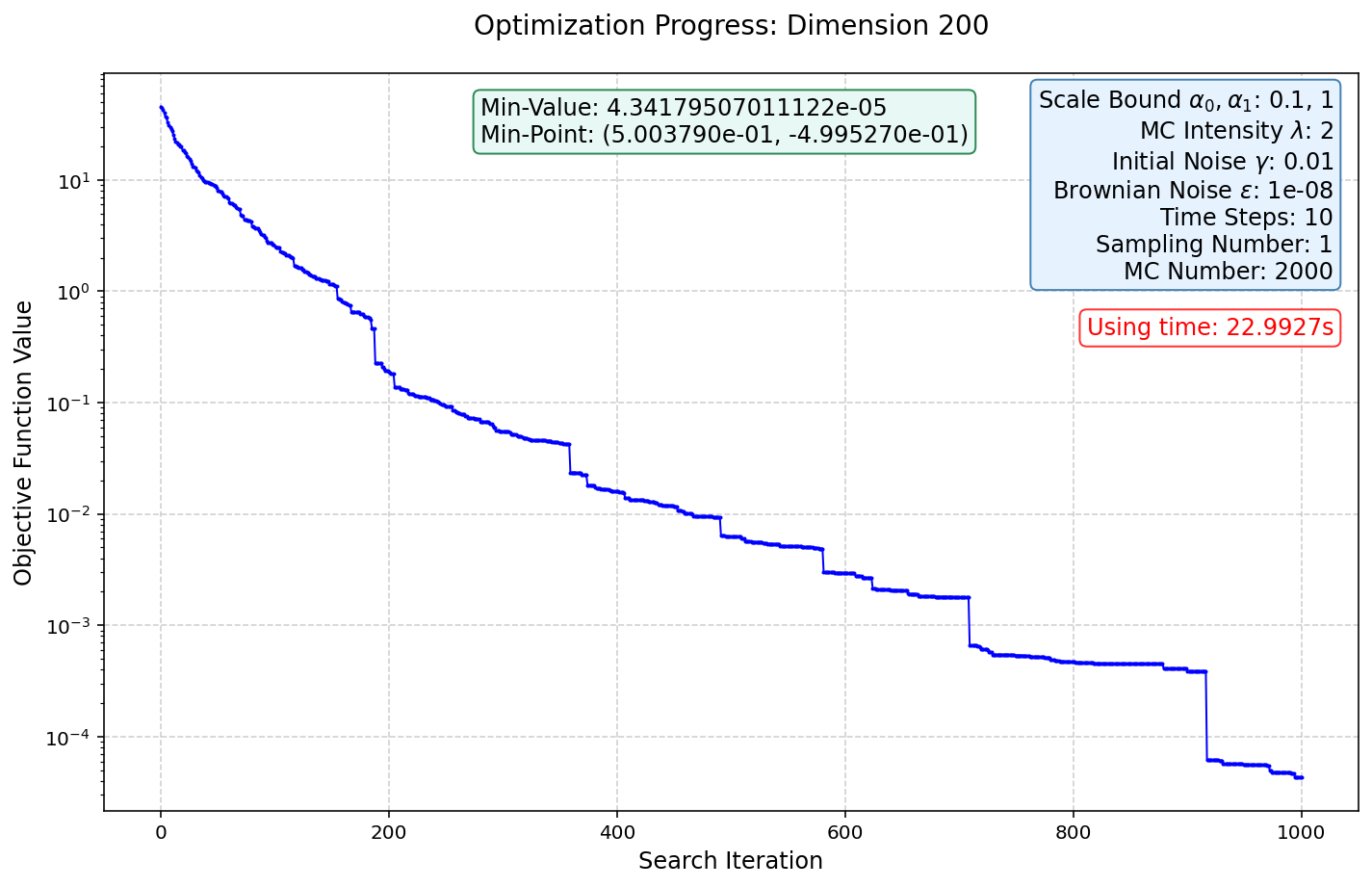}
\includegraphics[width=1.9in, height=1.4in]{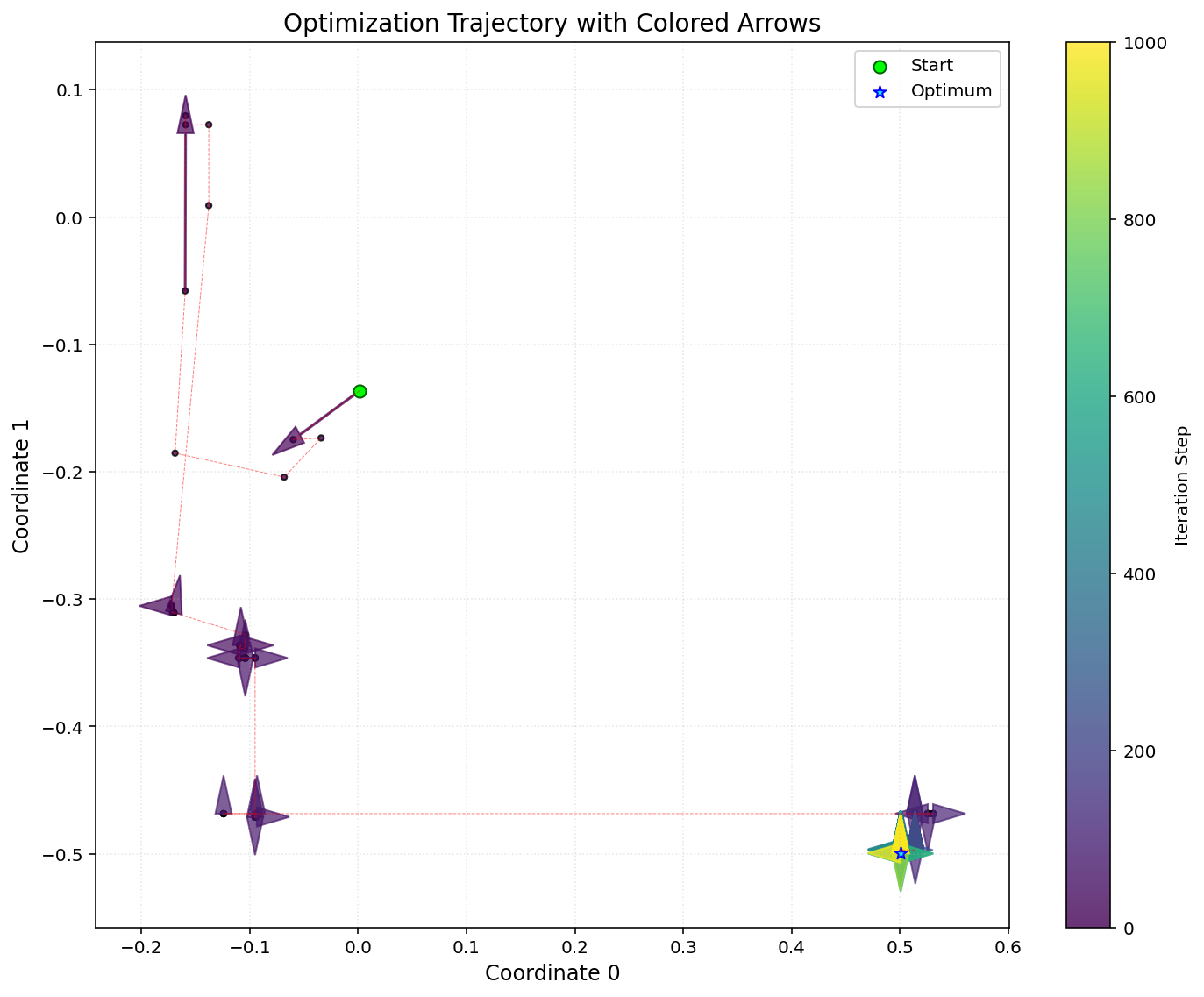}
\captionof{figure}{Discontinuous Function}
\end{minipage}
\end{center}
\ex

\bx\rm
\texttt{Artificial Function}: non-convex, separable, multimodal.

Consider the following artificially constructed benchmark function:
$$
U(x) = 10^5 \cdot \sum_{i=1}^d \sqrt{\left| \sin\left( |x_i - 0.1|^4 \right) \right|}, \quad x\in \mathbb{R}^d,
$$
with global minimum $\min_{x} U(x) = 0$ achieved when $x_i = 0.1$ for all $i = 1,\dots,d$.

This synthetic function was specifically designed to present significant challenges for conventional optimization algorithms. Key characteristics include:
\begin{itemize}
    \item High-frequency oscillations induced by the quartic exponent and absolute value operations.
    \item Non-differentiable points at $x_i = 0.1$ due to the absolute value term.
    \item Extreme scaling ($10^5$ multiplier) creating numerical stability issues.
\end{itemize}

As demonstrated in our experiments, gradient-based methods like L-BFGS-B consistently fail when initialized beyond an immediate neighborhood of the optimum. For the high-dimensional case ($d=100$), we observe the following results:
\begin{center}
\begin{minipage}{\textwidth}
\centering
\includegraphics[width=1.9in, height=1.4in]{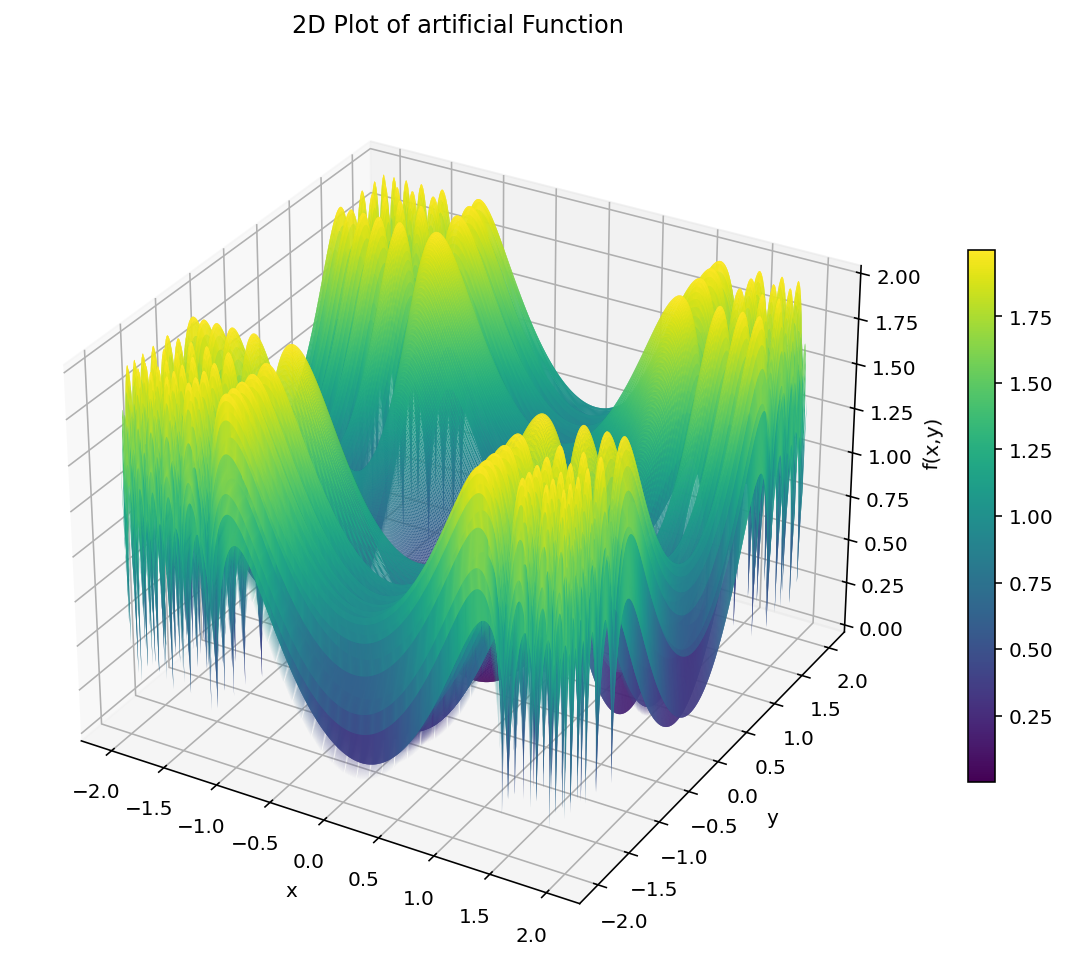}
\includegraphics[width=1.9in, height=1.4in]{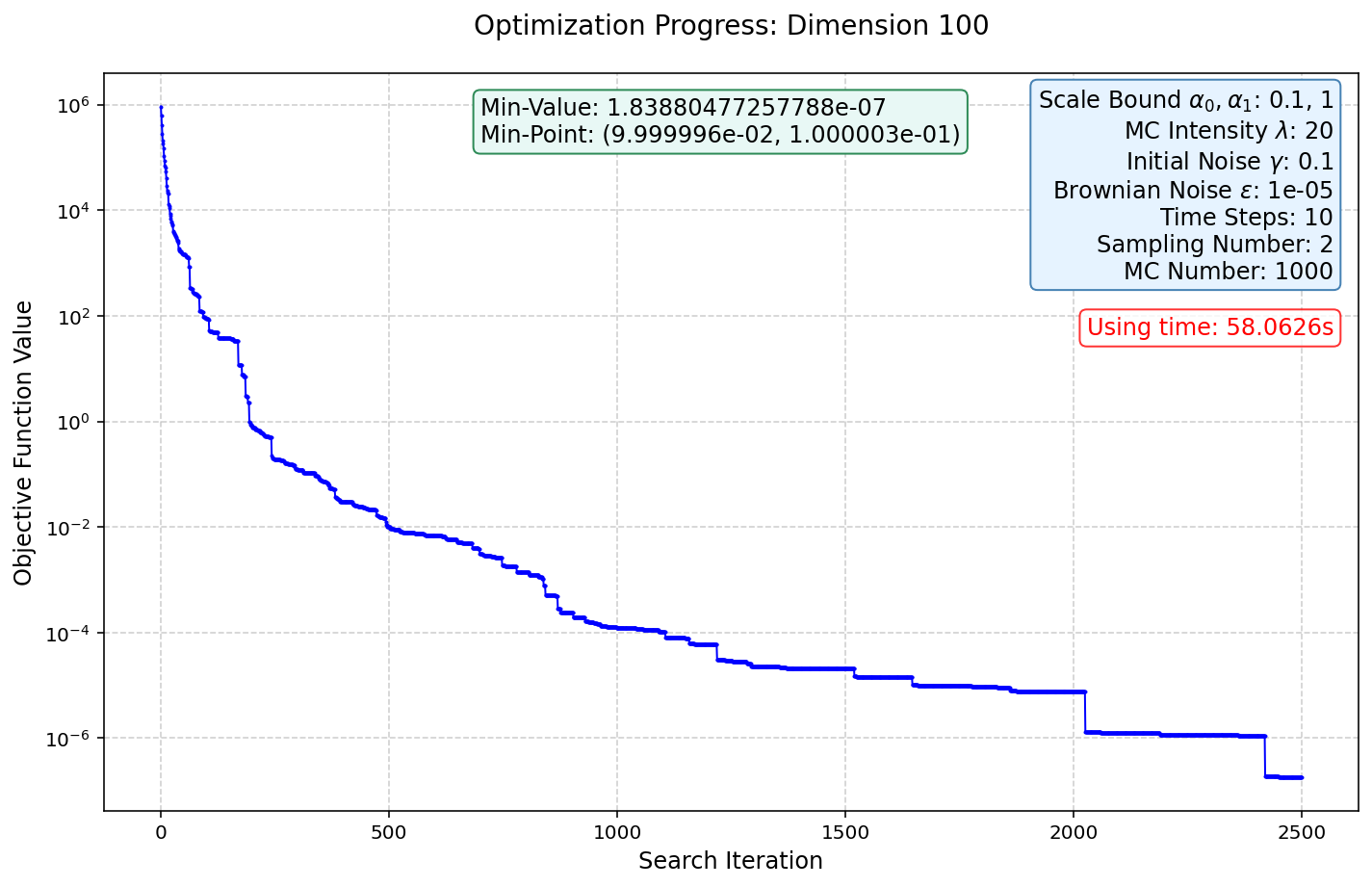}
\includegraphics[width=1.9in, height=1.4in]{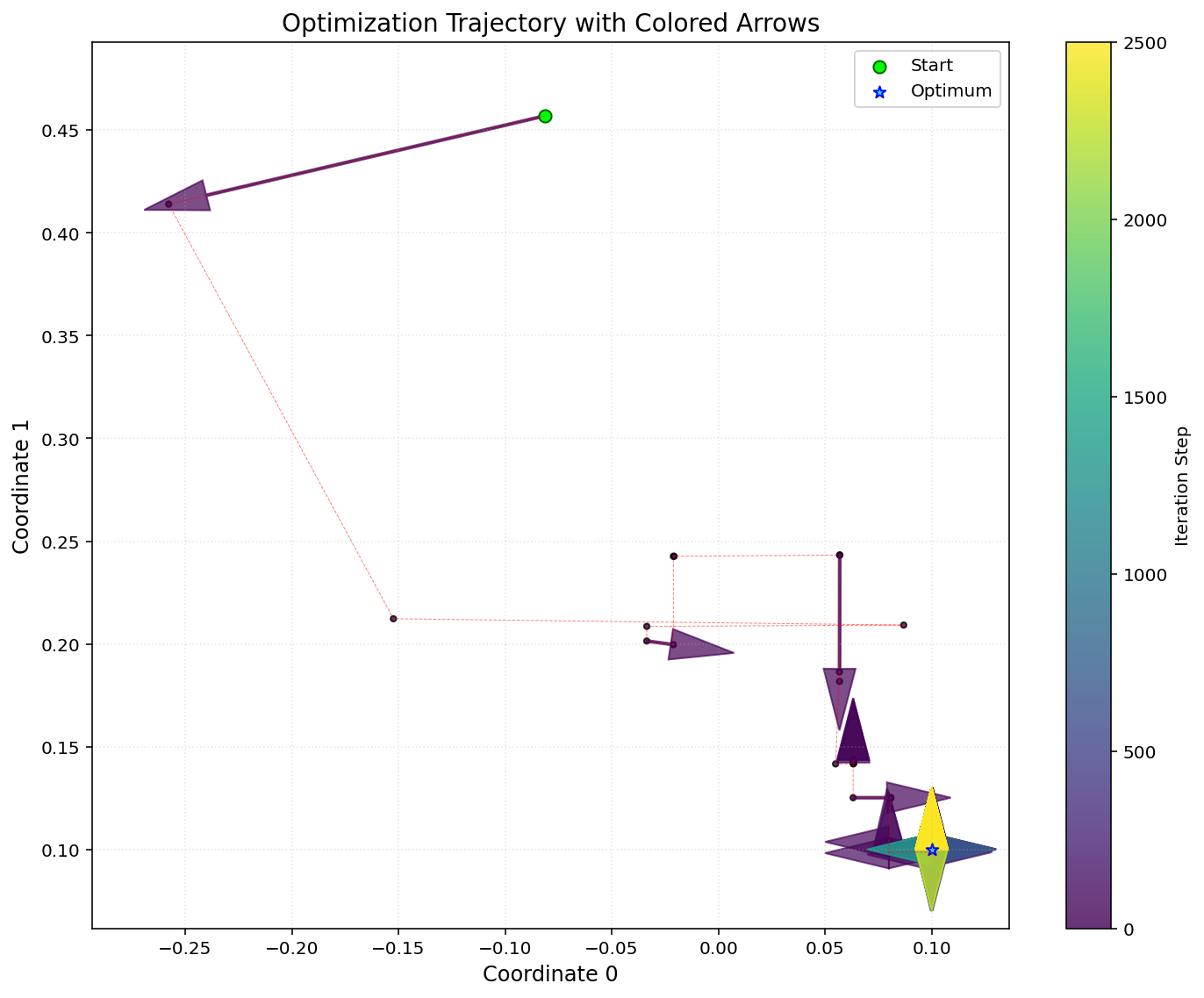}
\captionof{figure}{Artificial Function}
\end{minipage}
\end{center}
\ex

\subsection{ Comparative Analysis of Optimization Algorithms}

We present a systematic evaluation of six established optimization algorithms against our proposed sampling-based method. This analysis focuses on both solution quality and computational efficiency across diverse benchmark functions.
The six baseline methods considered are:

\begin{itemize}
    \item \textbf{Particle Swarm Optimization (PSO)} \cite{kennedy1995particle}: PSO simulates a swarm of particles (candidate solutions) that traverse the search space, influenced by both their personal best-known positions and the global best-known position. The collective dynamics guide the swarm toward promising regions.
    
    \item \textbf{Differential Evolution (DE)} \cite{storn1997differential}: DE evolves a population of candidate solutions via mutation, crossover, and selection. Mutations are based on the differences between randomly selected individuals, promoting diversity and exploration.
    
    \item \textbf{BFGS Method} \cite{nocedal2006numerical}: The Broyden-Fletcher-Goldfarb-Shanno (BFGS) method approximates the inverse Hessian matrix using only first-order gradient information. It is well-suited for smooth and convex optimization problems and typically exhibits rapid convergence.
    
    \item \textbf{Simulated Annealing (SA)} \cite{kirkpatrick1983optimization}: SA explores the solution space by probabilistically accepting worse solutions in the early stages to avoid local minima, with the acceptance probability decreasing over time according to a cooling schedule.
    
    \item \textbf{Stochastic Hill Climbing (SHC)} \cite{russell2010artificial}: SHC performs a local search by evaluating random neighbors of the current solution, accepting improvements. Although simple to implement, it can easily become trapped in local optima.
    
    \item \textbf{Adam Optimizer} \cite{kingma2014adam}: Adam combines momentum-based acceleration and adaptive learning rates using running averages of gradients and squared gradients. It is widely adopted in machine learning due to its effectiveness and efficiency.
\end{itemize}

Each of these algorithms serves a specific purpose: PSO and DE excel in global derivative-free optimization, BFGS is effective for smooth gradient-based problems, SA is robust on rugged landscapes, SHC offers a lightweight local search mechanism, and Adam is a go-to optimizer in deep learning contexts.

To evaluate the performance of these methods, we consider the following benchmark functions: Sphere function ($U(x) = \|x\|^2$), Schwefel function,  Ackley function, Griewank function, Rastrigin function and the above artificial function.

The experimental setup is as follows:
\begin{itemize}
    \item \textbf{Number of iterations}: 500
    \item \textbf{Problem dimension}: $d = 30$
    \item \textbf{Number of runs}: 10 independent trials per algorithm
\end{itemize}

Each comparative figure comprises two adjacent subplots:
\begin{itemize}
    \item \textbf{Left subplot} (Optimization Performance):
    \begin{itemize}
        \item Bars represent the mean of the minimum objective values obtained across 10 independent runs.
        The black line is the standard deviation.
        \item The first six blue bars correspond to the baseline optimization methods.
        \item The final green bar represents the performance of our proposed method.
    \end{itemize}
    
    \item \textbf{Right subplot} (Computational Time Cost):
    \begin{itemize}
        \item The first six red bars depict the average computational time required by the baseline algorithms.
        \item The final purple bar shows the computational time of our proposed approach.
    \end{itemize}
\end{itemize}

\begin{center}
\begin{minipage}{\textwidth}
\centering
\includegraphics[width=2.2in, height=1.2in]{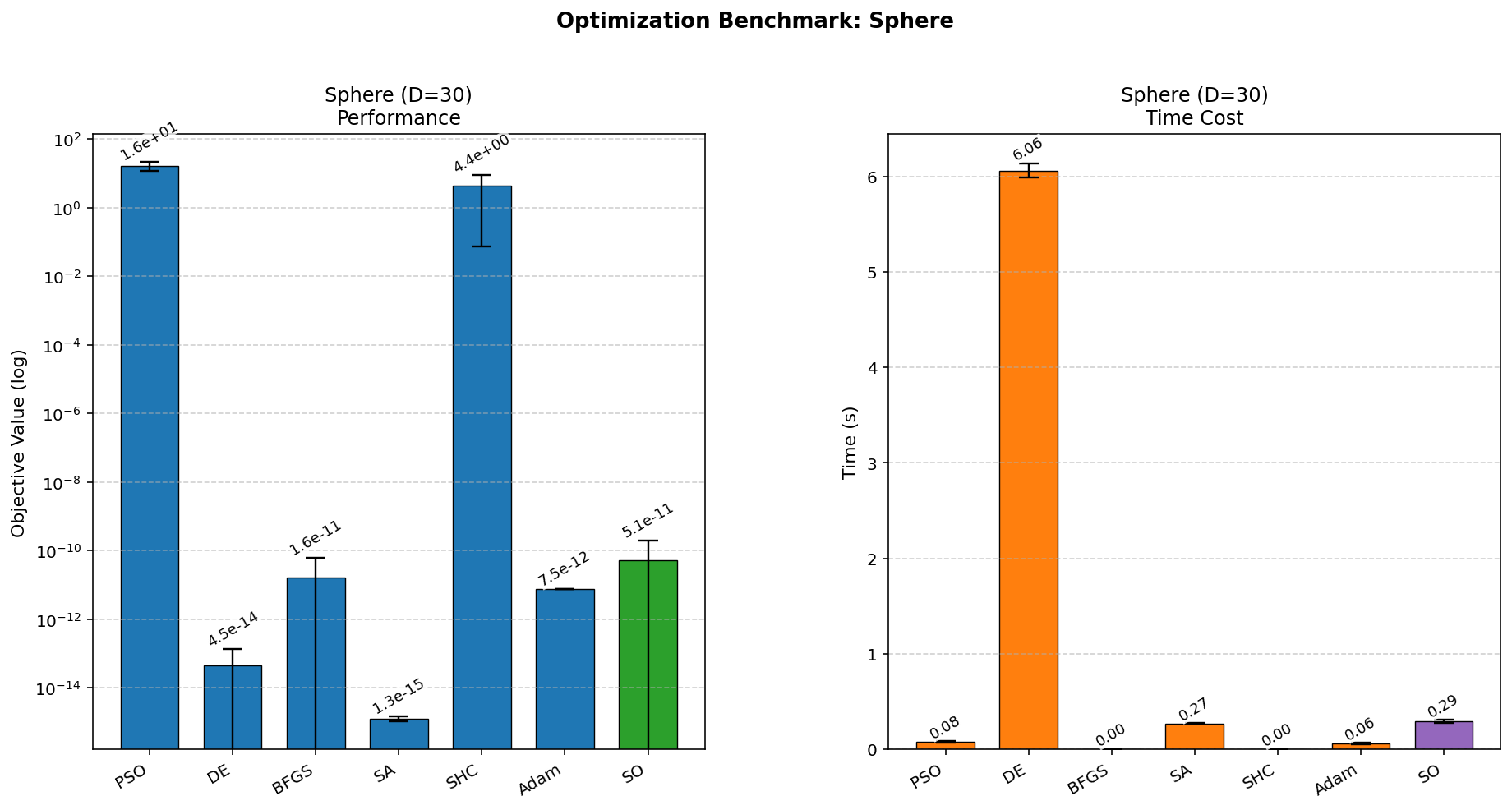}
\includegraphics[width=2.2in, height=1.2in]{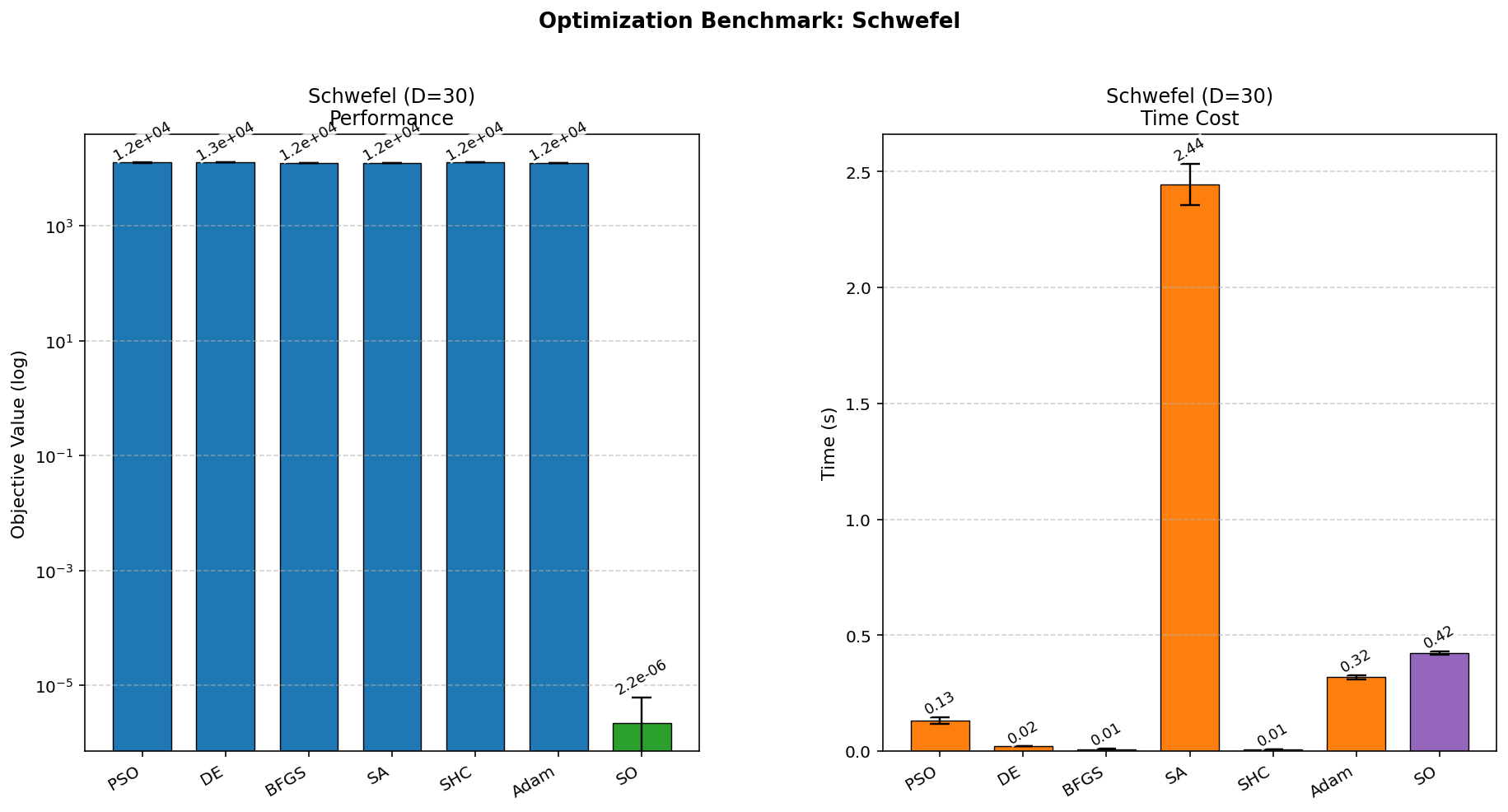}
\captionof{figure}{Sphere and Schwefel functions}
\end{minipage}
\end{center}

\begin{center}
\begin{minipage}{\textwidth}
\centering
\includegraphics[width=2.2in, height=1.2in]{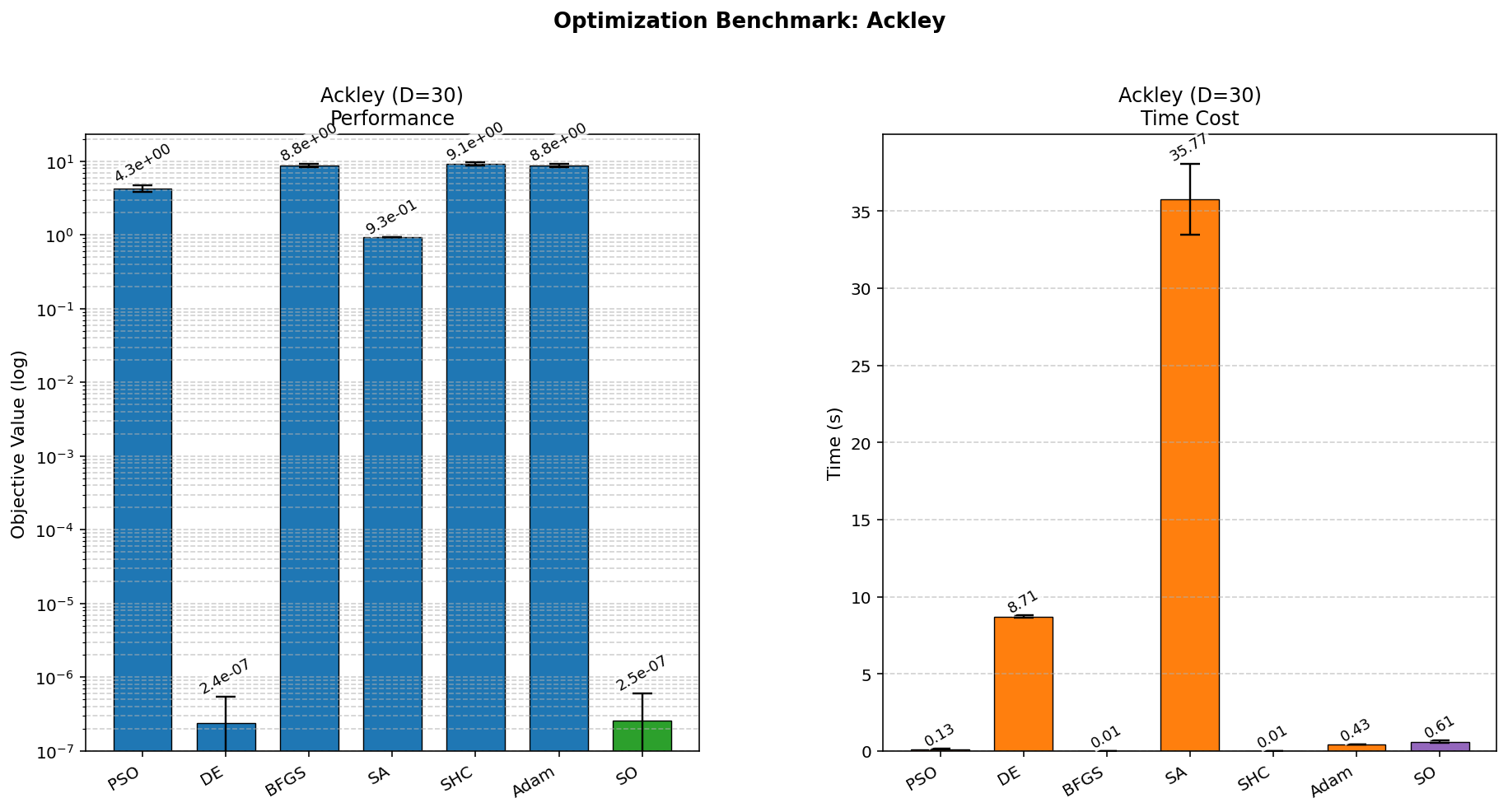}
\includegraphics[width=2.2in, height=1.2in]{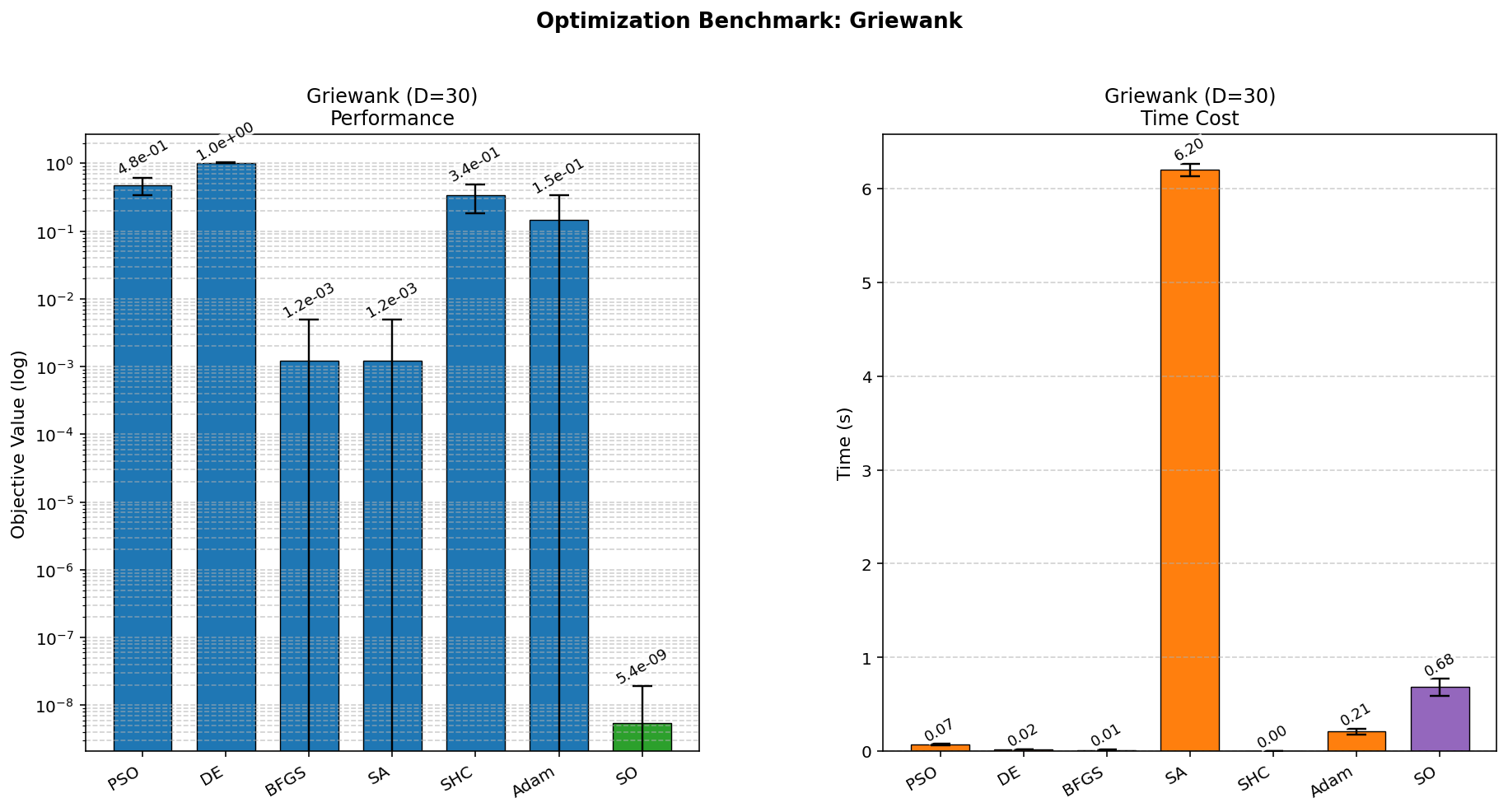}
\captionof{figure}{Ackley and Griewank functions}
\end{minipage}
\end{center}

\begin{center}
\begin{minipage}{\textwidth}
\centering
\includegraphics[width=2.2in, height=1.2in]{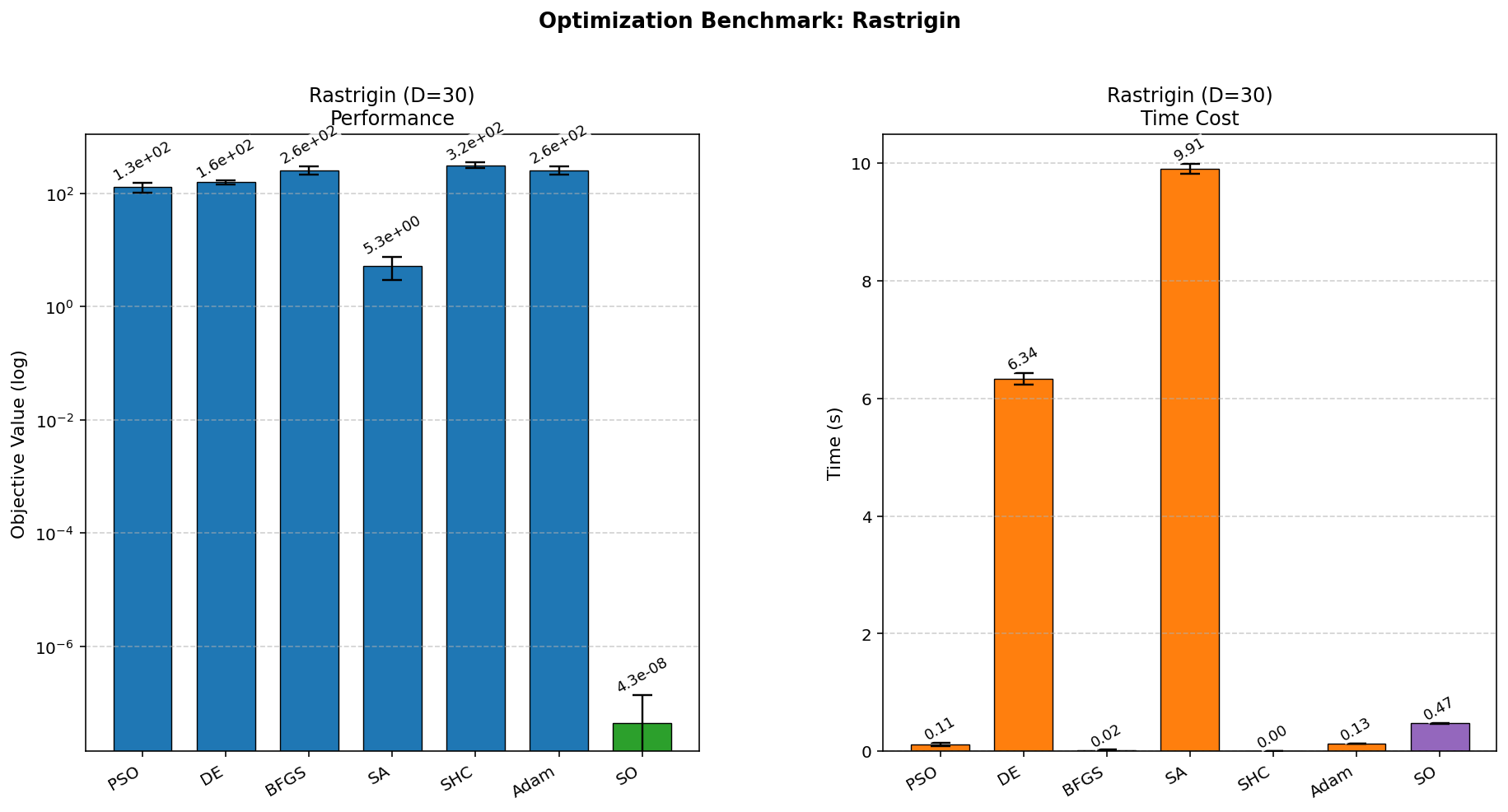}
\includegraphics[width=2.2in, height=1.2in]{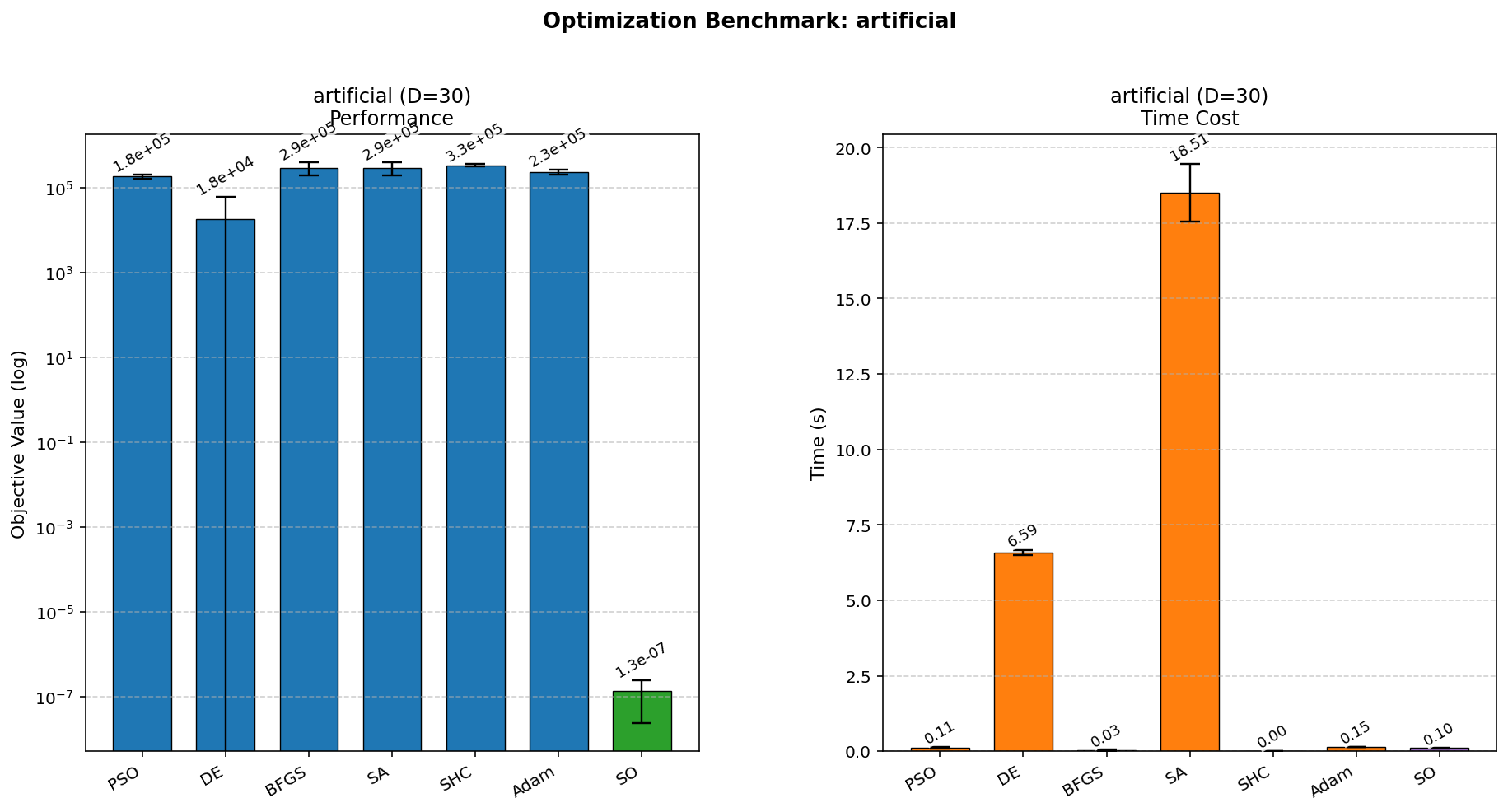}
\captionof{figure}{Rastrigin and Artificial functions}
\end{minipage}
\end{center}

Based on the experimental results, we draw the following key conclusions:
\begin{enumerate}
    \item Our algorithm consistently outperforms the baseline methods in terms of locating global minima when computational time is not the primary constraint.
    \item Even when time efficiency is taken into account, our method outperforms both Simulated Annealing (SA) and Differential Evolution (DE).
    \item The results highlight our algorithm’s ability to achieve a favorable balance between solution quality and computational efficiency.
\end{enumerate}

\subsection{Discussion}

The proposed Sampling-based Optimization (SO) algorithm consistently achieves superior or near-optimal objective values across a diverse set of challenging benchmark functions. It performs particularly well on complex, multimodal, and deceptive landscapes such as the Griewank, Rastrigin, and Schwefel functions. Although it is marginally slower than certain local optimization methods like BFGS on relatively simple functions (e.g., Sphere), the SO algorithm strikes a compelling balance between solution quality and robustness, positioning it as an effective and dependable optimization strategy.

The numerical experiments highlight the effectiveness and versatility of the proposed method across a broad spectrum of problem types. Several key observations emerge from our empirical study:

\begin{itemize}
    \item \textbf{Dimensional Scalability:} The SO algorithm maintains consistent performance across both low-dimensional and high-dimensional optimization problems, indicating strong scalability with respect to problem dimensionality.

    \item \textbf{Robustness to Landscape Complexity:} The algorithm demonstrates high resilience to multimodal, non-convex, and deceptive objective landscapes. It effectively mitigates the risk of premature convergence and consistently explores globally optimal regions.

    \item \textbf{Adaptability to Smoothness:} The method performs reliably on both smooth and non-smooth objective functions, including those that are discontinuous or non-differentiable, suggesting broad applicability across various classes of optimization problems.
\end{itemize}

\vspace{5mm}

{\bf Acknowledgement:} I would like to thank Zhenyao Sun and  Rongchan Zhu
for their valuable discussions.

\end{document}